\documentclass[10pt]{article}

\usepackage{amsmath, amssymb, amsthm, latexsym, bbm}

\usepackage{epsfig}
\usepackage{psfrag}
\usepackage{subfigure}

\allowdisplaybreaks
 \usepackage{hyperref}

\usepackage{color}
\def\loc{{\rm loc}}
\newtheorem{theorem}{Theorem}[section]
\newtheorem{corollary}[theorem]{Corollary}
\newtheorem{lemma}[theorem]{Lemma}
\newtheorem{proposition}[theorem]{Proposition}
\newtheorem{definition}[theorem]{Definition}

\newtheorem{remark}[theorem]{Remark}
\numberwithin{equation}{section}

\newcommand{\pf}{\noindent{\bf Proof.} }
\def\qed{{\hfill $\Box$ \bigskip}}
  
 \def\sE {{\cal E}} \def\sF {{\cal F}}
  
  \def\sL {{\cal L}}

\def\proof{{\medskip\noindent {\bf Proof. }}}

\def\RR{\mathbb{R}}
\def\br{\mathbb{R}}

\def\P{\mathbb{P}}

\def\r{\rho}
\def\b{\beta}
\def\aa{\alpha}

\def\ce{\mathcal{E}}
\def\cf{\mathcal{F}}
\def\EE{\mathcal{E}}
\def\FF{\mathcal{F}}

\def\R{{\mathbb R}}

\def\wt{\widetilde}
\def\eps{\varepsilon}

\def\bE {{\mathbb E}}
\def\bea{\begin{align*}}
\def\eea{\end{align*}}
\def\bee{\begin{equation}}
\def\eee{\end{equation}}
\def\bP {{\mathbb P}}

\def\nn{{\nonumber}}
\def\Lip{{\rm Lip}}
\def\FK{\mathrm{FK}}
\def\CS{\mathrm{CS}}

\def\RVD{\mathrm{RVD}}
\def\lam{{\lambda}}
\def\PI{\mathrm{PI}}
\def\<{\langle}
\def\>{\rangle}
\def\rf{{\rm ref}}
\def\lam{\lambda}

\def\1{\mathbbm {1}}

\def\bJ{{\bf J}}
\def\1{{\mathbbm 1}}

\makeatletter
\@addtoreset{equation}{section}

\makeatother

\begin{document}
\bibliographystyle{plain}

\title{\Large \bf
Heat kernels  for
 reflected diffusions with jumps \\ on
inner uniform domains}
\author{Zhen-Qing Chen \quad Panki Kim \quad Takashi Kumagai
\quad Jian Wang}

\maketitle

\begin{abstract}
In this paper, we study
sharp two-sided
heat kernel estimates for
a large class of symmetric reflected diffusions with jumps on the closure
of an inner uniform
domain $D$ in a length metric space.
The length metric  is
the intrinsic metric of a strongly local Dirichlet form.
When $D$ is an inner uniform domain
in the Euclidean space,  a prototype
for a special case of
the  processes under consideration
are
symmetric reflected diffusions
with jumps on $D$,  whose infinitesimal generators are non-local
 (pseudo-differential) operators $\sL $ on $D$ of the form
 $$
  \sL  u(x) =\frac12 \sum_{i, j=1}^d \frac{\partial}{\partial x_i} \left(a_{ij}(x)
  \frac{\partial u(x)}{\partial x_j}\right)
+ \lim_{\eps \downarrow 0} \int_{\{y\in D: \,
 \rho_D(y, x)>\eps\}}  (u(y)-u(x)) J(x, y)\, dy
 $$
 satisfying ``Neumann boundary condition".
    Here, $\rho_D(x,y)$ is the length metric on $D$, $A(x)=(a_{ij}(x))_{1\leq i,j\leq d}$
  is a measurable $d\times d$ matrix-valued function on $D$ that is
 uniformly
elliptic and bounded, and
$$
J(x,y):= \frac{1}{\Phi(\rho_D(x,y))}
\int_{[\alpha_1, \alpha_2]}  \frac{c(\alpha, x,y)} {\rho_D(x,y)^{d+\alpha}} \,\nu(d\alpha) ,
$$
where $\nu$ is a
finite measure on $[\alpha_1, \alpha_2] \subset (0, 2)$,
$\Phi$ is an increasing function on $[ 0, \infty )$ with $c_1e^{c_2r^{\beta}} \le \Phi(r) \le c_3 e^{c_4r^{\beta}}$
for some $\beta \in [0,\infty]$,
and  $c(\alpha , x, y)$ is  a jointly measurable function that
is bounded between two positive
constants and is symmetric in $(x, y)$.

\bigskip

\noindent\textbf{Keywords:}
reflected diffusions with jumps; symmetric Dirichlet form; inner uniform domain; heat kernel

\bigskip
\noindent {\bf AMS 2020 Mathematics Subject Classification}: Primary
 60J35, 60J76,
Secondary 31C25,
 35K08.
\end{abstract}

\section{Introduction and setup}\label{s:IS}
Let $E$ be a locally compact separable metric space,
and $m$ a $\sigma$-finite Radon measure with full support
on $E$. Let $(\EE^0, \FF^0)$ be a strongly local regular Dirichlet form
on $L^2(E; m)$, and  $\mu^0_{\<u\>}$ be the ($\EE^0$-) energy measure
of $u\in \FF^0$ so that
 $\EE^0 (u,u)= \frac12 \mu^0_{\<u\>}(E)$.
 We assume that the intrinsic metric
$\rho$  of $(\EE^0, \FF^0)$
defined by
$ \rho (x, y)=\sup\big\{f(x)-f(y): \ f\in \FF^0 \cap C_c(E)
\hbox{ with } \mu_{\<f\>}^0(dz) \leq m(dz)\big\}
$
is finite
for any $x,y\in E$
and
 induces the original topology on $E$, and
  that
  $(E, \rho)$ is a complete metric space.
In the above definition of $\rho (x, y)$, the condition
$\mu_{\<f\>}^0(dz) \leq m(dz)$ means  that  $\mu_{\<f\>}^0 (A) \leq m(A)$ for every $A\in {\mathcal B}(E)$.
We refer the reader to \cite[Theorem 2.11]{GS}
for  a summary of basic properties of $(E, \rho)$. In particular, $(E, \rho)$ is a geodesic length space; that is,
for any pair $x, y\in E$, there exists a continuous curve
$\gamma: \ [0, 1]\to E$
with $\gamma (0) = x$, $\gamma (1)=y$
and for every $s,  t \in [0, 1]$,
$\rho (\gamma(s), \gamma(t)) =|t-s|\, \rho (x, y)$.
Such a curve is called   a minimal geodesic (parameterized by
a multiple of arc length).
In the following,  we will always use the intrinsic metric $\rho$ for $E$, and use
$B(x, r)$ for  the ball centered at $x$ with
radius $r$ under $\rho$-metric.
Denote $m(B(x, r))$ by $V(x, r)$. We consider two
conditions as follows.

\begin{description}
\item{(VD)} ({\bf Volume doubling property}) There is a constant
  $c_1>0$ so that for every $x\in E$ and $r>0$, $V(x, 2r)\leq c_1\, V(x, r)$.

\item{(PI(2))}
({\bf Poincar\'e inequality}) There is a constant
 $c_2>0$
 so that for every $x\in E$, $r>0$ and $f\in \FF^0$,
 \begin{equation} \label{e:1.1}
 \min_{a\in\RR}\int_{B(x, r)}(f(y)-a)^2\, m(dy)
    \leq
    c_2 r^2 \, \mu_{\<f\>}^0(B(x, r)).
 \end{equation}
\end{description}
The Dirichlet form $(\EE^0, \FF^0)$ satisfying the above properties
is called a {\it Harnack-type Dirichlet space} in \cite[Chapter 2]{GS}.
The state space $E$ for Harnack-type Dirichlet space $(\EE^0, \FF^0)$  is connected and
the diffusion
process $Z^0$ associated with  $(\EE^0, \FF^0)$
is conservative (see \cite[Lemma 2.33]{GS}).

The following is Sturm's generalization to the strongly local Dirichlet form setting of
a celebrated
result by Grigor'yan and Saloff-Coste on Riemannian manifolds. (See \cite{St1, St2, St3}.)
   In this paper, we use $:=$ as a way of definition.

\begin{theorem}\label{T:1.1}
The following are equivalent
for the strongly local Dirichlet form
 $(\EE^0, \FF^0)$ on $L^2(E; m)$.
\begin{description}
\item{\rm (i)} {\rm(VD)} and {\rm(PI(2))} hold.
\item{\rm (ii)} The uniform parabolic Harnack inequality holds.
Namely, there exists a constant $C>0$ such that for every $x_0 \in E $, $t_0\ge 0$, $r>0$
 and for
every non-negative function $u=u(t,x)$ on $[0,\infty)\times E$ that is
caloric $($or space-time harmonic$)$ in
  cylinder $(t_0, t_0+4r^2)\times B(x_0,2r)$,
\[
\mathop{\rm ess \; sup \, }_{Q_- }u\le C \,\mathop{\rm ess \; inf \, }_{Q_+}u,
\]
where $Q_-:=(t_0+r^2,t_0+2r^2)\times B(x_0, r)$ and $Q_+:=(t_0+3r^2, t_0+4r^2)\times B(x_0, r)$.
\item{\rm (iii)} Aronson-type heat kernel estimates for the diffusion
process $Z^0$ associated with $(\EE^0, \FF^0)$ hold, i.e.,  $Z^0$ admits a jointly
continuous transition density $p^0(t, x, y)$ with respect to $m$ on
$(0, \infty)\times E\times E$ and there are
constants $c_1,c_2\ge 1$ so that
$$
\frac{1}{c_1V(x, \sqrt{t})} \exp\left(-\frac{c_2\rho (x, y)^2}{t}\right)
\leq p^0(t, x, y)\leq \frac{c_1}{V(x, \sqrt{t})} \exp\left(-\frac{\rho (x, y)^2}{c_2t}\right)
$$
for every $x, y\in E$ and $t>0$.
\end{description}
\end{theorem}

\smallskip

For a domain (i.e.\ connected open subset) $D$  of the length metric space
$(E, \rho)$, define for $x, y\in D$,
\begin{equation}\label{e:rhoD}
\rho_D (x, y)=\inf\{\hbox{length}(\gamma): \hbox{a continuous curve }
\gamma \hbox{ in } D \hbox{ with } \gamma (0)=x \hbox{ and }
\gamma (1)=y\}.
\end{equation}
Denote by ${\bar  D}$ the completion of $D$ under the metric $\rho_D$.
Note that
  $(\bar D,\rho_D)$ is a length metric space but may not be
  locally compact in general.
  For example, see \cite[Remark 2.16]{GS}.
   We extend the definition of $m|_D$ to ${\bar D}$ by setting
$m|_D ({\bar D}\setminus D)=0$.

For notational simplicity, by abusing the notation a little bit, we will use $m$ to denote this
measure $m|_D$.
Following \cite[Definition 3.6]{GS},
we say that $D$ is {\it inner uniform}
if there are constants
 $C_1, C_2 \in (0,\infty)$ such
that, for any $x, y \in D$, there exists a continuous curve
$\gamma_{x, y}: [0, 1] \to D$ with
$\gamma_{x, y}(0) = x$, $\gamma _{x, y}(1)=y$
 and satisfying the following two properties:
\begin{description}
\item{(i)} The length
of
$\gamma_{x,y}$ is at most $ C_1 \rho_D (x, y)$;
\item{(ii)} For any $z \in \gamma_{x,y}([0, 1])$,
$$
\rho (z, \partial D):=\inf_{w \in \partial D} \rho (z, w)   \geq C_2
\frac{\rho_D (z, x) \rho_D (z, y)}{\rho_D(x, y)}.
$$
\end{description}
We call the above constants $(C_1, C_2)$ the characteristics of the inner uniform domain $D$.
 Note that, when $D$ is inner uniform, $(\bar D,\rho_D)$
    is    locally compact (see \cite[Lemma 3.9]{GS}).
Let
$\FF^0_D=\{f\in \FF^0: \ f=0 \,\, \EE^0 \hbox{-q.e. on }
 D^c\}$.
It is well known that $(\sE^0, \FF^0_D)$ is the part Dirichlet form of $(\sE^0, \sF^0)$ on $D$,
or equivalently, it is the Dirichlet form on $L^2(D; m)$ of the subprocess of the diffusion process $X^0$
associated with $(\sE^0, \sF^0)$ killed upon leaving $D$.
A function $f$ is said to be locally in $\FF^0_D$, denoted as $f\in \FF^0_{D, \loc}$,
 if for every relatively compact subset $U$ of $D$, there is a function $g\in \FF^0_D$
 such that $f=g $ $m$-a.e. on $U$.
By \cite[Proposition 2.13]{GS}, we have  that for $x,y \in D$,
$$
 \rho_D (x, y)=\sup\left\{f(x)-f(y): \ f\in \FF^0_{D,{\rm loc}} \cap
 C_c(D)
\hbox{ with } \mu_{\<f\>}^0 (dz) \leq m(dz) \right\}.
$$

Define $\sF_D^{0, \rf} :=\{f\in \FF^0_{D,{\rm loc}}: \mu_{\<f\>}^0(D)<\infty\}$
and
\begin{equation}\label{EErf}
 \EE^{0, \rf} (f, f):=\frac12 \mu_{\<f\>}^0(D) \quad \hbox{for } f\in \FF_D^{0, \rf}.
\end{equation}
The bilinear form
 $(\EE^{0, \rf}, \sF_D^{0, \rf} \cap L^2(D; m) )$
is the active reflected Dirichlet form of $(\EE^0, \sF^0_D)$,
which is known to be a Dirichlet form on $L^2(\bar D; m)=L^2(D; m)$,
  see \cite[Chapter 6]{CF}.
Denote
$B_{{\bar D}}(x, r)$ by the ball centered at $x$ with
radius $r$ under $\rho_D$-metric, i.e., $
B_{{\bar D}}(x, r):=\{y\in {\bar D}: \rho_D (x, y)<r\}.
$
Denote $m(B_{{\bar D}}(x, r))$ by $V_D(x, r)$.
Let ${\rm Lip}_c( {\bar D})$ be the space of Lipschitz functions with compact support in ${\bar D}$.
The following is established in \cite[Lemma 3.9, Theorem 3.10,
Theorem 3.13 and  Corollary 3.31]{GS}.

\begin{theorem}\label{T:1.2}
 Suppose that $(\EE^0, \FF^0)$ is a strongly local
regular Dirichlet form on $L^2(E; m)$ that admits a carr\'e du champ operator
 $\Gamma_0$ $($that is, for every $u\in \FF$, $\mu_{\<u\>}^0(dx)=\Gamma_0 (u, u)\,m(dx)$ and $\Gamma_0 (u, u) \in L^1(E; m)$$)$. Assume that {\rm (VD)} and {\rm (PI(2))} hold for $(\EE^0, \FF^0)$
on $(E, \rho, m)$. Suppose that $D$ is  an inner uniform subdomain of $E$.
Then $(\EE^{0, \rf},
\FF_D^{0, \rf} \cap L^2(D; m))$
is a strongly local regular
Dirichlet form on $L^2(D; m)$ with core ${\rm Lip}_c( {\bar D})$ and that
the volume doubling property and the Poincar\'e inequality hold
for $(\EE^{0, \rf},
\FF_D^{0, \rf} \cap L^2(D; m))$
 on $({\bar D}, \rho_D, m)${\rm:}
\begin{description}
\item{{\rm (VD)}}
$(${\bf Volume doubling property on ${\bar D}$}$)$ There is a constant
$C_3>0$ so that for every $x\in {\bar D}$ and $r>0$, $V_D(x, 2r)\leq
C_3\, V_D(x, r)$.

\item{{\rm (PI(2))}}
$(${\bf Poincar\'e inequality on ${\bar D}$}$)$ There is a  constant
$C_4>0$ so that for every $x\in {\bar D}$, $r>0$ and
$f\in \FF_D^{0, \rf} \cap L^2(D; m)$,
\begin{equation} \label{e:1.3}
 \min_{a\in\RR}\int_{B_{{\bar D}}(x, r)}(f(y)-a)^2 \,m(dy)
    \leq
       C_4 r^2 \, \mu_{\<f\>}^0(B_{{\bar D}}(x, r)).
 \end{equation}
\end{description}

Consequently, the diffusion process associated with
$(\EE^{0, \rf}, \FF_D^{0, \rf})$ admits a jointly
continuous transition density function
$p_D^N(t, x, y)$ on
$(0, \infty)\times {\bar D}\times {\bar D}$, and there are
constants $c_3, c_4 \ge 1$ depending on
 $C_3,C_4$
 so that
\begin{equation}\label{e:1.4}
\frac{c_3^{-1}}{V_D(x, \sqrt{t})} \exp\left(-\frac{c_4\rho_D (x, y)^2}{t}\right)
\leq p_D^N(t, x, y)\leq \frac{c_3}{V_D(x, \sqrt{t})}
\exp\left(-\frac{\rho_D (x, y)^2}{c_4 t}\right)
\end{equation}
for every $x, y\in {\bar D}$ and $t>0$.
\end{theorem}

\medskip

Throughout this paper, we assume that $(\EE^0, \FF^0)$ and $D$ satisfy
the assumptions of
Theorem \ref{T:1.2}.
Recall that $\rho_D$ is defined in \eqref{e:rhoD},
 ${\bar D}$  is the completion of the metric space $(D, \rho_D)$, and (by abusing the notation) $m$ is the measure $m|_D$
 extended to $\bar D$ by setting
 $m({\bar D}\setminus D)=0$.
By Theorem \ref{T:1.2}, the volume doubling property  (VD)  holds for the inner uniform subdomain $D$ of $E$.
Note that,
 (VD) condition on $V_D(x,r)$ is equivalent to the following: there
exist constants
$c_2,d_2>0$ depending on $C_3$ only
such that for all $x\in {\bar D}$,
\begin{equation}\label{univd}
\frac{V_D(x,R)}{V_D(x,r)}\le c_2 \Big(\frac
Rr\Big)^{d_2}\quad \hbox{for } R\geq r>0.
\end{equation}
Since $D$ is connected,
the  reverse volume doubling property (RVD) holds for ${\bar D}$; that is, there exist constants $c_1,d_1>0$, depending on $C_1, C_3, C_4$  only, such that for
all $x\in {\bar D}$,
\begin{equation}\label{univd*}
 \qquad \frac{V_D(x,R)}{V_D(x,r)}\ge c_1 \Big(\frac
Rr\Big)^{d_1} \quad  \hbox{for } 0< r\le R \le 2\, \text{diam}(D).
\end{equation}
See \cite[Proposition 2.1 and  a paragraph before Remark 2.1]{YZ}.
We further note that, by \eqref{univd}, for all
  $x, y\in {\bar D}$  and $0<r\le R$,
\begin{equation}\label{e:1.8}
\frac{V_D(x,R)}{V_D(y,r)}\le  \frac{V_D(y,\rho(x,y)+R)}{V_D(y,r)}\le c_2\left(\frac{\rho(x,y)+R}{r}\right)^{d_2}.
\end{equation}
This in particular implies that, with $\wt c_2:= 2^{d_2} c_2$,
$$
 {\wt  c_2}^{-1} \left(\frac{\rho(x,y) }{r} \vee 1 \right)^{-d_2}
 \leq \frac{V_D(x,r)}{V_D(y,r)}   \le \wt  c_2\left(\frac{\rho(x,y) }{r} \vee 1 \right)^{d_2}
\quad \hbox{for all } x, y\in
\bar D
 \hbox{ and } r>0.
$$
  Here and in what follows, we use notations
 $a\wedge b:=\min \{a, b\}$ and $a\vee b:=\max\{a, b\}$ for  $a,b\in \R$.

\noindent {\bf Characteristic constants.}  Many estimates in this paper will depend on the characteristics $(C_1, C_2)$ of the inner uniform domain $D$ and on the constants $C_3$ and $C_4$
in (VD) and (PI(2)) of  Theorem \ref{T:1.2}. For convenience, we will call $(C_1, C_2, C_3, C_4)$ the {\it characteristic constants} of the domain
$D$ in this paper.

\medskip

 In this paper,  we are concerned with  Markov processes on $\bar D$ associated with the following type of non-local
symmetric Dirichlet forms $(\EE, \FF)$ on $L^2(D; m)$:
 \begin{equation}\label{e:sF}
 \FF=\FF_D^{0, \rf} \cap L^2(D; m),
 \end{equation}
 and, for $u\in \FF$,
 \begin{equation}\label{e:sE}
 \EE (u, u)  = \EE^{0,\rf}(u, u)+
 \frac{1}{2} \int_{D\times D}
 (u(x)-u(y))^2 J(x, y)\, m(dx)\,m(dy),
 \end{equation}
where  $J(x,y)$ is a non-negative
symmetric measurable function  on $D\times D \setminus {\rm diag}$ satisfying certain conditions to be specified below.
 Here and in what follows, ${\rm diag}$ is the diagonal of $D\times D$; that is, ${\rm diag}:=\{(x, x): x\in D\}$.

\medskip

Let
$\phi$ be a  strictly increasing function on $[0, \infty)$
such that
$\phi (0)=0$, $\phi (1)=1$  and
there exist constants $C_{5} \ge1$ and
 $0<\aa_* \leq \aa^*<2$
so that
 \begin{equation}\label{eqn:poly}
C_{5}^{-1}  \Big(\frac Rr\Big)^{\aa_*} \le
\frac{\phi (R)}{\phi (r)} \le C_{5}\Big(\frac
Rr\Big)^{\aa^*}
\quad \hbox{for every } 0<r<R<\infty.
\end{equation}
 Since the constant $\aa^*$ in  \eqref{eqn:poly} is strictly less than $2$,
there exists $c_3>0$ depending on
$C_5$ and $\aa^*$ such that
\begin{equation}\label{eqn:cond1}
\int_0^r\frac {s}{\phi (s)}\, ds \le \frac{ c_3 r^2}{\phi
(r)} \qquad   \hbox{for every } r>0.
\end{equation}

\begin{definition}\label{D:1.3}
\rm  Let  $\beta \in [0, \infty]$ and $\phi$
be a
strictly increasing function  on
 $[0, \infty)$  with $\phi (0)=0$ and $\phi (1)=1$
that satisfies the condition \eqref{eqn:poly} (with $0<\aa_* \leq \aa^*<2$).
For  a non-negative symmetric measurable function  $J(x, y)$ on $D\times D \setminus {\rm diag}$, we say
\begin{enumerate}

\item[(i)] condition  $(\bJ_{\phi, \beta, \leq})$   holds if there are positive constants $\kappa_1$ and $\kappa_2$ so that
 \begin{equation} \tag{$\bJ_{\phi, \beta, \leq}$}
J(x, y) \leq \frac{\kappa_1 }{V_D(x,  \rho_D(x, y))\phi  (\rho_D(x, y) )\exp (\kappa_2 \rho_D(x, y)^\beta)}
\quad \hbox{for   } (x, y) \in D\times D \setminus {\rm diag};
\end{equation}

\item[(ii)]  condition  $(\bJ_{\phi, \beta, \geq})$   holds if there are positive constants $\kappa_3$ and $\kappa_4$ so that
 \begin{equation} \tag{$\bJ_{\phi, \beta, \geq}$}
J(x, y) \geq  \frac{\kappa_3 }{V_D(x,  \rho_D(x, y))\phi  (\rho_D(x, y) )\exp (\kappa_4 \rho_D(x, y)^\beta)}
\quad \hbox{for   } (x, y) \in D\times D \setminus {\rm diag};
\end{equation}

\item[(iii)]  condition  $(\bJ_{\phi, \beta})$   holds if  both $(\bJ_{\phi, \beta, \leq})$ and    $(\bJ_{\phi, \beta, \geq})$  hold
 with  possibly  different constants $\kappa_i$ in the upper and lower bounds;

\item[(iv)]  condition  $(\bJ_{\phi,  0_+, \leq})$   holds if
there are positive constants $\kappa_5$ and $\kappa_6$ so that
  \begin{equation} \tag{$\bJ_{\phi, 0_+, \leq}$}
  \left\{
\begin{split}
&J(x, y) \leq  \frac{\kappa_5 }{V_D(x,  \rho_D(x, y)) \phi_*  (\rho_D(x, y) ) }
 \quad \hbox{for   } (x, y) \in D\times D \setminus {\rm diag}, \\
&  \sup_{x\in D}\int_{\{y\in D: \rho_D(x,y)>1\}} \rho_D(x,y)^2 J(x, y)\,m(dy)
 \leq \kappa_6<\infty,
\end{split}
\right.
\end{equation}
 where
 \begin{equation}\label{e:phi*}
 \phi_*(r):= \phi(r) \1_{\{r \leq 1\}} + r^2\1_{\{r>1\}} \quad \hbox{for } r\geq 0 ;
 \end{equation}

 \item[(v)]  condition  $(\bJ_{\phi,  0_+, \leq 1})$   holds if
there are positive constants $ \kappa_7$ and $ \kappa_8$ so that
  \begin{equation} \tag{$\bJ_{\phi, 0_+, \leq 1}$}
  \left\{
\begin{split}
&J(x, y) \leq  \frac{ \kappa_7 }{V_D(x,  \rho_D(x, y)) \phi  (\rho_D(x, y) ) }
 \quad \hbox{for   } x\not= y \hbox{ in } D \hbox{ with } \rho_D (x, y)\leq 1,  \\
&  \sup_{x\in D}\int_{\{y\in D: \rho_D(x,y)>1\}} \rho_D(x,y)^2 J(x, y)\,m(dy)
\leq  \kappa_8<\infty;
\end{split}
\right.
\end{equation}

\item[(vi)] condition  $(\bJ_{\phi, \leq 1})$  holds  if
there are positive constants $\kappa_9$ and $\kappa_{10}$ so that
 \begin{equation} \tag{$\bJ_{\phi, \leq 1}$}
\left\{
\begin{split}
 & J(x, y) \leq \frac{\kappa_9 }{V_D(x,  \rho_D(x, y))\phi  (\rho_D(x, y) ) }
\quad \hbox{for   } x\not= y \hbox{ in }   D \hbox{ with }
\rho_D (x, y) \leq 1 ,  \\
&\sup_{x\in D} \int_{\{y\in D: \rho_D(x, y) >1\}} J(x, y) \,m(dy)
\leq \kappa_{10}<\infty.
\end{split}
\right.
\end{equation}
 \end{enumerate}
\end{definition}

 Clearly,  $(\bJ_{\phi,  0_+, \leq})\Longrightarrow  (\bJ_{\phi, 0, \leq})
\Longrightarrow   (\bJ_{\phi,  \leq 1})$ and
\begin{equation}\label{e:Jcom}
(\bJ_{\phi, \beta, \leq}) \Longrightarrow (\bJ_{\phi,  0_+, \leq})\Longrightarrow   (\bJ_{\phi,  0_+, \leq 1})
\Longrightarrow   (\bJ_{\phi,  \leq 1})
\quad  \hbox{for any  }\beta \in (0, \infty] .
\end{equation}
When $\beta =0$,
we will simply write $(\bJ_{\phi,  \leq})$, $(\bJ_{\phi,  \geq})$ and  $(\bJ_{\phi })$
for $(\bJ_{\phi, 0, \leq})$, $(\bJ_{\phi, 0,  \geq})$ and $(\bJ_{\phi, 0 })$, respectively.
When $\beta =\infty$, conditions $(\bJ_{\phi, \infty, \leq})$ and $(\bJ_{\phi, \infty, \geq})$ are equivalent, respectively,  to
\begin{equation}  \label{e:1.12}
J(x, y) \leq \frac{\wt \kappa_1 }{V_D(x,  \rho_D(x, y))\phi  (\rho_D(x, y) ) }  \1_{\{\rho (x, y) \leq 1\}}
\quad \hbox{for   } (x, y) \in D\times D \setminus {\rm diag}
\end{equation}
 and
  \begin{equation}  \label{e:1.13}
J(x, y) \geq  \frac{\wt \kappa_3 }{V_D(x,  \rho_D(x, y))\phi  (\rho_D(x, y) ) }  \1_{\{\rho (x, y) \leq 1\}}
\quad \hbox{for   } (x, y) \in D\times D \setminus {\rm diag}.
\end{equation}

\medskip

We will see from  Proposition \ref{l:regular}  below
that, under  condition
 $(\bJ_{\phi, \leq 1})$,
$(\EE, \FF)$ is a regular Dirichlet form
on $L^2(D; m)$, and
there is a symmetric Hunt process $Y$ on ${\bar D}$ associated with it that can start from
$\EE$-quasi-everywhere in ${\bar D}$. Moreover, the process $Y$ is conservative;
that is, $Y$ has infinite lifetime
almost surely.  For notational convenience, we regard $0_+$ as   an ``added" or ``extended" number
and declare that it is larger than $0$ but  smaller than any positive real number.
With this notation, we can write, for instance, $(\bJ_{\phi, \beta, \leq})$
for $\beta\in [0, \infty] \cup \{0_+\}$.

 The  jumping intensity kernel $J(x, y)$ determines a
L\'evy system of $Y$, which describes the jumps of the process $Y$: for any non-negative
measurable function $f$ on $\br_+ \times {\bar D}\times {\bar D}$
with $f(s, x, x)=0$ for all $s>0$ and $x\in {\bar D}$,
and stopping time $T$ (with respect to the filtration of
$Y$),
\begin{equation}\label{levy-1}
\bE_x \left[\sum_{s\le T} f(s,Y_{s-}, Y_s) \right]=
\bE_x \left[ \int_0^T \left( \int_{D} f(s,Y_s, y) J(Y_s,y) m(dy)
\right) ds \right].\end{equation}
See, for example, \cite[proof of Lemma 4.7]{CK} and \cite[Appendix A]{CK2}.

\medskip

The goal of this paper is
to derive global two-sided sharp estimates on the heat kernel of Dirichlet form $(\EE, \FF)$ (i.e., the transition density function of the process $Y$)
 under the assumption that $J(x, y)$ satisfies $(\bJ_{\phi_1, \beta_*, \leq})$ and $(\bJ_{\phi_2, \beta^*, \geq})$ for some strictly increasing
 functions $\phi_1$ and $\phi_2$
 satisfying
 $\phi_i(0)=0$, $\phi_i(1)=1$ and \eqref{eqn:poly} (with $\phi_i$ in place of $\phi)$ for $1\le i\le 2$,
  and for $\beta_* $ and $\beta^*$ in $[0, \infty] \cup \{0_+\}$  in the following
 two cases:
  \begin{itemize}
\item $\beta_*=\beta^*=0$ and $\phi_1=\phi_2$; and
\item $0_+\leq \beta_*\le \beta^*\le \infty$ excluding $\beta_*=\beta^*=0_+$.
\end{itemize}
See Subsection \ref{Case3} in the Appendix of this paper for brief discussions on the case under conditions
$(\bJ_{\phi_1, \leq })$ and $(\bJ_{\phi_2, \beta^*, \geq})$ with $\beta^*\in (0, \infty]$.

\smallskip

We start with the $\beta_*=\beta^* = 0 $ case.
 Define for $x\in \bar D$, $t>0$ and $r\geq 0$,
\begin{equation}\label{e:pc}
  p^{(c)}(t,x, r) :=  \frac1{V_D(x,\sqrt{t})}\, \exp ( - r^2/t),
 \end{equation}
and for a strictly increasing function $\phi$ on $[0, \infty)$,
\begin{equation}\label{eqn:5}
 p_{\phi}^{(j)}(t,x, r) :=  \frac1{V_D(x,\phi^{-1}(t))} \wedge
\frac{t}{V_D(x, r) \phi (r)},
\end{equation}
where $  \phi^{-1}(r)$ denotes the inverse function of $\phi$.
Set
 \begin{equation}\label{e:1.20}
   H_{\phi, 0} (t, x, r) :=
 \frac{1 }{   V_D(x,\phi^{-1}(t) \vee \sqrt{t}\,)}
\wedge  \left( p^{(c)}(t,x, r)+p_{\phi}^{(j)}(t, x, r) \right).
 \end{equation}

\begin{theorem} \label{T:1.4}
 {\rm  (The $\beta_*=\beta^* = 0 $ case.)}
Assume that
condition  $(\bJ_{\phi })$  holds for a strictly increasing function $\phi $ satisfying \eqref{eqn:poly}.
 Then $(\sE, \sF)$
 defined by \eqref{e:sF}--\eqref{e:sE} is a regular Dirichlet form on $L^2(D; m)$
 and there is a  conservative  Feller  process $Y$ associated with
   it  that starts from every point in $\bar D$.
 Moreover, $Y$ has a jointly H\"older continuous transition density function $q(t, x, y)$
   that   enjoys  the following
    two-sided estimates. There
   are positive constants $c_i$, $1\leq i \leq 4$,
depending only on the characteristic constants
$(C_1, C_2, C_3, C_4)$ of $D$ and
 the   constant parameters in  \eqref{eqn:poly}  and   $(\bJ_{\phi})$,
such that
for every $t>0$ and $x, y\in {\bar D}$,
\begin{equation}\label{e:1.12a}
  c_1H_{\phi,0}   ( t, x, c_2 \rho_D(x, y)  )  \leq q(t, x, y) \leq c_3 H_{\phi, 0} (t, x, c_4 \rho_D(x, y)).
 \end{equation}
\end{theorem}

Note that, in the present setting,
the fact that
the underlying state space $D$ has infinite diameter is equivalent to the fact that $D$ has infinite
 volume; see \cite[Corollary 5.3]{GH}.
 So, when ${\rm diam}(D)=\infty$, we can deduce Theorem \ref{T:1.4} from
  \cite[Theorem 1.13]{CKW4}. Furthermore, as seen from \cite{GHH} and \cite[Remark 1.19]{CKW2},
all results of the paper \cite{CKW4} continue to hold  for bounded state space with obvious localized versions.
Thus when ${\rm diam}(D)<\infty$, according to \cite[Theorem 1.13]{CKW4} again,
for any $T_0>0$, one can obtain estimates
of $q(t,x,y)$ for $t\in (0, T_0]$ with constants $c_i$ further dependent on $T_0$.
So for  Theorem \ref{T:1.4}, it remains to prove the statement for $t\ge T_0$ and ${\rm diam}(D)<\infty$.
In fact, when $D$ is bounded, it holds that
$V_D(x, \sqrt{t}) \simeq 1$ for all $x\in D$ and $t\geq 1$,
and the large
time estimates \eqref{e:1.12a} (i.e., $t\in [1,\infty)$) of the heat kernel $q(t,x,y)$  are
simply
\begin{equation}\label{e:1.27}
q(t,x,y)\simeq1  \quad \hbox{for } x,y\in \bar D \hbox{ and }  t\ge1.
\end{equation}
Here and in what follows,
$f(t,x, r)\simeq g(t,x, r)$
means that there exist constants $c_{1},c_{2}>0$ such that
$
c_{1}g(t,x,r)\leq f(t,x,r)\leq c_{2}g(t,x,r)
$
for the specified range of the argument $(t,x,r)$.
The  full proof  of Theorem \ref{T:1.4} will be  given in  Section \ref{S:2}.

Note also that when $D$ is bounded, conditions $( \bJ_{\phi, \beta, \leq})$ and
 $ (\bJ_{\phi, \beta, \geq})$ with $\beta \in \{0_+ \} \cup (0, \infty]$ are reduced to $ (\bJ_{\phi,   \leq})$ and $( \bJ_{\phi,    \geq})$, respectively.
 In view of the above theorem,
the main task of the paper is to consider  the case where $D$ is unbounded and
the jumping density kernel $J(x,y)$ satisfying conditions   $(\bJ_{\phi_1, \beta_*,  \leq})$
and  $(\bJ_{\phi_2, \beta^*, \geq})$    with $0_+ \leq \beta_*\le \b^*\le \infty$, which is
harder  than the case $\beta_*= \beta^*=0$ case.

 \medskip

To state our main result in this direction, we need some notations.
Recall that   $ p^{(c)}(t,x, r)$ and
$ p_{\phi}^{(j)}(t,x, r)$ are  defined by  \eqref{e:pc} and \eqref{eqn:5}, respectively.
For $\beta \in [0, \infty]$ and
  a strictly increasing function $\phi$ on $[ 0, \infty )$  with
   $\phi (0)=0$
   and $\phi (1)=1$,
  set for $x\in \bar D$, $t>0$ and $r\geq 0$,
  \begin{equation}\label{eqn:4}
 p^{(j)}_{\phi, \beta } (t,x, r) :=
  \frac1{V_D(x,\phi^{-1}(t))} \wedge
 \frac{t}{V(x, r) \phi (r)  \exp (r^\beta) }.
\end{equation}
In particular,
 $p^{(j)}_{\phi, 0} (t,x, r) \simeq p^{(j)}_{\phi} (t,x, r)$.
Define $\beta \in (0, 1]$,
  $$
 H_{\phi,\beta} (t, x, r) :=
 \begin{cases} \displaystyle
 \frac{1}{V_D(x,\sqrt{t})} \wedge \left(  p^{(c)}(t, x, r)+p^{(j)}_{\phi, \beta} (t, x, r)\right)
 \quad & \hbox{if } t \in (0,1] ,  \smallskip \\
 \displaystyle   \frac{1}{V_D(x,  \sqrt{t} ) } \exp  \left( - \left( r^\beta \wedge  ({r^2}/{t})\right) \right)
    &\hbox{if } t \in (1,\infty);
    \end{cases}
 $$
 for $\beta\in (1,\infty)$,
$$
 H_{\phi,\beta} (t,x,r) :=
 \begin{cases} \displaystyle
 \frac{1}{V_D(x,\sqrt{t})} \wedge \left( p^{(c)}(t,x, r)+p_{\phi, \beta }^{(j)}(t, x, r) \right)
  \quad & \hbox{if } t \in (0,1] \hbox{ and }  r \leq 1,  \smallskip \\
   \displaystyle  \frac{t}
    {V_D(x, r )\phi(r)}
     \exp \left(  - \left(r(1+\log ^+ (r/t))^{(\beta-1)/\beta}\right)\wedge r^\beta \right)
  & \hbox{if } t \in (0,1] \hbox{ and }    r> 1 , \smallskip \\
  \displaystyle   \frac{1}{V_D(x,\sqrt{t})}\exp \left( - \left(  r\, ( 1+\log^+ (r/t) )^{(\beta-1)/\beta} \right) \wedge (r^2/t)\right)
  &\hbox{if } t \in (1,\infty)  ; \end{cases}
$$
where $\log^+(x):=\log (x\vee 1)$,
and $H_{\phi,\infty} (t,x,r) := \lim_{\beta \to \infty} H_{\phi,\beta} (t, x, r) $ for $\beta =\infty$, that is,
$$
 H_{\phi,\infty} (t,x,r) :=   \begin{cases}  \displaystyle
 \frac{1}{V_D(x,\sqrt{t})} \wedge \left( p^{(c)}(t,x, r)+p_{\phi, \beta }^{(j)}(t, x, r)\right)
 \qquad & \hbox{if } t \in (0,1] \hbox{ and }  r \le1,
  \smallskip \\
 \displaystyle    \frac{t}
    {V_D(x, r )\phi(r)}
     \exp \left(  - r(1+\log^+ (r/t) ) \right)
  & \hbox{if } t \in (0,1]  \hbox{ and }
r  >1 ,  \smallskip \\
 \displaystyle    \frac{1}{V_D(x,\sqrt{t})}\exp \left( -  \left( r\, \big( 1+\log^+ (r/t)  \big)\right)  \wedge (r^2/t)\right)
  &\hbox{if } t \in (1,\infty)  .
  \end{cases}
$$
We further define for $x\in \bar D$, $t>0$ and $r\geq 0$,
$$
 H_{\phi,0_+} (t,x,r) :=  \frac{1}{V_D(x,
\sqrt{t}
 ) } \wedge \left( p^{(c)}(t,x, r)+ p_{\phi_*}^{(j)}(t, x,  r) \right)  ,
$$
where $\phi_*$ is given by \eqref{e:phi*}.

\medskip

 We can write the above functions in a more explicit way.
 For $\beta\in (0,1]$,
$$
 H_{\phi,\b} (t,x,r)\asymp
 \begin{cases} \displaystyle
 \frac{1}{V_D(x,\sqrt{t})} \wedge \left(
 p^{(c)}(t, x, r)+p_{\phi, \beta }^{(j)}(t, x, r)\right)
 \quad & \hbox{if } t \in (0,1] \hbox{ and }
  r\leq 1,  \smallskip \\
 \displaystyle   \frac{t}{V_D(x,r)
     \phi (r)} \exp(- r^{\beta} )
 \quad & \hbox{if } t \in (0,1] \hbox{ and }
 r \ge 1,  \smallskip \\
  \displaystyle     \frac{1}{V_D(x,\sqrt{t})}\, \exp(-   r^2 /t )
\quad &\hbox{if } t \in (1,\infty) \hbox{ and }t \ge r^{2-\beta},  \smallskip \\
 \displaystyle\frac{1}{V_D(x,\sqrt{t})}\, \exp(-   r^{\beta} )
\quad &\hbox{if } t \in (1,\infty) \hbox{ and }t \le r^{2-\beta};\end{cases}
$$ and for $\beta\in (1,\infty]$,
$$
 H_{\phi,\b} (t,x,r)\asymp
 \begin{cases}
 \displaystyle \frac{1}{V_D(x,\sqrt{t})} \wedge \left(
 p^{(c)}(t,x, r)+p_{
 \phi,\beta }^{(j)}(t, x, r)\right)
 & \hbox{if } t \in (0,1] \hbox{ and }
  r \leq 1,  \smallskip \\
   \displaystyle  \frac{t}
    {V_D(x,r)
    \phi (r)}
     \exp \left(-r\left(1+\log^+ (r/t)\right)^{(\beta-1)/\beta} \right)
  & \hbox{if } t \in (0,1],
  r \ge 1\hbox{ and } t\ge r\exp(-r^\beta),  \smallskip \\
 \displaystyle  \frac{t}{V_D(x,r)
    \phi (r)}  \exp(-r^\beta)
  & \hbox{if } t \in (0,1],
  r \ge 1\hbox{ and } t\le r\exp(-r^\beta),  \smallskip \\
   \displaystyle    \frac{1}{V_D(x,\sqrt{t})}\, \exp(-   r^2 /t )
  &\hbox{if } t \in (1,\infty) \hbox{ and }t \ge r,  \smallskip \\
 \displaystyle\frac{1}{V_D(x,\sqrt{t})}\, \exp\left(- r\, \left( 1+\log^+ (r/t)  \right)^{(\beta-1)/\beta} \right)
  &\hbox{if } t \in (1,\infty) \hbox{ and }t \le r.\end{cases}
$$Here and in what follows, we write
 $f(t,x,r)\asymp g(t,x,r)$,
if there exist constants $c_k>0$, $k=1, \cdots, 4$, such that
$
c_1g( t,x, c_2r)\leq f(t,x,r)\leq c_3 g(t,x, c_4 r)
$
for the specified range of $(t,x,r)$.

The next theorem
 gives the existence and its upper bound for the transition density function of $Y$ associated with the Dirichlet
form $(\sE, \sF)$ of \eqref{e:sF}--\eqref{e:sE} under condition $({\bf J}_{\phi, 0_+,\le})$.
Its proof will be given in Section \ref{S:3}.
 Note that it holds on any inner uniform domain $D$, regardless whether it is bounded or not.

    \begin{theorem} \label{T:2.4}
Suppose that condition $({\bf J}_{\phi, 0_+, \leq})$ holds for  some  strictly increasing function $\phi$
on $[0, \infty)$ satisfying the condition \eqref{eqn:poly}.
Then the  Dirichlet form $(\sE,\sF)$ of \eqref{e:sF}--\eqref{e:sE} is regular on $L^2(D; m)$ and there is a
 conservative Feller process $ Y$ on $\bar D$ associated with it
 that starts from every point in $\bar D$. Moreover, $  Y$ has a jointly H\"older continuous transition density function
 $  q(t, x, y)$ on $(0, \infty) \times \bar D \times \bar D$ with respect to the measure $m$, and
 there are constants $c_1,c_2>0$
 such that
  \begin{equation}\label{e:hk00}
   q(t,x,y)\le  c_1  H_{\phi,0_+} (t,x, c_2 \rho_D(x, y) )
   \quad \hbox{for all } x,y\in {\bar D} \hbox{ and } t>0.
     \end{equation}
 Furthermore,
 for every  $c_3\in (0,1/2)$, there are positive constants  $c_4$ and $c_5 $ so
   that for any $x\in {\bar D} $  and
   $r \in (0, c_3${\rm diam}$(D))$,
\begin{equation}\label{e:2.14}
 c_4 r^2\le \bE_x \left[  \tau_{B_{{\bar D}}(x,r)} \right] \le c_5r^2,
\end{equation}
where
$ \tau_{B_{{\bar D}}(x,r)}:=\inf\{  t> 0:   Y_t \notin B_{{\bar D}}(x,r)\}$
is the first exit time from the ball $B_{{\bar D}}(x,r)$ by  the process $  Y$.
The positive constants
$c_1,  c_2$ and $c_4, c_5$
depend  only on the characteristic constants
 $(C_1, C_2, C_3, C_4)$ of $D$ and
 on the constant parameters in $({\bf J}_{\phi, 0_+, \leq })$
 and \eqref{eqn:poly} for the function $\phi$,
with constants $c_4$ and $c_5$, in addition,  on $c_3$ as well.
 \end{theorem}

The following is the main result of this paper on the two-sided heat kernel estimates of $Y$.

\begin{theorem}\label{T1}
{\rm  (The
$\beta_*\le \beta^* \le \infty $ in $\{0_+\}\cup (0, \infty]$ case.)}
Suppose that $D$ is unbounded.
Assume that $J(x,y)$ satisfies conditions  $({\bf J}_{\phi_1,  \beta_*, \leq })$ and $({\bf J}_{\phi_2, \beta^*, \geq })$
for some strictly increasing functions
 $\phi_1,\phi_2$
 satisfying $\phi_i(0)=0$, $\phi_i(1)=1$ and \eqref{eqn:poly} $($with $\phi_i$ in place of $\phi$$)$ for $i=1,2$,
 and for $\beta_*\leq \beta^*$ in $ \{0_+\} \cup (0, \infty] $   excluding   $\beta_*=\beta^*=0_+$.
Then the transition density function $q(t,x,y)$ of the conservative Feller process $Y$ associated with
$(\ce, \cf)$ has  the
 following estimates:  for every $t>0$  and $x, y \in {\bar D}$,
 \begin{equation}\label{e:1.23}
 c_1
 H_{\phi_2, \beta^*}   ( t, x, c_2 \rho_D(x, y)  )  \leq q(t, x, y) \leq c_3 H_{\phi_1, \beta_*} (t, x, c_4 \rho_D(x, y)),
  \end{equation}
where $c_i$, $1\leq i \leq 4$, are positive constants that
 depend only on the characteristic constants
 $(C_1, C_2, C_3, C_4)$ of $D$ and
 the constant parameters in  $({\bf J}_{\phi_1,  \beta_*, \leq })$ and  $({\bf J}_{\phi_2, \beta^*, \geq })$ as well as
 in \eqref{eqn:poly} for $\phi_1$ and $\phi_2$, respectively.
  \end{theorem}

\begin{remark} \rm
\begin{itemize}
\item[(i)]
A prototype of the model
considered in this paper is the
following. Suppose that $D$ is a Lipschitz domain in $\RR^d$,
and that $m(dx)=dx$ is the Lebesgue measure on $\RR^d$.
Let
\begin{equation}\label{e:ddf}\begin{split}
\EE (f, f)&=
\frac{1}{2} \int_{D} \nabla u(x)\cdot A(x) \nabla
v(x)\,dx+\frac{1}{2}\int_{D \times D} (f(x)-f(y))^2 \, J(x, y)\, dx\, dy, \\
\FF &= W^{1,2}(D):=\left\{f\in L^2(D;dx): \nabla f\in L^2(D;dx)\right\}.
\end{split}\end{equation}
Here,
$A(x):=(a_{ij}(x))_{1\leq i,
j\leq d}$ is a measurable $d\times d$ matrix-valued function on $D$ that is
uniformly elliptic and
bounded, and $J(x, y)$ is a symmetric
function
on $\br^d\times \br^d \setminus  {\rm diag}$
defined by
\begin{equation}\label{e:prot}
J(x,y)=
\int_{[\alpha_1, \alpha_2]}  \frac{c(\alpha, x,y)} {|x-y|^{d+\alpha} \, \Phi(|x-y|)}\,
\nu(d\alpha),
\end{equation}
where $\nu$ is a probability measure on $[\alpha_1, \alpha_2] \subset (0, 2)$, $\Phi$ is
  an increasing function on $[ 0, \infty )$ with
   \begin{equation}\label{e:1.26}
  c_1e^{c_2r^{\beta}} \le\Phi(r) \le c_3 e^{c_4r^{\beta}}  \quad \hbox{for some } \beta \in [0,\infty],
  \end{equation}
  and  $c(\alpha , x, y)$ is  a jointly measurable function that
is symmetric in $(x, y)$ and is bounded between two positive
constants.
It is easy to check that $(\EE, \FF)$ is a regular Dirichlet form on $L^2(D; dx)$.
When $\beta =0$ in \eqref{e:1.26} for $\Phi$
and $D=\R^d$, two-sided heat kernel estimates for
the Dirichlet form
\eqref{e:ddf} were obtained in \cite{CK3}.

\item[{\rm (ii)}]
When a Dirichlet
form $(\EE, \FF)$ is of the form \eqref{e:ddf} with
$D=\br^d$ and $A(x)\equiv 0$,
and $J(x, y)$ a symmetric
function defined on $\br^d\times \br^d \setminus
 {\rm diag}$
 that satisfies
\[\frac{1}{c|x-y|^d\phi   (|x-y|^d)}
\le J(x, y) \leq \frac{c}{|x-y|^d\phi   (|x-y|^d)}
\quad \hbox{for   } (x, y) \in \br^d\times \br^d \setminus
 {\rm diag},
\]
where $c\ge1$, and $\phi $ satisfies  $\phi (0)=0$, $\phi (1)=1$ and  condition
 \eqref{eqn:poly},
it is established in \cite{CK2} that the transition density function $q(t, x, y)$ of the associated pure jump process
is jointly H\"older continuous and
has the following two-sided sharp estimates: there is a constant $C>1$ so that
$$
C^{-1} p_0^{(J)} (t, |x-y|) \leq q(t, x, y) \leq C p_0^{(J)} (t, |x-y|)
\quad \hbox{for } (t, x, y) \in (0, \infty) \times \R^d\times \R^d,
$$
 where
 $$
      p_0^{(J)} (t, r):= \frac{1}{\phi^{-1}(t)^d} \wedge \frac{1}{r^d \phi (r)} .
 $$
 When a Dirichlet
form $(\EE, \FF)$ is of the form
\eqref{e:ddf} with $D=\br^d$, $A(x)\equiv 0$ and
$J(x,y)$   a symmetric  function
  defined on $\br^d\times \br^d \setminus  {\rm diag}$
that  satisfies the condition
$({\bf J}_{\phi,  \beta})$ for some $\beta\in (0,\infty]$
as well as
the {\bf UJS} condition
when $\beta \in (0,\infty)$  (see Subsection \ref{S:4.2} below for its definition, and \cite{BBK,CKK2,CKK3} for  discussions on  {\bf UJS}),
it is shown in \cite{CKK3} that the transition density function $q(t, x, y)$ of the associated pure jump process
is jointly  continuous and has the two-sided estimates
$$
q(t, x, y) \asymp p_\beta^{(J)} (t, |x-y|)
\quad \hbox{for } (t, x, y) \in (0, \infty) \times \R^d\times \R^d,
$$
where
 for $\b\in (0,1]$,
$$
p_\beta^{(J)} (t, r):=
\begin{cases}
\frac{1}{\phi^{-1}(t)^d}\wedge  \frac{t}{r^d \phi (r) \Phi (r)}&\text{ if } t\le 1,\\
 t^{-d/2}\exp \left(-
r^{\beta} \wedge  ({r^2}/t) \right)&\text{ if  } t>1;
\end{cases}
$$
and for $\b \in (1,\infty]$,
$$
p_\beta^{(J)} (t, r):=
\begin{cases}
\frac{1}{\phi^{-1}(t)^d}\wedge  \frac{t}{r^d \phi (r) \Phi (r)}&\text{ if }   t\le 1, r<1,\\
t \exp\left(- ( r (1 + \log^+({r}/{t}))^{(\beta-1)/{\beta}})\wedge  r^\beta\right)  &\hbox{ if } t\le 1, r \ge 1,\\
t^{-d/2}\exp\left(- ( r (1 + \log^+({r}/{t}))^{(\beta-1)/{\beta}})\wedge   ({r^2}/t)\right)  &\hbox{ if } t> 1.
\end{cases}
$$
See \cite[Theorems 1.2 and 1.4]{CKK3} for more details.
Such processes with tempered jumps at infinity arise
in  statistical  physics to model   turbulence
 as well as in mathematical finance to model stochastic volatility;  see, for example, \cite{MS, Wu}.
When a Dirichlet form $(\EE, \FF)$ is of the form \eqref{e:ddf} with
$D=\br^d$, $A(x) $  uniformly elliptic and bounded,
and $J(x, y)$  a symmetric
function defined on $\br^d\times \br^d \setminus  {\rm diag}$ that  satisfies the condition  $({\bf J}_{\phi,  \beta})$ for some $\beta\in (0,\infty]$
as well as  the {\bf UJS} condition,     one can easily verify by checking the expressions of $p^{(J)}_\beta (t, r)$ in each case
that the estimates in
 Theorem \ref{T1} can be stated as
  \begin{equation}\label{e:1.27a}
\begin{split}
c_0^{-1} & \Big({t^{-d/2}}
\wedge    ( p^{(c)}(t, c_1\rho_D(x,y))+
p_\beta^{(J)}
(t, c_1\rho_D(x,y))  )\Big)\\
&\le  q(t,x,y)\le c_0\Big({t^{-d/2}}
 \wedge   (p^{(c)}(t, c_2\rho_D(x,y))+
p_\beta^{(J)}(t, c_2\rho_D(x,y))) \Big),
\end{split}
\end{equation}
where
$$
p^{(c)}(t, r) :=  t^{-d/2}\, \exp ( - r^2/t).
$$
See  \cite{CKW4} for the corresponding results for symmetric diffusions with jumps on metric measure spaces
where $\beta=0$ but the diffusion part can be sub-diffusive and the upper growth exponent   $\alpha^*$ for $\phi$ in \eqref{eqn:poly} can be possibly
larger than 2.
 We
 further remark that
 the two-sided  heat kernel estimates for diffusions with jumps  of the form \eqref{e:1.20}--\eqref{e:1.12a}
 was first given in \cite{CK3}.
\end{itemize}
\end{remark}

We now
make some comments on the approach of the
main result Theorem \ref{T1}.
Two-sided heat kernel estimates for symmetric diffusions with stable-like jumps on metric measure spaces were obtained in the recent paper \cite{CKW4} by
three of the authors of this paper.
In that paper,
some powerful tools for the study of stability of heat kernel estimates and parabolic Harnack inequalities for symmetric jump processes developed in \cite{CKW1, CKW2} are adapted; on the other hand, a new self-improvement argument for upper bounds via exit time probability estimates is proposed, see \cite[Section 4.3]{CKW4} for details. However, as for symmetric diffusions with
(sub- or super-)exponential
decay jumps in the present paper, the approach of \cite{CKW4} (in particular, the self-improvement argument for upper bounds as mentioned above) does not work, due to light tails of jumps for the associated process.
The novelty of the proof for Theorem \ref{T1} is twofold.

\begin{itemize}
\item Concerning upper bounds, we apply the stability result from \cite{CKW4} to obtain on-diagonal heat kernel estimates, and
obtain
precise tail probability estimates for
truncated processes, the processes with   jumps of  size larger than some  suitable level removed,
to achieve off-diagonal estimates.
Our main technical tool to obtain
tail probability estimates for truncated processes is a new generalization of Davies' method
from our recent work
\cite{CKKW} for Dirichlet heat kernels of the truncated process. The Meyer's construction and the relations among the original process, the truncated process and the killed truncated processes are also fully
 used.

\item
While we rely heavily on techniques from \cite{CKW1, CKW2, CKW4} to study lower bounds
for the heat kernel $q(t, x, y)$,
 a new observation here is that parabolic
 Harnack inequalities for full ranges do not
hold.
In fact under conditions    $({\bf J}_{\phi_1,  \beta_*, \leq })$ and $({\bf J}_{\phi_2, \beta^*, \geq })$ with $\b_* <\b^*$ in $ \{0_+\} \cup (0, \infty] $,
the jumping kernel $J(x,y)$ may not satisfy the  {\bf UJS} condition, see Subsection \ref{S:4.2}. Thus, it follows from (the proof of) \cite[Proposition 3.3]{CKW2} that parabolic Harnack inequalities for full ranges  do not
hold.
Our results in Theorem \ref{T1} in particular give  a family of Feller processes that
satisfy global two-sided heat kernel estimates, but the associated parabolic Harnack inequalities  for full ranges
fail, which is of independent interest. We shall mention that, under condition
$(\bJ_{\phi, 0_+, \leq })$,
we always have the joint H\"older continuity for the heat kernel  $q(t,x,y)$
so that
we can establish two sided estimates for $q(t,x,y)$
for every $t>0$ and $x, y\in {\bar D}$ without introducing any exceptional set.
\end{itemize}

We further
emphasize
that our approach to Theorem \ref{T1} is quite robust
in the sense that, we can
deal with heat kernel estimates for diffusions with jumps
in the
general (VD) setting, and most importantly the case that the jumping kernel $J(x,y)$
  is bounded by different
  weighted functions (in particular,
   with possibly different increasing functions $\phi_1$ and $\phi_2$ and
   possibly different indexes $\beta_*$ and $\beta^*$ in the exponential terms).
   In fact, we have more.
 Roughly speaking,   we have  upper bound estimates $H_{\phi_1, \beta_*}$ for the heat kernel $q(t,x, y)$ under condition
$({\bf J}_{\phi_1, \beta_*, \leq})$ with $\beta_*\in   \{0_+\} \cup (0, \infty] $;
 and we have
 lower bound estimates
   $H_{\phi_2, \beta^*}$ under conditions $({\bf J}_{\phi_1, 0_+, \leq})$
 and $({\bf J}_{\phi_2, \beta^*, \geq})$ with $\beta^*\in (0, \infty]$.
 See Theorems \ref{T1u} and \ref{general} below for more details.
We further mention that, in Theorems \ref{T:1.4} and \ref{T:2.4}, we do not require  the underlying state space $D$ to have infinite volume
 (equivalently,   $D$ to have infinite diameter).
 To  handle the case that diam$(D)$ is bounded, we use localized versions of results in \cite{CKW1, CKW2, CKW4} to obtain heat kernel estimates
 in bounded time intervals; we then
combine a modification of Doeblin's celebrated result with the small time heat kernel bounds to get estimates for large times.

\medskip

The rest of the paper is organized as follows. In the next section, we will give the full proof of Theorem \ref{T:1.4}. In particular, we present the details for the case that the underlying state space $D$ has finite diameter.
In Section \ref{S:3}, we prove Theorem \ref{T:2.4}.
Under
$ ({\bf J}_{\phi, 0_+, \le})$, we show
in Subsection \ref{S:2.1} that
the
(local version of) Faber-Krahn inequality
$\FK (2)$
and the cut-off Sobolev inequality $\CS (2)$
hold for $(\sE, \sF)$. This allows us to get on-diagonal upper bounds
for the heat kernel
$q(t, x, y)$ of $(\sE, \sF)$
and the parabolic H\"{o}lder regularity for caloric functions
 from a recent work \cite{CKW4}.
 In particular, the process associated with $(\sE, \sF)$ can start from every point in $\bar D$
 and is in fact a Feller process having  the strong
Feller property.
By using the results from \cite{CKW2,CKW4}, we further present near-diagonal lower bounds for the heat kernel $q(t, x, y)$ of $(\sE, \sF)$ in Subsection \ref{S:3.1}.
In Section \ref{S:4}, we study off-diagonal heat kernel upper bounds for $(\sE, \sF)$  whose jumping kernel $J(x, y)$ satisfies
 condition
  $({\bf J}_{\phi, \beta_*, \leq})$ with $\beta_* \in \{0_+\} \cup (0, \infty]   $
  when the underlying state space $D$ has infinite diameter.
 This section is technically the most important section in this paper.
 Here, we consider
  the heat kernel
 $  q^{\langle \lambda \rangle} (t, x, y)$
of the finite range Dirichlet form $(\sE^{(\lambda)}, \sF)$ on $L^2(D; m)$
 whose jumping kernel $    J^{\langle \lambda  \rangle }
 (x, y):=   J(x, y) {\1}_{\{\rho_D (x, y) \leq \lambda\}}$.
Upper bound estimates for $    q^{\langle \lambda \rangle}
(t, x, y)$ are derived
in Theorem  \ref{T1u}
through the upper bound estimates
of the Dirichlet heat kernel $  q^{\langle \lambda \rangle, B}
(t, x, y)$ in balls.  The latter uses a recent result from \cite{CKKW}
based on a generalization of Davies' method.
We then derive   upper bound estimates of $  q(t, x, y)$ from that of $
   q^{\langle \lambda \rangle}
 (t, x, y)$ through Meyer's
method of adding/removing jumps.
Lower bound estimates for
$  q(t,x,y)$ are derived in
 Theorem \ref{general} of Section \ref{S:5} under conditions $({\bf J}_{\phi_1,  0_+, \leq})$
 and $ ({\bf J}_{\phi_2, \beta^*, \geq})$ with $\beta^*\in (0, \infty]$ and possibly different
 strictly increasing functions $\phi_1$ and $\phi_2$.
  Its proof
uses the approaches from \cite{CKW1, CKW2, CKW4}. The main result of this paper, Theorem \ref{T1},  follows directly by combining Theorem \ref{T1u} with Theorem \ref{general}
 as well as Theorem \ref{T:2.4}.
Some preliminary integral estimates are given in the
Appendix (Subsection \ref{S:4.1})
of this paper, as
 well as
a brief discussion in Subsection \ref{S:4.2} on
a property concerning the {\bf UJS} condition.

\medskip

\noindent{\bf Notations} \,\, Throughout this paper, we use $c_1, c_2, \cdots$ to denote generic
constants, whose exact values are not important and can  change from
one appearance to another. The labeling of the constants $c_1, c_2,
\cdots$ starts anew in the statement of each result.
For $p\in [1,  \infty]$, we will use $\| f\|_p$ to denote the $L^p$-norm in $L^p(D, m)$.
For any open subset  $U$ of ${\bar D}$, $C_c(U)$ is the space of continuous functions on $U$ having compact support
with respect to the $\rho_D$-metric.

 \section{The case   $\beta_*=\beta^*=0$}\label{S:2}

We first give a general statement that holds under $(\bJ_{\phi, \leq 1})$ given
in Definition \ref{D:1.3}(vi).

\begin{proposition}\label{l:regular}
Suppose that condition  $(\bJ_{\phi, \leq 1})$ holds and $(\sE, \sF)$ is the bilinear form defined by \eqref{e:sF}--\eqref{e:sE}.
Then  there exists  a  constant
 $c>0$,
depending on constants in   $(\bJ_{\phi, \leq 1})$  and \eqref{eqn:poly} for $\phi$  such that
$$
\sE_1 (f,f) \,\le \, c\,\EE^{0, \rf}_1 (f,f) \quad \text{ for every } f \in \FF,
$$
where $\sE_1(f, f):= \sE(f, f)+\|f\|_2^2$ and $\EE^{0, \rf}_1 (f,f):= \EE^{0, \rf} (f,f) + \|f\|_2^2$.
Consequently, $(\sE, \sF)$ is a regular Dirichlet form on $L^2({D}; m)$ and so there is a symmetric Hunt process $Y$
on $\bar D$ associated with it. Moreover, it is conservative in the sense that the Hunt process $Y$ has infinite lifetime.
 \end{proposition}

\pf
Recall that  $Z:=(Z_t)_{t\ge0}$ is the diffusion process  associated with
 the regular Dirichlet form   $(\EE^{0, \rf}, \FF )$  on $L^2(D; m)$,
 and that $\alpha^* \in (0,2)$ is the exponent in \eqref{eqn:poly} for the strictly increasing function $\phi$.
Let $S:=(S_t)_{t\ge0}$  be an
$(\aa^*/2)$-stable subordinator  that is independent of $Z$.
The  subordinated diffusion $ \{Z_{S_t}, \ge0\}$ is $m$-symmetric and we use   $j_{\aa^*}(x,y)$  to denote its   jumping  kernel.
By  condition  $(\bJ_{\phi, \leq 1})$ and Lemma \ref{p:Ialpha1}
from the Appendix of this paper,
we have that for $f \in \FF$,
\begin{align*}
 \sE_1 (f,f)
&=  \EE^{0, \rf} (f,f)+\frac{1}{2}\int_{{D}} \int_{{D}} (f(y)-f(x))^2
   J(x, y) \1_{\{ \rho_D(x,y) \le 1\}} \, m(dx)\, m(dy) \\
 &\quad +\frac{1}{2}\int_{{D}} \int_{{D}} (f(y)-f(x))^2
   J(x, y) \1_{\{ \rho_D(x,y) > 1\}}  \, m(dx)\, m(dy)+
 \|f\|^2_{2}
   \\
 &\le     \EE^{0, \rf} (f,f)+\frac{1}{2}\int_{{D}} \int_{{D}} (f(y)-f(x))^2
   J(x, y) \1_{\{ \rho_D(x,y) \le 1\}} \, m(dx)\, m(dy) \\
 &\quad + \left(\int_{{D}} f(x)^2\,m(dx)\right) \left(\sup_{x\in {D}} \int_{{D}}
   J(x, y) \1_{\{ \rho_D(x,y) > 1\}} \, m(dy)\right)+
 \|f\|^2_{2}  \\
 &\le   \EE^{0, \rf} (f,f)+c_1 \int_{{D}} \int_{{D}} (f(y)-f(x))^2
   j_{\aa^*} (x, y)
  \, m(dx)\, m(dy)+ c_1
 \|f\|^2_{2}.\end{align*}
This along with Lemma \ref{p:Ialpha1--}
in the Appendix   yields that
$$
\sE_1 (f,f) \le c_2 \left(  \EE^{0, \rf} (f,f) +
 \|f\|^2_{2}\right)= c_2 \EE^{0, \rf}_1 (f,f).
$$
The regularity of $(\EE^{0, \rf},
\sF  )$ on $L^2(D; m)$  implies  the regularity of $(\sE, \sF)$ on $L^2(D; m)$.

Let    $  Y$ be the  symmetric Hunt process on ${\bar D}$ associated with the regular Dirichlet form
$(\sE, \sF)$ on $L^2(D; m)$  that can start from
$\sE$-quasi-everywhere in ${\bar D}$. By \cite[Theorem 1]{MUW}, under
(VD)  and condition $(\bJ_{\phi, \leq 1})$, the process $  Y$ is conservative.
 \qed

\medskip

  \noindent{\bf Proof of Theorem  \ref{T:1.4}}.  Recall that, by \cite[Corollary 5.3]{GH},
the underlying state space $D$ has infinite diameter is equivalent to $D$ has infinite volume. Then the proof is divided into two parts.

\smallskip

\noindent
(1) The case that diam$(D)=\infty$.
In this case, $D$ has infinite volume and the desired result follows from the main result of \cite{CKW4}, which we explain below.
In the proof of this part, we
use the Poincar\'e inequality and cut-off Sobolev inequality for general Dirichlet forms. The readers are referred to \cite{AB, CKW1, CKW2, CKW4, GHL} for these definitions (see also Subsection \ref{S:2.1} for some details). First, we write
 \begin{equation}\label{e:prof1.4-1}\EE (u, u)  = \EE^{0,\rf}(u, u)+ \EE^{(J)}(u,u),\end{equation} where
 $$\EE^{(J)}(u,u)=
 \frac{1}{2} \int_{D\times D}
 (u(x)-u(y))^2 J(x, y)\, m(dx)\,m(dy).$$
Recall that  $\Gamma_{0}$ is the
 carr\'e du champ operator corresponding to $(\EE^0,\sF^0)$. Let
 $$\frac{d\Gamma_{J}(u,u)}{dm}(x)=\int_{{D}} (u(y)-u(x))^2
J(x, y) \, m(dy),\quad u\in \sF.$$ Denote by $\Gamma$ the
 carr\'e du champ operator corresponding to $(\EE,\FF)$. Then,
 $$\Gamma(u,u)=\Gamma_{0}(u,u)+\Gamma_{J}(u,u),\quad u\in \sF.$$

According to Theorem \ref{T:1.2} and \cite[Remark 1.11(4)]{AB} (or \cite[p.\ 1942]{GHL}), the Poincar\'e inequality PI$(2)$ and the cut-off Sobolev inequality CS$(2)$ hold for $(\EE^{0,\rf}, \FF)$;
that is, there is a  constant
$c_1>0$ so that for every $x\in {\bar D}$, $r>0$ and
$f\in \FF$,
\begin{align*}
 \min_{a\in\RR}\int_{B_{{\bar D}}(x, r)}(f(y)-a)^2 \,m(dy) \leq & \
 c_1 r^2          \int_{B_{{\bar D}}(x, r)}\Gamma_0 (f, f)(x)\,m(dx);
\end{align*}  and there is  a constant $c_2>0$ such that for every
$0<r\le R$, $x_0\in {\bar D}$ and any $f\in \sF$,
\begin{equation}\label{e:sc-d}
\int_{B_{\bar D}(x_0,R)} f^2 \, d\Gamma_0 (\psi,\psi) \le
\frac{c_2}{r^2}
\int_{B_{\bar D}(x_0,R)} f^2  \,dm,
\end{equation} where $\psi(x)=h(\rho_D(x_0,x))$ is a non-negative cut-off function on $B(x_0,r)\subset B(x_0,R)$ with $h\in C^1([0,\infty))$ being such that $0\le h\le 1$, $h(s)=1$ for all $s\le R$, $h(s)=0$ for $s\ge R+r$ and $|h'(s)|\le 2/r$ for all $s\ge0$.

Furthermore, by
$(\bJ_{\phi})$,
\eqref{eqn:poly} with $\aa^*<2$ and \cite[Remark 1.7]{CKW1}, the following cut-off Sobolev inequality CS$(\phi)$ holds for $(\EE^{(J)}, \FF)$:
for every
$0<r\le R$, $C_0\in (0,1]$, $x_0\in {\bar D}$ and any $f\in \sF$,
$$\int_{B_{\bar D}(x_0,R+(1+C_0)r)} f^2 \, d\Gamma_J (\psi,\psi) \le  \frac{c_2}{\phi(r)}  \int_{B_{\bar D}(x_0,R+(1+C_0)r)} f^2  \,dm,$$ where $\psi(x)=h(\rho_D(x_0,x))$ is a non-negative cut-off function on $B(x_0,r)\subset B(x_0,R)$ as in \eqref{e:sc-d}.
(Note that $\CS(2)$ that appeared before is a special case of $\CS(\phi)$ when $\phi(r)=r^2$
for the corresponding Dirichlet form.)
On the other hand, according to \cite[Theorem 1.13]{CKW1} and \cite[Corollary 1.3 and Theorem 1.18]{CKW2}, we know that the following Poincar\'e inequality PI$(\phi)$ is  also satisfied for $(\EE^{(J)}, \FF)$: there is a constant $c_4>0$ so that for every $x\in {\bar D}$, $r>0$ and
$f\in \FF$,
\begin{align*}
 \min_{a\in\RR}\int_{B_{{\bar D}}(x, r)}(f(y)-a)^2 \,m(dy)
         \leq & \  c_4\phi(r) \int_{B_{{\bar D}}(x, r)}\int_{B_{{\bar D}}(x, r)}(f(x)-f(y))^2J(x,y)\,m(dx)\,m(dy).
\end{align*}
($\PI(2)$ that appeared before is a special case of $\PI(\phi)$ when $\phi(r)=r^2$
for the corresponding Dirichlet form.)
 Here we note that, thanks to Proposition \ref{l:regular} again, $(\EE^{(J)}, \FF)$ is a regular Dirichlet form.

  Combining all the conclusions with \eqref{e:prof1.4-1}, we know that the Poincar\'e inequality PI$(\bar\phi)$ and the cut-off Sobolev inequality CS$(\bar\phi)$ hold for $(\EE, \FF)$, where $\bar \phi(r)
  =r^2\wedge \phi(r)$. That is, there
is a  constant
$c_5>0$ so that for every $x\in {\bar D}$, $r>0$ and
$f\in \FF$,
\begin{align*}
 \min_{a\in\RR}\int_{B_{{\bar D}}(x, r)}(f(y)-a)^2 \,m(dy)
         \leq & \  c_5 \bar\phi(r)\bigg(
     \int_{B_{{\bar D}}(x, r)}\Gamma_0 (f, f)(x)\,m(dx)\\
    &\qquad\qquad +\int_{B_{{\bar D}}(x, r)}\int_{B_{{\bar D}}(x, r)}(f(x)-f(y))^2J(x,y)\,m(dx)\,m(dy)\bigg);
\end{align*} and for every
$0<r\le R$, $C_0\in (0,1]$, $x_0\in {\bar D}$ and any $f\in \sF$,
$$\int_{B_{\bar D}(x_0,R+(1+C_0)r)} f^2 \, d\Gamma (\psi,\psi) \le  \frac{c_6}{\bar\phi(r)}  \int_{B_{\bar D}(x_0,R+(1+C_0)r)} f^2  \,dm,$$ where $\psi(x)=h(\rho_D(x_0,x))$ is a non-negative cut-off function on $B(x_0,r)\subset B(x_0,R)$ as in \eqref{e:sc-d}.
This together with condition $(\bJ_{\phi})$
and \cite[Theorem 1.13 and Theorem 1.18]{CKW4} yields the desired assertion.

\smallskip

\noindent
(2) The case that diam$(D)<\infty$. In particular, $\bar D$ is compact.
(This is because, by (VD), bounded sets are totally bounded, and
totally bounded closed sets on
 complete metric spaces
are compact.)
As mentioned in \cite[Remark 1.19]{CKW2}, with minor adjustments of the proofs
localized versions of the results
 in \cite{CKW1,CKW2,CKW4} should continue to hold.
 In particular, the heat kernel  $q(t, x, y)$ of  $(\EE, \FF)$ is jointly H\"{o}lder continuous on $(t,x,y)$, and \eqref{e:1.12a} holds for all $t\in (0, t_0]$ and $x,y\in \bar D$, where $t_0=\bar\phi(r_0)$ for some $r_0\in (0, {\rm diam}(D)]$. The readers are referred to \cite{GHH} for heat kernel estimates for
stable-like
Dirichlet forms with the Ahlfors $d$-set condition in both bounded and
unbounded cases.  Next, we will claim that for all $t\ge t_0$ and $x,y\in \bar D$,
\begin{equation}\label{e:prof-1.4.3}
q(t,x,y)\simeq1.\end{equation}

Since $\bar D$ is compact,
thanks to the joint H\"{o}lder continuity of $q(t,x,y)$ and the fact that \eqref{e:1.12a} holds for $t\in (0, t_0]$ and $x,y \in \bar D$, the process $Y$ associated with
$(\ce, \cf)$ is irreducible and has the strong Feller property. Hence, by a modification of Doeblin's celebrated result (see \cite[p.\ 365, Theorem 3.1]{BLP}), the process $Y$ is strongly exponentially ergodic in the sense that
 there exist positive constants $\lambda_1$ and
$c_0$ so that
\begin{equation}\label{e:1.28}
\sup_{x\in \bar D}|\bE_xf(Y_t)-\bar m (f)|\le c_0e^{-\lambda_1
t}\|f\|_\infty   \quad \hbox{for every } t>0 \hbox{  and }  f\in B_b(\bar D),
\end{equation}
where $\bar m=m(D)^{-1}m$ is the unique
 invariant probability measure of the process $Y$. For any $y\in \bar D$, applying $f(z)=q(t_0,z,y)$ into \eqref{e:1.28}, we find that for any $t>0$,
 \begin{equation}\label{e:prof1.4-2}\left| q(t+t_0,x,y)- \frac{1}{m(D)}\right|\le c_1e^{-\lambda_1t},\end{equation} where
 $c_1=c_0\sup_{z,u\in \bar D} q(t_0,z,u)<\infty,$ due to the fact that \eqref{e:1.12a} holds for $t\in (0, t_0]$ and $x,y \in \bar D$. Here, we used the semigroup property of $q(t,x,y)$ and the conservativeness of the process $Y$. According to \eqref{e:prof1.4-2}, for all $t>0$ and $x,y\in \bar D$,
    \begin{equation}\label{e:prof1.4-2-1}
 q(t+t_0,x,y)\le \frac{1}{m(D)}+ c_1e^{-\lambda_1t}\le c_1 +\frac{1}{m(D)}=:c_2
 \end{equation}
 and
 $$ q(t+t_0,x,y)\ge \frac{1}{m(D)}- c_1e^{-\lambda_1t}.$$
In particular, we obtain that for all $t\in (t_0,\infty)$ and $x,y\in \bar D$, $q(t,x,y)\le c_2$, and that for all $t\in [t_1,\infty)$ and $x,y\in \bar D$, $q(t,x,y)\ge \frac{1}{2m(D)}$, where $t_1>t_0$ satisfies that $2m(D)c_1\le e^{\lambda_1(t_1-t_0)}$.
Furthermore,  for any $t\in (t_0, t_1]$, we take $n\ge 1$ such that $t':=t/n\in[t_0/2, t_0]$. Then, by the semigroup property of $q(t,x,y)$, for any $x,y\in \bar D$,
\begin{align*}
q(t,x,y)=&\int_{\bar D}\cdots \int_{\bar D} q(t',x,z_1)q(t',z_1,z_2)\cdots q(t',z_{n-1},y)\,
m(dz_1)\,m(dz_2)\,\cdots m(d z_{n-1})\\ \ge& \left(  m(D)
\Big( \inf_{s\in [t_0/2,t_0]}\inf_{z,u\in \bar D} q(s,z,u)\Big)\right)^{n-1}\\
\ge & \left( 1\wedge \Big( m(D)
\Big( \inf_{s\in [t_0/2,t_0]}\inf_{z,u\in \bar D} q(s,z,u)\Big)\Big)\right)^{2t_1/t_0}=:c_3>0,
\end{align*}
where
$c_3$ may depend on $m(D)$. Here
in the last inequality, we used again the fact that \eqref{e:1.12a} holds for all $t\in (0,t_0]$ and $x,y\in \bar D$.
Putting all the estimates above together,
 we get \eqref{e:prof-1.4.3}, and, consequently, the desired assertion when diam$(D)<\infty$.
 \qed

\section{General framework under condition   $({\bf J}_{\phi, 0_+, \leq })$} \label{S:3}

In this section,  $(\sE, \sF)$  is the bilinear form defined by \eqref{e:sF}--\eqref{e:sE} with $J(x,y)$ satisfying condition $({\bf J}_{\phi, 0_+, \leq })$
of Definition \ref{D:1.3}(iv),
 where $\phi$ is a strictly increasing function on $[0, \infty)$ satisfying  $\phi  (0)=0$, $\phi (1)=1$ and \eqref{eqn:poly}.
 By  \eqref{e:Jcom} and Proposition \ref{l:regular},
 $(\sE, \sF)$ is a regular Dirichlet form on $L^2(D; m)$ and there is a conservative symmetric Markov process $Y$
 associated with it.  We will show   the existence of
 the heat kernel
$q(t,x,y)$ for the Hunt process $  Y$ and its joint H\"older regularity, and obtain on-diagonal upper bound estimates as stated in Theorem \ref{T:2.4}
as well as its near diagonal lower bound estimates for $q(t,x,y)$.

\smallskip

\subsection{On-diagonal upper bound estimates}\label{S:2.1}
Recall that  $\Gamma_{0}$ is the
 carr\'e du champ operator corresponding to
  the strongly local regular Dirichlet form  $(\EE^0,\sF^0)$ on $L^2(E; m)$.
  The inequality \eqref{e:1.3} in particular gives the following
Poincar\'e inequality PI($2$) of \eqref{pinow}  for the regular Dirichlet form
 $(\sE, \sF)$ on $L^2(D; m)$, which  in turn implies  (see, e.g.\,\cite[Proposition 7.6]{CKW1})
  the local Nash-type inequality \eqref{E:Na}  for  $(\sE, \sF)$.

 \begin{lemma}\label{L:local-nash}
  Suppose that  $({\bf J}_{\phi, 0_+, \leq })$ holds. Then there
 is a  constant
$c_1>0$ so that for every $x\in {\bar D}$, $r>0$ and
$f\in \FF$,
\begin{align}
 \min_{a\in\RR}\int_{B_{{\bar D}}(x, r)}(f(y)-a)^2 \,m(dy) \leq & \
          c_1 r^2\bigg(
     \int_{B_{{\bar D}}(x, r)}\Gamma_0 (f, f)(x)\,m(dx)\nonumber\\
    &\quad\quad +\int_{B_{{\bar D}}(x, r)}\int_{B_{{\bar D}}(x, r)}(f(x)-f(y))^2J (x,y)\,m(dx)\,m(dy)\bigg).
     \label{pinow}
\end{align}
Consequently, there is a constant $c_2>0$ such that for any $r>0$ and  $u\in \sF\cap L^1({D};m)$,
\begin{equation}\label{E:Na}
\|u\|_2^2\le c_2\Big(\frac {\|u\|_1^2}{\inf_{z\in {\rm supp}[u]}V_D(z,r)}
+r^2\sE (u,u)\Big).
\end{equation} \end{lemma}

Note that under
general (VD)   condition  on the metric measure space $(\bar D, \rho_D, m)$,
one can not expect
the global
Nash-type  inequality such as that
of  \cite[Theorem 2.5]{CK2} to hold. In the present setting,
 we can only expect the local
 Nash-type inequality \eqref{E:Na} as shown in Lemma \ref{L:local-nash}.

\smallskip

 We next introduce a localized version of Faber-Krahn inequality.
For any open set $U \subset {\bar D}$, let $\sF_U$   be the
$\sE_1$-closure in $\sF$ of  $\sF\cap C_c(U)$.
Define
$$ \lam_1(U)
= \inf \left\{ \sE(f,f):  \,  f \in \sF_U \hbox{ with }  \|f\|_2 =1 \right\} .
 $$

\begin{definition} \rm We say that the
 (localized version of) {\em Faber-Krahn inequality} $\FK(2)$
holds for
the regular Dirichlet form $(\sE, \sF)$ on $L^2(D; m)$
if there exist
positive constants $\sigma\in (0,1/2)$,
$C$ and
$\nu$ such that
 for any $r \in (0, \sigma $diam$(D))$,  any ball $B_{{\bar D}}(x,r)$
and any open set $U \subset
B_{{\bar D}}(x,r)$, $$
 \lam_1 (U) \ge
 \frac{C}{r^2}
 (V_D(x,r)/m(U))^{\nu}.$$
\end{definition}

In view of
($\RVD$) of \eqref{univd*} and the  Poincar\'e inequality $\PI(2)$ of \eqref{pinow}
  for $(\sE, \sF)$,
we  conclude by
 the proof of  \cite[Proposition 7.4]{CKW1}
(see also the proof of \cite[Lemma 3.5]{GHH})
that  $\FK(2)$ holds  for $(\sE, \sF)$.

\smallskip

We next show that  the cut-off Sobolev inequality  {\rm CS($2$)}  holds for $(\sE,\sF)$;
see \cite[Definition 1.6]{CKW4} for its definition.
Recall again that $\Gamma_{0}$ is the
 carr\'e du champ operator corresponding to $(\EE^0,\sF^0)$. A function $f$ is said to be locally in $\FF^0$, denoted as $f\in \FF^0_{\loc}$,
 if for every relatively compact subset $U$ of $E$, there is a function $g\in \FF^0$
 such that $f=g $ $m$-a.e. on $U$.
It is well known that the following chain rule holds (see \cite[Theorem 4.3.7]{CF} and \cite[Theorem 3.2.2]{FOT}):
for every $C^1$ function
$h:\R\to \R$
and bounded function $f \in
\FF^0_{\loc}$,
 $$
\Gamma_0 ( h(f ), h(f ))=    (h'(f ))^2  \, \Gamma_0(f , f ).
$$
We remark  that the assumption $h(0)=0$ in \cite[Theorem 4.3.7]{CF} and \cite[Theorem 3.2.2]{FOT}
 is not needed since, by the strongly local property of $(\sE^0, \sF^0)$ and the fact that $1 \in
 \FF^0_{\loc}$,
$$ \Gamma_0 ( h(f )-h(0), h(f )-h(0))=\Gamma_0 ( h(f ), h(f )) .
$$
In particular, for every  bounded function $f \in
\FF^0_{\loc}$,
$
\Gamma_0 ( e^f, e^f)=e^{2f} \Gamma_0 (f, f).$
 Let
 $$\frac{d\Gamma_{J }(u,u)}{dm}(x)=\int_{{D}} (u(y)-u(x))^2
J (x, y) \, m(dy),\quad u\in \sF.$$ Denote by $\Gamma$ the
 carr\'e du champ operator corresponding to $(\sE,\sF)$. Then,
 $$\Gamma(u,u)=\Gamma_{0}(u,u)+\Gamma_{J }(u,u),\quad u\in \sF.$$
 For any fixed $x_0\in \bar D$ and $r,R>0$, we choose a non-negative cut-off function $\psi(x)=h(\rho_D(x_0,x))$, where $h\in C^1([0,\infty))$ such that $0\le h\le 1$, $h(s)=1$ for all $s\le R$, $h(s)=0$ for $s\ge R+r$ and $|h'(s)|\le 2/r$ for all $s\ge0$. It holds that for all $x\in {\bar D}$,
 $$
 \frac{d\Gamma_{0}(\psi,\psi)}{dm}(x)=h'(\rho_D(x_0,x))^2\frac{d\Gamma_{0}(\rho_D(x_0,\cdot),\rho_D(x_0,\cdot))}{dm}(x) \le 4/r^2.
 $$
  On the other hand,
under condition $({\bf J}_{\phi, 0_+, \leq})$,
 \begin{align*}
\frac{d\Gamma_{J } (\psi, \psi)}{dm } (x)&=\int_{{D}} (\psi(x)-\psi(y))^2 J (x,y)\,m(dy)
\le \frac{4}{r^2}\int_{D} \rho_D(x,y)^2 J (x,y) \,m(dy) \le \frac{c_1}{r^2}.
\end{align*}Hence,
$$ \frac{d\Gamma(\psi, \psi)}{dm }
\leq  \frac{c_2}{r^2 }.
$$
This yields that for every
$0<r\le R$, $C_0\in (0,1]$, $x_0\in {\bar D}$ and any $f\in \sF$,
$$
\int_{B_{\bar D}(x_0,R+(1+C_0)r)} f^2 \, d\Gamma (\psi,\psi) \le
\frac{c_2}{r^2}
\int_{B_{\bar D}(x_0,R+(1+C_0)r)} f^2  \,dm.
$$
The above inequality in particular implies that the cut-off Sobolev inequality  {\rm CS($2$)}  holds for $(\sE,\sF)$.
 The cut-off Sobolev inequality CS($2$) guarantees the existence of good cut-off functions and it plays an important role
 when the walk dimension  of the strongly local Dirichlet form $(\EE^0, \FF^0)$ on $L^2(E; m)$ is strictly larger than 2.
 In the framework of this paper,    the walk dimension of $(\EE^0, \FF^0)$ is 2 by Theorem \ref{T:1.1}(ii)
 and it is pretty easy to establish  {\rm CS($2$)}  of $(\sE,\sF)$ as shown above.
 This will enable us to take advantage of using results recently obtained in \cite{CKW4} in our study of the heat kernel of $(\sE,\sF)$.
We refer the reader to   \cite[Remark 1.7]{CKW4}  and \cite[Definition 1.5 and Remark 1.6]{CKW1}   for  more information
about the cut-off Sobolev inequality for symmetric diffusion with jumps and for  pure jump Dirichlet forms, respectively.

\medskip

Recall that $\phi_*$ is the strictly increasing function on $[0, \infty)$ defined by \eqref{e:phi*},
and  $p^{(c)}(t,x,r)$ and $ p^{(j)}_{\phi_*}(t,x,r)$ are the functions defined by \eqref{e:pc} and \eqref{eqn:5} (with $\phi_*$ in place of $\phi$), respectively.

\medskip

 \noindent {\bf Proof of Theorem \ref{T:2.4}.}
 Since $(\sE, \sF)$ is a regular Dirichlet form on $L^2(D; m)$, there is an associated
 symmetric Hunt process $  Y$ that starts from every point from $\bar D \setminus  {\cal N}$
  for a   properly $\sE$-exceptional set ${\cal N}\subset  \bar D $.
 As   verified in the above,
   the  localized version of  $\FK(2)$ and
 $\CS(2)$ hold for the Dirichlet form $(\sE,\sF)$.
As  $(\bJ_{\phi, 0_+, \leq})$ implies  $(\bJ_{\phi, \le})$, we have by the proof of \cite[Theorem 1.14]{CKW4}, \cite[Remark 1.19]{CKW2} and \cite{GHH}
that
 the Hunt process $  Y$,  by enlarging the $\sE$-exceptional set ${\cal N}\subset {\bar D}$ if needed,
  has a transition density function $  q(t,x, y)$  defined on $(0, \infty) \times
 (\bar D \setminus {\cal N} ) \times (\bar D \setminus {\cal N} )$ so that
 \eqref{e:hk00} holds for every $x, y \in \bar D \setminus {\cal N}$ and
 $t>0$ ($t \in (0, t_0]$ for some $t_0>0$ if ${\rm diam}(D) <\infty$)
 and that \eqref{e:2.14} holds for every $x, y \in \bar D \setminus {\cal N}$ and
  $r\in (0, {\rm diam}(D)]$.
 Here, we used the fact that the scaling function the process $ Y$ defined by \cite[(1.14)]{CKW4} is comparable
 to $r^2$.

 On the other hand,    by $({\bf J}_{\phi, 0_+, \leq})$ and
  \cite[Theorem 3.1]{CKK2}, we know that
 any bounded  caloric function  $u(t, x)$ of
$  Y$  (locally in $t>0$) is  jointly  H\"{o}lder continuous.
Thus $   q(t,x, y)$ can be refined to be a jointly H\"older continuous function defined
 pointwise on
$(0, \infty) \times  \bar D  \times  \bar D  $ and  the upper bound estimate
\eqref{e:hk00} holds without the exceptional set ${\cal N}$.
Consequently, the Hunt process $  Y$ can be refined to be a Feller process having strong Feller property
that starts from every point on $\bar D$, and \eqref{e:2.14} holds for every $x\in \bar D$ and
$r\in (0, c_3 {\rm diam}(D)]$
with any $c_3\in (0, 1/2)$.

 When ${\rm diam}(D) <\infty$,  the upper bound in  \eqref{e:hk00} for $t > t_0$ is simply a positive constant, which follows from the same argument as the one in
 \eqref{e:1.28}--\eqref{e:prof1.4-2-1}.
   \qed

\subsection{Near diagonal lower bound estimates}\label{S:3.1}
In this subsection,  we
present the sketch of the proof of near diagonal lower bound estimates for the heat kernel $  q(t,x,y)$.
As shown in Subsection \ref{S:2.1},
the symmetric Hunt process $Y$, or equivalent, its associated regular Dirichlet form $(\sE, \sF)$ on $L^2(D; m)$,
satisfies $\PI(2)$, $\CS(2)$ and
 $(\bJ_{\phi, \le})$. Thus,  by \cite[Proposition 5.2]{CKW4} and its proof,
 the elliptic H\"{o}lder regularity (EHR) holds for the process $ Y$; that is, there exist constants $c>0$, $\theta\in (0,1]$ and $\varepsilon\in (0,1)$ such that for every $x_0\in {\bar D}$,
  $r \in (0, $\rm diam$(D)/3)$
  and every
  bounded measurable function $u$ on ${\bar D}$ that is harmonic in $B_{{\bar D}}(x_0,r)$
 with respect to the Hunt process $Y$,
 $$|u(x)-u(y)|\le c\left(\frac{\rho_D(x,y)}{r}\right)^\theta {\rm ess\,sup}_{{\bar D}}|u|$$ for any $x,y \in B_{{\bar D}}(x_0,\varepsilon r) $.
 On the other hand,
 two-sided exit time estimates   \eqref{e:2.14}  hold for $ Y$.
 With those two
  properties,
  we can follow
 the proof of \cite[Proposition 4.10]{CKW2}
 to get the following proposition.
 Recall that $ q(t,x,y)$ is the transition density function of $Y$.
  For any open set $U\subset {\bar D}$, denote by $ q^U(t,x,y)$ the
  transition density function of the subprocess $Y^U$ of $Y$ killed up exiting $U$, which is the
  (Dirichlet) heat kernel associated with the part Dirichlet form $(\sE, \sF_U)$.

\begin{proposition}\label{P:3.1}
Under  $({\bf J}_{\phi, 0_+, \leq})$, there exist $\varepsilon\in (0,1)$ and $c_1>0$ such that for any $x_0\in {\bar D}$,
 $r \in (0, $\rm diam$(D)/3)$
 and $0<t\le
 r^2$,
$$ q^{B_{{\bar D}}(x_0,r)}(t,x,y)\ge \frac{c_1}{V_D(x_0,\sqrt{t})},\quad x,y\in B_{{\bar D}}(x_0,\varepsilon t^{1/2}).$$
 {\it In particular, there exist
$c_2, c_3>0$ such that for any
   $t \in (0,  \, {\rm diam} (D)^2/9)$
 and any $x,y\in {\bar D}$ with $\rho_D(x,y)\le c_2t^{1/2}$,
 \begin{equation}\label{e:ee-} q(t,x,y)\ge \frac{c_3}{V_D(x,\sqrt{t})}.\end{equation}}
\end{proposition}

\begin{remark}\label{Add:oneremark}\rm
Note that PI($2$) of \eqref{pinow} follows directly from \eqref{e:1.3} as $\sE(u,u)\ge \EE^{0,\rf}(u, u)$,
so it holds under $(\bJ_{\phi,  0_+, \leq 1})$.
Under condition $(\bJ_{\phi,  0_+, \leq 1})$ of  Definition \ref{D:1.3}(v),
 we can easily see that $\CS(2)$ is satisfied. With these at hand,
 all assertions in this section still hold true under condition  $(\bJ_{\phi,  0_+, \leq 1})$
except that \eqref{e:hk00} holds for all $t\in (0, 1]$ and $x, y\in \bar D$ with $\rho_D(x, y)
\leq 1$, \eqref{e:2.14} holds for all $x\in \bar D$ and
 $r\in (0, 1 \wedge $\rm diam$(D)/3$), and \eqref{e:ee-} holds for all
 $t \in (0, 1\wedge $\rm diam$(D)^2/9)$
 and any $x,y\in {\bar D}$ with $\rho_D(x,y)\le c_2t^{1/2}$.
\end{remark}

\section{Upper bound estimates   under $(\bJ_{\phi_1, \beta_*, \leq})$    with
 $\beta_* \in (0, \infty]$
}\label{S:4}

 In  this   section, the inner uniform domain $D$ is assumed to be unbounded.
  We will  derive  full   upper bounds of    the heat kernel
  associated with the Dirichlet form $(\sE,
\sF)$ of \eqref{e:sF}--\eqref{e:sE}  (or equivalently, of its associated Hunt    process $  Y$)
under the assumption $(\bJ_{\phi_1, \beta_*, \leq})$ for some $\beta_* \in (0, \infty]$.
Recall that $\phi _1$ is a strictly increasing function on $[0, \infty)$
satisfying $\phi_1  (0)=0$, $\phi_1  (1)=1$
 and  \eqref{eqn:poly}.
 In view of \eqref{e:Jcom},
   all the assertions in Section \ref{S:3} hold in this section under $(\bJ_{\phi_1, \beta_*, \leq})$.

  To study off-diagonal upper bounds
 for
  the transition density function  $ q(t,x,y)$ of $ Y$,
 we will use the following Meyer's construction  \cite{Mey75}
 which is to decompose the process by removing (large) jumps.
For any $\lambda>0$, set
\[
 J^{\langle \lambda \rangle}(\xi,\eta):= \1_{\{ \rho_D(\xi,\eta)\leq
\lambda \} } \,  J (\xi, \eta)~~\mbox{ and }~~
 J_{\langle \lambda \rangle}(\xi,\eta):= \1_{\{ \rho_D(\xi,\eta) > \lambda
\} } \,  J (\xi, \eta) .
\]
One can remove the jumps of $ Y$ of size larger than $\lambda$
to obtain a new process $ Y^{(\lambda)}$ as follows.
For each $x\in \bar D $, start a copy
 $ Y^1$
 of the process $ Y$ with initial position $x$.
Run it
until the stopping time
$$
T_1:=\inf \big\{t>0: \rho_D(  Y^{1}_{t-}, Y^{1}_t
)>\lambda \big\},
$$
and define
$   Y^{(\lambda)}_t :=  Y^{1}_t
 $ for $t\in [0, T_1)$.
If $T_1<\infty$,  let $ Y^{2}$
 be a copy of $ Y$ that is independent of  $ Y^{1}$
  and starts from  the point $ Y_{T_1-}$.
Let
$$
T_2:=\inf \big\{t>0: \rho_D(
 Y^{2}_{t-}, Y^{2}_t
)>\lambda \big\},
$$
and define  $   Y^{(\lambda)}_{T_1+t} :=
  Y^{2}_t
 $ for $t\in [0, T_2)$.
Repeat this procedure until one of $T_n$'s is  infinity or countably many times if all $T_n<\infty$.
Let $\zeta^{(\lambda)}=\infty$ in the first case and $\zeta^{(\lambda)}=\sum_{k=1}^\infty T_k$ in the second case.
Meyer \cite{Mey75} showed that
the resulting process $ Y^{(\lambda)}$ is a
Hunt process having    jumping kernel $ J^{\langle \lambda \rangle}(x,y)$ and lifetime $\zeta^{(\lambda)}$
that   can start from every point in ${\bar D}$.
It is easy to see that the Dirichlet form
 of $ Y^{(\lambda)}$
  on $L^2(D; m)$ is $( \sE^{(\lambda )}, \sF)$, where
$$
    \sE^{(\lambda)} (v, v)=
    \frac{1}{2}\int_{D} \Gamma_0 (v, v) (x)\,m(dx)+\frac{1}{2}\int_{D} \int_{D} (v(\xi)-v(\eta))^2\,
 J^{\langle \lambda \rangle} (\xi, \eta) \,
 m(d\eta)  \, m(d\xi).
$$
Note that by
$(\bJ_{\phi_1, \beta_*, \leq})$
and Lemma \ref{intelem}(i),
$$ \sup_{\eta\in {D}}\displaystyle \int_{D}  J_{\langle \lambda \rangle}(\eta, \xi)\, m(d\xi)
\le c_1(\lambda)<\infty.$$
  Thus, for $v\in \sF$,
$$
 0\leq \sE(v, v)- \sE^{(\lambda)} (v, v)
  \leq
  2 \int_{D} v(\xi )^2  \int_{D}   J_{\langle \lambda \rangle}
   (\xi, \eta) \,m(d\eta)   \,m(d\xi)
   \leq
  2 c_1(\lambda)
   \, \int_{D} v(\xi )^2  \,m(d \xi)
   $$
   and so
$$
 \big( 1+
  2 c_1(\lambda)\big)^{-1}
 \sE_1 (v, v)\leq  \sE_1^{(\lambda )} (v, v) \leq \sE_1 (v, v)
\quad \hbox{for every } v\in \sF.
$$
Thus a set $A\subset \bar D$ is  $ \sE^{(\lambda )}$-polar if and only if it is $\sE$-polar.

According to \eqref{e:hk00},
  $(\bJ_{\phi_1, \beta_*, \le})$
and the proof of \cite[Lemma 5.2]{CKW1},
the process $ Y^{(\lambda)}$ admits the transition density function
$ q^{\langle \lambda \rangle}(t, x, y)$ defined on $(0, \infty)\times {\bar D}
\times {\bar D}$, which satisfies that for all
$t>0$
 and all $x,y\in {\bar D}$,
$$
 q^{\langle \lambda \rangle}(t,x,y)\le c_2\left(\frac{1}{V_D(x,\sqrt{t})}+\frac{1}{V_D(y,\sqrt{t})}\right)\exp\left(\frac{c_3t}{\lambda^2}-\frac{c_4\rho_D(x,y)}{\lambda}\right).
$$
Consequently, by \eqref{e:1.8},   for all
$t>0$
 and $x,y\in {\bar D}$,
\begin{equation}\label{e:ondia-u--}
 q^{\langle \lambda \rangle}(t,x,y)\le  \frac{c_5}{V_D(x,\sqrt{t})}\left(1+\frac{\rho_D(x,y)}{\sqrt{t}}\right)^{d_2}\exp\left(\frac{c_3t}{\lambda^2}-\frac{c_4\rho_D(x,y)}{\lambda}\right).
\end{equation}
On the other hand, there is a close relation between $ q(t,x,y)$ and $ q^{\langle \lambda \rangle}(t,x,y)$ in the viewpoint of Meyer's decomposition. In particular, according to
$(\bJ_{\phi_1, \beta_*, \leq})$
and \cite[(4.34), Proposition 4.24
and
its proof]{CKW1}, for any $t,\lambda>0$ and any $x,y\in {\bar D}$.
\begin{equation}\label{e:up-01}
  q(t,x,y)\le  q^{\langle \lambda \rangle}(t,x,y)+ \frac{c_6 t}{V_D(x,\lambda)\phi_1(\lambda)
   \exp (c_7\lambda^{\beta_*})}.
 \end{equation}

\medskip

In order to obtain sharp upper bounds of $ q(t,x,y)$, we need to
refine the estimate \eqref{e:ondia-u--} for $ q^{\langle \lambda \rangle}(t,x,y)$. For this, we will use the following bounds for the Dirichlet heat kernel $ q^{\langle \lambda \rangle, B_{{\bar D}}(x_0,\lambda)}(t,x,y)$
of $ Y^{(\lambda)}$ in $B_{{\bar D}}(x_0,\lambda)$,
which are based on a generalized Davies' method recently developed in \cite{CKKW}.
For an open set
$U\subset {\bar D}$, let $ q^{\langle \lambda \rangle, U}(t,x,y)$ be the (Dirichlet) heat kernel of
the subprocess $ Y^{(\lambda), U}$ of $ Y^{(\lambda)}$ killed up exiting $U$; that is,
$$
\bE_x \big[ f(Y_t^{(\lambda)}):  \, t<    \tau_{U}^{(\lambda)} \big]
=\int_{U} f(y) q^{\langle \lambda \rangle, U}(t,x,y)\,m(dy),\quad t>0, \ x\in U, \ f\in L^2(U;m),
$$
where $ \tau_U^{(\lambda)}$ is the first exit time from $U$
by
the process $ Y^{(\lambda)}$. Recall that in this section, we assume that diam$(D)=\infty$. Then, by
\eqref{univd}, \eqref{univd*}, \eqref{e:hk00} and \cite[Theorem 5.1]{CKKW},
for any $\beta_* \in (0,  \infty]$ and $l\ge 2$,
there exists
a constant
$c_0 >0$ such that for any $x_0\in {\bar D}$,
 $\lambda > 0 $,
 $f \in \Lip_{c}({\bar D})$, $t>0$ and any $x,y\in B_{{\bar D}}(x_0,l\lambda)$,
\begin{equation}\label{e:3.2qw1}
 q^{\langle \lambda \rangle, B_{\bar D}(x_0,l\lambda)}(t,x,y)\leq
\frac{c_0}{V_D(x_0,\lambda)}\,\left(\left(\frac{\lambda}{\sqrt{t}}\right)^{d_1}\vee \left(\frac{\lambda}{\sqrt{t}}\right)^{d_2}\right)\;\exp\left( -|f(y)-f(x)|+2
\Lambda^{\langle \lambda \rangle}(f)^2\; t\right),\end{equation}
where $d_1, d_2>0$ are the constants in \eqref{univd*} and  \eqref{univd} respectively, and
$$
\Lambda^{\langle \lambda \rangle} (f)^2  =\|e^{-2f} \Gamma_{\langle \lambda \rangle}(e^{f}) \|_\infty \vee
\|e^{2f}  \Gamma_{\langle \lambda \rangle}(e^{-f})\|_\infty,$$
with $$
 \Gamma_{\langle \lambda \rangle}(f ) (\xi)  =  \Gamma_0 (f, f) (\xi) + \int_{D} \left( {f(\xi)-f(\eta)}\right)^2
 J^{\langle \lambda \rangle} (\xi, \eta)\,
 m(d\eta).$$

\begin{proposition}\label{P:2.41} Suppose that
$(\bJ_{\phi_1, \beta_*, \leq})$ holds  for some $\beta_*\in (0, \infty]$.
Then the following hold.
\begin{itemize}
\item[{\rm(i)}] If $\beta_* \in (0,  1]$, then, for any $l\ge 2$, there exist $c_1,c_2$ and $C_*>0$ such that for all $x_0\in {\bar D}$, $\lambda>0$, $t>0$ and $x, y\in B_{{\bar D}}(x_0,l\lambda)$,
\begin{align*}
 q^{\langle \lambda \rangle, B(x_0,l\lambda)} (t,x,y)\le
\frac{c_1}{V_D(x_0,\lambda)}\,\left(\left(\frac{\lambda}{\sqrt{t}}\right)^{d_1}\vee\left(\frac{\lambda}{\sqrt{t}}\right)^{d_2}\right)
\times \begin{cases}
e^{-c_2\rho_D(x,y)^2/t}&
\mbox{if }\,C_*\lambda^{1-\beta_*}\rho_D(x,y) <t,\\ e^{-c_2\lambda^{\beta_*-1}\rho_D(x,y)} &\mbox{if }\,
C_*\lambda^{1-\beta_*}\rho_D(x,y)  \ge t.
\end{cases}
\end{align*}

\item[{\rm(ii)}] If $\beta_*\in
 (1, \infty)$, then, for every $l\ge 2$ and
 $C_* \in (0,1)$,
 there exist $c_1,c_2, c_3>0$ such that for all $x_0\in {\bar D}$, $\lambda>0$, $t>0$ and $x, y\in B_{{\bar D}}(x_0,l\lambda)$ with $ \rho_D(x,y) \leq t/C_*$,
$$
 q^{\langle \lambda \rangle, B(x_0,l\lambda)} (t,x,y)\le
\frac{c_1}{V_D(x_0,\lambda)}\,\left(\left(\frac{\lambda}{\sqrt{t}}\right)^{d_1}\vee \left(\frac{\lambda}{\sqrt{t}}\right)^{d_2}\right)
 \exp \left(-\frac{c_2\rho_D(x,y)^2}t \right);
$$
and that for all $x_0\in {\bar D}$, $\lambda>0$, $t>0$ and $x, y\in B_{{\bar D}}(x_0,l\lambda)$ with $ \rho_D(x,y) > t/C_*$,
\begin{align*}
 q^{\langle \lambda \rangle, B(x_0,l\lambda)} (t,x,y)&\le
\frac{c_1}{V_D(x_0,\lambda)}\,\left(\left(\frac{\lambda}{\sqrt{t}}\right)^{d_1}\vee \left(\frac{\lambda}{\sqrt{t}}\right)^{d_2}\right)
 \exp \left(-c_3\rho_D(x,y) \left(\log\frac{\rho_D(x,y)}{t}\right)^{\frac{\beta_*-1}{\beta_*}} \right).
\end{align*}

\item[{\rm(iii)}]
Suppose that the jumping kernel $ J(x,y)$ satisfies
$(\bJ_{\phi_1, \infty, \leq})$, or equivalently,  \eqref{e:1.12} with $\phi_1$ in place of $\phi$ there.
Then, for every $l\ge2$ and
 $C_* \in (0,1)$,
 there exist $c_1,c_2, c_3>0$ such that for all $x_0\in {\bar D}$, $\lambda>0$, $t>0$ and $x, y\in B_{{\bar D}}(x_0,l\lambda)$ with $ \rho_D(x,y) \leq t/C_*$,
$$
 q^{\langle \lambda \rangle, B(x_0,l\lambda)} (t,x,y)\le
\frac{c_1}{V_D(x_0,\lambda)}\,\left(\left(\frac{\lambda}{\sqrt{t}}\right)^{d_1}\vee \left(\frac{\lambda}{\sqrt{t}}\right)^{d_2}\right)
 \exp \left(-\frac{c_2\rho_D(x,y)^2}t \right);
$$
and  that for all $x_0\in {\bar D}$, $\lambda>0$, $t>0$ and $x, y\in B_{{\bar D}}(x_0,l\lambda)$ with $ \rho_D(x,y) > t/C_*$,
\begin{align*}
 q^{\langle \lambda \rangle, B(x_0,l\lambda)} (t,x,y)&\le
\frac{c_1}{V_D(x_0,\lambda)}\,\left(\left(\frac{\lambda}{\sqrt{t}}\right)^{d_1}\vee \left(\frac{\lambda}{\sqrt{t}}\right)^{d_2}\right) \exp \left(-c_3\rho_D(x,y) \left(\log\frac{\rho_D(x,y)}{t}\right)\right).
\end{align*}
\end{itemize}
\end{proposition}

\pf
Let $l\ge2$, $x_0\in {\bar D}$ and $\lambda>0$.
 For fixed $x, y\in
 B_{{\bar D}}(x_0,l\lambda)$, let
 $R:=\rho_D(x,y)$ and $ f(\xi):=s \left(  \rho_D(\xi,x) \wedge R \right)$ for  $\xi \in {\bar D}$,
where $s>0$ is a constant to be chosen
 later for each case.
Then, $f \in \Lip_{c}({\bar D})$  and
$ |f (\xi)-f (\eta)|\leq s \rho_D(\xi,\eta) $ for every $\xi,\eta \in {\bar D}.$
Thus, by \cite[Proposition 2.21]{GS},
$
 \left\| \Gamma_0 (f , f )  \right\|_\infty \le s^2.
$
This
together with
the elementary inequality that
 $|e^{r}-1|^{2}\leq
r^{2}e^{2|r|}$ for all $r\in \R$ yields that
\begin{equation}\label{e:betainf}\begin{split}
 e^{-2f (\xi)} \Gamma_{\langle \lambda \rangle}(e^f ) (\xi) &=  e^{-2f(\xi)} \Gamma_0 (e^f , e^f ) (\xi)+ \int_{\{ \eta \in { D}: \rho_D(\xi,\eta) \le \lambda \}}(e^{f (\xi)-f (\eta)}-1)^{2} J(\xi, \eta)\,m(d\eta)
\\
&\le   \Gamma_0 (f , f ) (\xi)+ s ^{2}  \int_{\{ \eta \in { D}: \rho_D(\xi,\eta) \le \lambda \}}  \rho_D(\xi,\eta)^{2}e^{2 s \rho_D(\xi,\eta) }  J(\xi, \eta)\,m(d\eta)
\\
&\leq
s^2+
 s ^{2}   \int_{\{ \eta \in { D}: \rho_D(\xi,\eta) \le \lambda \}}  \rho_D(\xi,\eta)^{2}e^{2s\rho_D(\xi,\eta) } J(\xi, \eta)\,m(d\eta) \\
&\leq s^2+ c_0\kappa_1 s ^{2}  \int_{{D}} \rho_D(\xi,\eta)^{2}\frac{e^{2s\rho_D(\xi,\eta)-
 \kappa_2
 \rho_D(\xi,\eta)^{\beta_*}} }{V_D(\xi,\rho_D(\xi,\eta)) \phi_1(\rho_D(\xi,\eta))}\,m(d\eta),
 \end{split}
 \end{equation}
where the last inequality is due to
$(\bJ_{\phi_1, \beta_*, \leq})$.

(i)
Suppose that the jumping kernel $ J(x,y)$ satisfies
 $(\bJ_{\phi_1, \beta_*, \leq})$
   with $\beta_* \in (0,  1]$.
 Take $s=
\kappa_2  \lambda^{\beta_*-1}/4$. We have for $u\in [0, \lambda]$,
 $$2su-  \kappa_2 u^{\beta_*}  =
\kappa_2\lambda^{\beta_*-1} u/2-
\kappa_2 u^{\beta_*}
 =  \kappa_2 u^{\beta_*} \left((u/\lambda)^{1-\beta_*}/2  -1\right) \le -
\kappa_2 u^{\beta_*} /2 .
$$
 Thus by \eqref{e:betainf} and
Lemma \ref{intelem}(ii),
$$
e^{-2f (\xi)} \Gamma_{\langle \lambda \rangle}(e^f ) (\xi) \le C s^2,
\qquad \hbox{where } C:=1+c_1\int_{0}^{\infty}
\frac{u}{ \phi_1(u)} e^{-  \kappa_2 (u/2)^{\beta_*} /2}\,
d u<\infty.
$$
The same estimate holds for
$e^{2f (\xi)} \Gamma_{\langle \lambda \rangle}(e^{-f} ) (\xi)$.
Hence
$\Lambda^{\langle \lambda \rangle} (f)^2 \leq Cs^2$.
Thus
$$
-|f(y)-f(x)|+2
\Lambda^{\langle \lambda \rangle}(f)^2\,
 t \,\le\,
 s  \left( - R+2C\, t s \right)  = \frac{ \kappa_2}4
 \lambda^{\beta_*-1} \left(-
 R +
   \frac{C \kappa_2}{2} \left(\frac{t}{\lambda^{2-\beta_*}}\right) \lambda \right).
$$
In particular, letting $C_* :=
 (C \kappa_2)^{-1}$, it holds that
\begin{equation}\label{e:3.26}
-|f(y)-f(x)|+\;
 2\Lambda^{\langle \lambda \rangle}(f)^2\, t
\leq - \frac{ \kappa_2}{8} \lambda^{\beta_*-1}R \quad \hbox{for  }
0<t\leq C_* \lambda^{1-\beta_*}R.
\end{equation}
On the other hand, if $t > C_* \lambda^{1-\beta_*}R$, then,
with  $\displaystyle s:=\frac{R}{4Ct} =
\frac{ \kappa_2 C_*R}{4t} \le \frac{
\kappa_2}{4} \lambda^{\beta_*-1}$, we have
$$
2su-  \kappa_2 u^{\beta_*}\le
\frac{ \kappa_2}2\lambda^{\beta_*-1} u-
\kappa_2 u^{\beta_*}
=-  \kappa_2 u^{\beta_*} \left(-\frac{1}2(u/\lambda)^{1-\beta_*}  +1\right)
 \le -  \kappa_2 u^{\beta_*} /2$$ for all $ u \le \lambda .$
Hence, for
$t > C_* \lambda^{1-\beta_*}R$,
\begin{equation}\begin{split}
-|f(y)-f(x)|+
2\Lambda^{\langle \lambda \rangle}(f)^2\, t \,
\le s  \left( -R+2C\, t s \right) = \frac{R}{4Ct} \left(-R +\frac{R}2 \right)
=  -\frac{1}{8C} \frac{R^2}{t} .\label{e:3.27}
\end{split}\end{equation}
Now, applying \eqref{e:3.26} and \eqref{e:3.27} to \eqref{e:3.2qw1}, we
get
the desired assertion.

\medskip

(ii)
Suppose next that the jumping kernel $ J(x,y)$ satisfies
$(\bJ_{\phi_1, \beta_*, \leq})$
with $\beta_*\in  (1, \infty)$.
When $0<s \le  \kappa_2/4$, we have  by \eqref{e:betainf}
\begin{align*}
e^{-2f (\xi)} \Gamma_{\langle \lambda \rangle}(e^f ) (\xi)
 \le & s^2+
c_0 \kappa_1 e^{ \kappa_2/2}s ^{2}
 \int_{\{\eta \in {D} :\rho_D(\xi,\eta) \le 1\}} \frac{\rho_D(\xi,\eta)^{2} }{V_D(\xi,\rho_D(\xi,\eta)) \phi_1(\rho_D(\xi,\eta))}\,m(d\eta)\\
 &  +
 c_0\kappa_1 s ^{2}
 \int_{\{\eta \in {D} :\rho_D(\xi,\eta) \ge 1\}} \frac{\rho_D(\xi,\eta)^{2}e^{-
\kappa_2 \rho_D(\xi,\eta)/2   } }{V_D(\xi,\rho_D(\xi,\eta)) \phi_1(\rho_D(\xi,\eta))}\,m(d\eta).
\end{align*}

Thus by (i) and (ii) of  Lemma \ref{intelem}
and \eqref{eqn:cond1},
 $$e^{-2f (\xi)} \Gamma_{\langle \lambda \rangle}(e^f ) (\xi)\leq  s^2+ c_2 s ^{2}  \left(\int_{0}^{1} \frac{u}{ \phi_1(u)}\,du+1
\right) \leq c_3s^2/2 \quad \text{for all }
 0< s \le  \kappa_2/4, $$
which implies that
\begin{equation}
-|f(y)-f(x)|+
2\Lambda^{\langle \lambda \rangle}(f)^2\, t \leq -s R+c_{3}ts^2= s (-R+c_3 t s)
 \quad \text{for all }
 0< s \le  \kappa_2/4.
\label{add1130}\end{equation}
Choose $c_3$ larger if necessary so that
$c_3  \,\ge\, \frac{2}{C_*  \kappa_2}.$
For each
$R \le t/C_*$, take $\displaystyle s:=\frac {R}{2c_3t} \le \frac1 {2c_3C_*} \le \frac{ \kappa_2}{4}$.
Then it follows  from \eqref{add1130} that
$$-|f(y)-f(x)|+
2\Lambda^{\langle \lambda \rangle}(f)^2\, t   \le   -\frac{R^2}{4c_3t}.$$
Putting this into \eqref{e:3.2qw1}, we obtain the assertion
for $R \le t/C_*$.

 \medskip

On the other hand, with $c_4:=(4/ \kappa_2)^{1/(\beta_*-1)}$,  we have by
\eqref{e:betainf}
and (ii) and (iii) of Lemma \ref{intelem}  that for all $s>0$,
\begin{align*}
 e^{-2f (\xi)} \Gamma_{\langle \lambda \rangle}(e^f ) (\xi)
& \le  s^2+ c_0 \kappa_1 s ^{2}
 \int_{\{\eta \in D :\rho_D(\xi,\eta) \le c_4 s^{1/(\beta_*-1)}\}} \frac{\rho_D(\xi,\eta)^{2} e^{2s\rho_D(\xi,\eta)} }{V_D(\xi,\rho_D(\xi,\eta))\phi_1(\rho_D(\xi,\eta))}\,m(d\eta)\\
 &\quad +  c_0\kappa_1 s^{2}
 \int_{\{\eta \in  D :\rho_D(\xi,\eta) \ge c_4s^{1/(\beta_*-1)}\}} \frac{\rho_D(\xi,\eta)^{2}
 e^{- \kappa_2 \rho_D(\xi,\eta)^{\beta_*}/2   } }{V_D(\xi,\rho_D(\xi,\eta)) \phi_1(\rho_D(\xi,\eta))}\,m(d\eta)\\
 & \le  s^2+    c_5 s ^{2}  \int_{0}^{c_4 s^{1/(\beta_*-1)}}
\frac{u\,e^{4s u}}
{ \phi_1(u)}\,du+    c_5\kappa_1 s ^{2} \int_{D}\frac{\rho_D(\xi,\eta)^{2}
e^{- \kappa_2 \rho_D(\xi,\eta)^{\beta_*}/2   } }{V_D(\xi,\rho_D(\xi,\eta)) \phi_1(\rho_D(\xi,\eta))}m(d\eta)\\
&\le c_6\left(s ^{2}
 \exp(4 c_4 s^{{\beta_*}/({\beta_*-1})})
\int_{0}^{c_4 s ^{1/(\beta_*-1)}}\frac{u}{ \phi_1(u)}\,du+
  s^{2}\right).\end{align*}
According to
\eqref{eqn:poly} and  \eqref{eqn:cond1},
the above is less than or equal to
\begin{eqnarray*} c_7 s ^{2} \left(
\frac{s^{2/(\beta_*-1)}}{ \phi_1(s^{1/(\beta_*-1)})}
\,\exp(4c_4 s^{{\beta_*}/({\beta_*-1})})
 + 1 \right)
\leq  2^{-1} c_8\,s^2 \exp(c_8 s^{{\beta_*}/({\beta_*-1})}), \quad s>0.
\end{eqnarray*}
So,
$$
-|f(y)-f(x)|+2
\Lambda^{\langle \lambda \rangle}(f)^2\, t \leq -s R+c_{8}ts^2 \exp(c_8 s^{{\beta_*}/({\beta_*-1})})
= s R\left(c_8 s (t/R)\exp(c_8 s^{{\beta_*}/({\beta_*-1})})
-1\right).
$$
Choose $c_9 \le (2c_8)^{-(\beta_*-1)/\beta_*}$ small so that
$$
c_8 c_9 \sup_{a \ge  C_*^{-1}}\left(\log  a \right)^{(\beta_*-1)/\beta_*} a^{-1/2} <\frac12$$
and take
$
s:=c_9 \left(\log (R/t)\right)^{(\beta_*-1)/\beta_*}.
$
Then,
$$
(t/R)\exp\big(c_8 s^{{\beta_*}/({\beta_*-1})}\big)
=(t/R)(R/t)^{c_8 c_9^{\beta_*/(\beta_*-1)}}\le (t/R)(R/t)^{1/2}= (R/t)^{-1/2},
$$
and so, for $R/t >C_*^{-1} >1$,
\begin{align*}
-|f(y)-f(x)|+2
\Lambda^{\langle \lambda \rangle}(f)^2\, t
&\le s R \left( c_8 c_9 \left(\log (R/t)\right)^{(\beta_*-1)/\beta_*}  (R/t)^{-1/2}  -1 \right)\\
&\le  -\frac12 s R = -\frac12 c_9 R\left(\log (R/t)\right)^{(\beta_*-1)/\beta_*}.\end{align*}
Putting this into \eqref{e:3.2qw1}, we obtain
the desired  assertion for $R> t/C_*$.
\medskip

  (iii)
Now suppose that the jumping kernel $ J(x,y)$ satisfies
$(\bJ_{\phi_1, \infty, \leq})$, or equivalently,  \eqref{e:1.12} with $\phi_1$ in place of $\phi$ there.
In this case, the argument is similar to that of (ii).
Since $ J(\xi, \eta)=0$ on $\rho_D(\xi,\eta) >1$, it follows from
\eqref{e:betainf}  and Lemma \ref{intelem}(iii)  that
\begin{equation}\label{e:3.177}\begin{split}
e^{-2f (\xi)} \Gamma_{\langle \lambda \rangle}
(e^f ) (\xi)
&\leq   s^2+ c_{10} s ^{2}  \int_{0}^{1} \frac{u\,e^{4su }}{ \phi_1(u)}\,du .\end{split}
\end{equation}

\medskip

When $0<s \le 1$, by \eqref{e:3.177},
\begin{equation}
-|f(y)-f(x)|+2
\Lambda^{\langle \lambda \rangle}(f)^2\, t\le  -s R+
c_{11} ts^2= s \left(-R+c_{11} t s \right).
\label{add11311}\end{equation}
Let
$c_{12} = c_{11}  \vee \frac{1}{2C_*}$.
For each
$R \le t/C_*$, take
 $ \displaystyle
s:=\frac {R}{2c_{12}t} \le \frac1 {2 c_{12}C_*} \le 1$.
We have from \eqref{add11311} that
$$-|f(y)-f(x)|+2
\Lambda^{\langle \lambda \rangle}(f)^2\, t
\leq s \left(-R+ c_{12} t s \right)
\le   -\frac{R^2}{4 c_{12}t}.
$$
Putting this into \eqref{e:3.2qw1}, we obtain the desired assertion
for $\rho_D(\xi,\eta) \le t/C_*$.

On the other hand, by \eqref{e:3.177}, for all $\lambda>0$,
\begin{equation}
-|f(y)-f(x)|+2
\Lambda^{\langle \lambda \rangle}(f)^2\, t\leq -sR+
c_{13}ts^2 e^{ 4 s}
= sR\left(c_{13} s (t/R)e^{ 4 s}
-1\right).
\label{add113011}\end{equation}
Choose $c_{14} \in (0,  1/8)$ small so that
$$
c_{13} c_{14} \, \sup_{a \ge  1/C_*} \frac{\log  a}{\sqrt{a}   }   <\frac12 .
$$
For  any $x, y\in  B_{\bar D}(x_0, l\lambda)$ and  $t>0$ with $R/t >C_*^{-1} >1$,
  take $s:=c_{14} \left(\log (R/t)\right)>0$
and we have
$$
(t/R)e^{ 4 s} =(t/R)(R/t)^{ 4 c_{14}}\le (t/R)(R/t)^{1/2}= (R/t)^{-1/2} .
$$
Thus, from \eqref{add113011}, we get
\begin{align*}
-|f(y)-f(x)|+2
\Lambda^{\langle \lambda \rangle}(f)^2\, t
\le  sR \left( c_{13} c_{14} \left(\log (R/t)\right)  (R/t)^{-1/2}  -1 \right)
\le   -\frac12 s R = -\frac12 c_{14} R\log (R/t).\end{align*}
Putting this into \eqref{e:3.2qw1}, we obtain
the desired assertion for $R> t/C_*$ as well.
\qed

The estimates for the following short time region require more sophisticated choices of test functions
in order to obtain the right polynomial exponents.
\begin{proposition}\label{P:2.4c}
Suppose that
$(\bJ_{\phi_1, \beta_*, \leq})$ holds for some $\beta_*\in (0, \infty]$.
Then the following hold.
\begin{itemize}
\item[{\rm(i)}] If $\beta_*\in (0,1]$, then, for every $l\ge2$ and $a>0$, there exist $c_1,c_2>0$ and  $C^* \in (0,1)$
 such that for every
   $x_0\in {\bar D}$,
  $\lambda>a/(4C^*(d_2+2))$,
$x, y\in B_{{\bar D}}(x_0,l\lambda) $ with $\rho_D(x,y)\geq t/C^*$,
$$
 q^{\langle \lambda \rangle, B(x_0,l\lambda)} (t,x,y) \le
\frac{c_1}{V_D(x_0,\lambda)}\,\left(\left(\frac{\lambda}{\sqrt{t}}\right)^{d_1}\vee \left(\frac{\lambda}{\sqrt{t}}\right)^{d_2}\right) \left(\frac{\rho_D(x,y)}{t}\right)^{-\rho_D(x,y)/(8\lambda)} e^{-c_2\lambda^{\beta_*-1 }\rho_D(x,y)}.
$$
\item[{\rm(ii)}]  If $\beta_*\in (1,\infty)$, then, for every $l\ge 2$ and $a>0$, there exist $c_1,c_2>0$ and
 $C^* \in (0,1)$
 such that for every $x_0\in {\bar D}$, $\lambda>a/(4C^*(d_2+2))$,
   $x, y\in B_{{\bar D}}(x_0,l\lambda) $ with $\rho_D(x,y)\geq t/C^*$,
\begin{align*}
 q^{\langle \lambda \rangle, B(x_0,l\lambda)} (t,x,y) \le&
\frac{c_1}{V_D(x_0,\lambda)}\,\left(\left(\frac{\lambda}{\sqrt{t}}\right)^{d_1}\vee \left(\frac{\lambda}{\sqrt{t}}\right)^{d_2}\right) \left(\frac{\rho_D(x,y)}{t}\right)^{-\rho_D(x,y)/(8\lambda)}\\
&\times \exp\left(-c_2\rho_D(x,y) \left( \log (\rho_D(x,y)/t)\right)^{(\beta_*-1)/\beta_*}\right).
\end{align*}

\item[{\rm(iii)}]  Suppose that the jumping kernel $ J(x,y)$ satisfies
$(\bJ_{\phi_1, \infty, \leq})$, or equivalently,  \eqref{e:1.12} with $\phi_1$ in place of $\phi$ there.
Then, for every $l\ge 2$ and $a>0$, there exist $c_1,c_2>0$ and   $C^* \in (0,1)$
 such that for every $x_0\in {\bar D}$, $\lambda>a/(4C^*(d_2+2))$,
  $x, y\in B_{{\bar D}}(x_0,l\lambda) $ with $\rho_D(x,y)\geq t/C^*$,
\begin{align*}
 q^{\langle \lambda \rangle, B(x_0,l\lambda)} (t,x,y) \le&
\frac{c_1}{V_D(x_0,\lambda)}\,\left(\left(\frac{\lambda}{\sqrt{t}}\right)^{d_1}\vee \left(\frac{\lambda}{\sqrt{t}}\right)^{d_2}\right) \left(\frac{\rho_D(x,y)}{t}\right)^{-\rho_D(x,y)/(8\lambda)}\\
&\times \exp\big(-c_2\rho_D(x,y)
  \log (\rho_D(x,y)/t)\big).
\end{align*}
\end{itemize}
\end{proposition}

\pf
Let $l\ge2$, $x_0\in {\bar D}$ and $\lambda>0$. Fix $x, y\in B_{{\bar D}}(x_0,l\lambda)$ with $\rho_D(x,y)\geq a/C^*$,
where $C^*>0$ is a constant to be chosen  later for each cases.
Set $R:=\rho_D(x,y)$ and
$
f(\xi):=\frac {s+v}3
\left( \rho_D(\xi, x) \wedge R \right)$ for all $\xi \in {\bar D},$
where $s$ and $v$ are two positive constants to be chosen later.
Since $|f(\eta)-f(\xi)|\leq\frac
{s+v}3  \rho_D(\xi,\eta)$ for all $\xi,\eta\in {\bar D}$,
 by the same argument as that for \eqref{e:betainf},
\begin{align*}
e^{-2f (\xi)} \Gamma_{\langle \lambda \rangle}(e^f ) (\xi)
&\leq {c_0} \, (s+v)^2 \,\left(1+\int_{\{ \eta \in {D}: \rho_D(\xi,\eta) \le \lambda \}} \rho_D(\xi,\eta)^{2}\frac{e^{2(s+v)\rho_D(\xi,\eta)/3-
\kappa_2 \rho_D(\xi,\eta)^{\beta_*}} }{V_D(\xi,\rho_D(\xi,\eta)) \phi_1(\rho_D(\xi,\eta))}\,m(d\eta)\right)\\
&\leq {c_0} \, (s+v)^2 \,\left(1+e^{2s\lambda/3} \int_{\{ \eta \in {D}: \rho_D(\xi,\eta) \le \lambda \}} \rho_D(\xi,\eta)^{2}\frac{e^{2v\rho_D(\xi,\eta)/3-
\kappa_2 \rho_D(\xi,\eta)^{\beta_*}} }{V_D(\xi,\rho_D(\xi,\eta)) \phi_1(\rho_D(\xi,\eta))}\,m(d\eta)\right).
\end{align*}

(i)
Suppose that the jumping kernel $ J(x,y)$ satisfies
$(\bJ_{\phi_1, \beta_*, \leq})$
 with $\beta_*\in (0,1]$.
Take $v= \kappa_2 \lambda^{\beta_*-1}$.
Since
$2vu/3- \kappa_2 u^{\beta_*} \le
- ( \kappa_2/3) u^{\beta_*}$ for $u \le \lambda$,
we have by Lemma \ref{intelem} (i) and (ii) that, for $\lambda>0$,
 \begin{align*}  &e^{-2f (\xi)} \Gamma_{\langle \lambda \rangle}(e^f ) (\xi)\\
 &\leq  {c_0} \, (s+v)^2 \,\bigg[1+e^{2s\lambda/3}\bigg( \int_{\{ \eta \in {D}: \rho_D(\xi,\eta) \le 1 \}} \frac{\rho_D(\xi,\eta)^{2} }{V_D(\xi, \rho_D(\xi,\eta)) \phi_1(\rho_D(\xi,\eta))}\,m(d\eta)
 \\& \qquad\qquad\qquad\qquad\qquad\qquad+\int_{\{ \eta \in {D}: 1 \le \rho_D(\xi,\eta) <\infty \}} \frac{\rho_D(\xi,\eta)^{2}e^{-(
\kappa_2/3) \rho_D(\xi,\eta)^{\beta_*}} }{V_D(\xi,\rho_D(\xi,\eta)) \phi_1(\rho_D(\xi,\eta))}\,m(d\eta)\bigg)\bigg]\\
  &\leq   \frac{c_1}{6}
 \, \left(s+ \kappa_2\lambda^{\beta_*-1}\right)^2 \, e^{2s\lambda/3}
\end{align*}
for every $\xi\in {\bar D}$.
 The same estimate holds for
$\Gamma_{\langle \lambda \rangle} (-f)(\xi)$.
Thus we have $\Lambda^{\langle \lambda \rangle}(f)^2 \leq \frac{c_1}{6}
   \left(s+ \kappa_2\lambda^{\beta_*-1}\right)^2 \, e^{2s\lambda/3}$ and so
\begin{equation}\label{eqn:22512}
-|f(y)-f(x)|+2 \Lambda^{\langle \lambda \rangle}(f)^2\, t \,\le\,\frac {s+ \kappa_2\lambda^{\beta_*-1}}3
 R \left( -1+ c_1 \,
 \left(s+ \kappa_2\lambda^{\beta_*-1}\right)\, (t/R)e^{2s\lambda/3}
\right).
\end{equation}
Now let $s=\frac {3}{4\lambda}\log (\frac {R}t)$, and then $ \displaystyle
(t/R)e^{2s\lambda/3}= (t/R)e^{\frac 12\log
\frac{R}t}=\sqrt{ {t}/{R} }.
$
Choose
$C^*\in (0, 1)$ such that
$$
\frac{3C^*(d_2+2)}{a}
\left(\sup_{0<v\leq C^*} \sqrt{v} \log \frac1{v} \right) +  \kappa_2
\left(\frac{4C^*(d_2+2)}{a}\right)^{1-\beta_*}
\sqrt{ C^*}
 <(2c_1)^{-1},
$$
 where $c_1>0$ is the constant in (\ref{eqn:22512}).
Then for
 $R\geq t/C^*$ and $\lambda \ge\frac{a}{4C^*(d_2+2)} $,
we have  by \eqref{eqn:22512},
\begin{align*}
& -|f(y)-f(x)|+ 2\Lambda^{\langle \lambda \rangle}(f)^2\, t \\
&\le\frac {s+ \kappa_2\lambda^{\beta_*-1}}3
 R \left( -1+  \,
 c_1 \frac{3}{4\lambda}\sqrt{\frac tR}\log \Big(\frac {R}t\Big) +c_1 \kappa_2\lambda^{\beta_*-1}\sqrt{\frac tR}
\right) \\
&\le\frac{s+ \kappa_2\lambda^{\beta_*-1}}3 \, R\, \left(-1+
 c_1 \frac{3C^*(d_2+2)}{a}
 \,
\sqrt{\frac{t}{R}}
\log\left(\frac{R}t\right)+c_1  \kappa_2
\Big(\frac{4C^*(d_2+2)}{a}\Big)^{1-\beta_*}
\sqrt{\frac tR}\right) \\
&\leq - \frac{s+ \kappa_2\lambda^{\beta_*-1}}{6} \, R
\le
- \frac{R}{8\lambda} \log \Big(\frac{R}t \Big)-
 \frac{ \kappa_2}{6} \lambda^{\beta_*-1} \, R.
\end{align*}
 Applying the estimate above into \eqref{e:3.2qw1}, we get the desired assertion.

\medskip

(ii)
Suppose next that the jumping kernel $ J(x,y)$ satisfies
$(\bJ_{\phi_1, \beta_*, \leq})$
with $\beta_*\in (1,\infty)$.   Let $c_* :=(4/3 \kappa_2)^{1/(\beta_*-1)}$.
  Then, since $2vu/3- \kappa_2 u^{\beta_*} \le - \kappa_2 u^{\beta_*} /2$ for
   $c_*v^{1/(\beta_*-1)} \le u$,
   by Lemma \ref{intelem}   (ii) and (iii),
 for every $\xi \in {\bar D}$,
\begin{align*}
& e^{-2f (\xi)} \Gamma_{\langle \lambda \rangle}(f)(\xi) \\
&\leq  {c_0} \, (s+v)^2 \,\Bigg[1+e^{2s\lambda/3}\bigg( \int_{\{ \eta \in {D}: \rho_D(\xi,\eta) \le c_*v^{1/(\beta_*-1)}  \}} \rho_D(\xi,\eta)^{2}\frac{e^{
 2v\rho_D(\xi,\eta)/3}
  }{V_D(\xi,\rho_D(\xi,\eta)) \phi_1(\rho_D(\xi,\eta))}\,m(d\eta)\\&\qquad \qquad\quad \quad\qquad\qquad\qquad +\int_{\{ \eta \in {D}: c_*v^{1/(\beta_*-1)} \le \rho_D(\xi,\eta) <\infty \}} \rho_D(\xi,\eta)^{2}\frac{e^{2v\rho_D(\xi,\eta)/3
  - \kappa_2 \rho_D(\xi,\eta)^{\beta_*}} }{V_D(\xi, \rho_D(\xi,\eta)) \phi_1(\rho_D(\xi,\eta))}\,m(d\eta)\bigg)\Bigg]\\
 \\
&\leq  {c_0} \, (s+v)^2 \,\Bigg[1+e^{2s\lambda/3}\bigg(
e^{c_* v^{\beta_*/(\beta_*-1)}} \int_{\{ \eta \in { D}: \rho_D(\xi,\eta) \le c_*v^{1/(\beta_*-1)}  \}} \frac{\rho_D(\xi,\eta)^{2}
e^{-v\rho_D(\xi,\eta) /3 }
}{V_D(\xi,\rho_D(\xi,\eta)) \phi_1(\rho_D(\xi,\eta))}m(d\eta)\\&\qquad\qquad\quad\quad\qquad\qquad\qquad+
\int_{\{ \eta \in {D}: c_*v^{1/(\beta_*-1)} \le \rho_D(\xi,\eta) <\infty \}}
\frac{ \rho_D(\xi,\eta)^{2}
e^{- \kappa_2 \rho_D(\xi,\eta)^{\beta_*}/2} }{V_D(\xi,\rho_D(\xi,\eta)) \phi_1(\rho_D(\xi,\eta))}\,m(d\eta)\bigg)\Bigg] \smallskip \\
   &\leq   {c_0} \, (s+v)^2\,\left(1+ c_1 e^{2s\lambda/3} e^{c_* v^{\beta_*/(\beta_*-1)}} \right)  \\
   &\leq   \frac{c_2}6
\, (s+v)^2 \, e^{2s\lambda/3} e^{(4c_3)^{-1} v^{\beta_*/(\beta_*-1)}},
\end{align*}
where $c_3:=1/(4c_*)$.
The same estimate holds for
$\Gamma_{\langle \lambda \rangle} (-f)(\xi)$.
Hence we have $$\Lambda^{\langle \lambda \rangle}(f)^2 \leq \frac{c_2}6
  (s+v)^2 \, e^{2s\lambda/3} e^{(4c_3)^{-1} v^{\beta_*/(\beta_*-1)}}.$$
Take
 $v:=\left(c_3 \log \frac{R}{t}\right)^{(\beta_*-1)/\beta_*}$.
Then for
$t >0$,
\begin{align*}
 -|f(y)-f(x)|+ 2\Lambda^{\langle \lambda \rangle}(f)^2\, t
&\le\,\frac {s+v}3 \left[-R+
c_2\,  \left(s+ \left(c_3\log \frac{R}{t}\right)^{(\beta_*-1)/\beta_*}\right)\, e^{2s\lambda/3} \left(\frac{R}{t}\right)^{1/4} \, t
\right]\\
& =\frac {s+v}3 \left[R\left(-1+
c_2\,  \left(s+ \left(c_3\log \frac{R}{t}\right)^{(\beta_*-1)/\beta_*}\right)\, e^{2s\lambda/3} \left(\frac{t}{R}\right)^{3/4} \right)
\right].
\end{align*}
Next we take $s:=
\frac {3}{4\lambda}\log (\frac {R}t)$.
  Then,
$$e^{2s\lambda/3}\left(\frac{t}{R}\right)^{3/4}= \exp\left(\frac 1{2}\left(\log
\frac{R}t\right)\right)\left(\frac{t}{R}\right)^{3/4}=\left(
\frac{t}R\right)^{1/4}.
$$
 Choose
$C^*\in (0, 1 )$ such
that
$$ c_2 \,  \sup_{0<v<C^*}\left(\frac
{3C^*}{a}(d_2+2)\log \frac1{v}+ \left(c_3\log \frac1{v}\right)^{(\beta_*-1)/\beta_*}\right)\, v^{1/4} <1/2.
$$
Then, for $R\geq t/C^*$ and  $\lambda \ge
\frac{a}{4C^*(d_2+2)}$, it holds that
\begin{align*}
&-|f(y)-f(x)|+2\Lambda^{\langle \lambda \rangle}(f)^2\, t
\\&\le \frac{s+v}3 \left[R\left(-1+c_2 \,  \left(\frac
{3C^*}{a}(d_2+2)\log\left (\frac {R}t\right)+\left(c_3\log \frac{R}{t}\right)^{(\beta_*-1)/\beta_*}\right)\, \left(\frac{t}{R}\right)^{1/4} \right)
\right]
\\
&\le-\frac{s+v}6 R
\,\le\, -\frac{R}{8\lambda } \log \frac{R}t - \frac{1}{6} \left(c_3\log \frac{R}{t}\right)^{(\beta_*-1)/\beta_*}
 \, R.
\end{align*} Applying the estimate above into \eqref{e:3.2qw1}, we get the desired assertion.

\medskip

(iii)
Now suppose  the jumping kernel $ J(x,y)$ satisfies
$(\bJ_{\phi_1, \infty, \leq})$, or equivalently,  \eqref{e:1.12} with $\phi_1$ in place of $\phi$ there.
 In this case,
for every $\xi \in {\bar D}$
 we have by \eqref{eqn:cond1}
 and  Lemma \ref{intelem} (iii),
\begin{align*}
 e^{-2f (\xi)} \Gamma_{\langle \lambda \rangle}(f)(\xi)
&\leq  {c_0} \, (s+v)^2 \,\Bigg[1+e^{2s\lambda/3} \int_{\{ \eta \in {\bar D}: \rho_D(\xi,\eta) \le 1  \}} \rho_D(\xi,\eta)^{2}\frac{e^{2v\rho_D(\xi,\eta)/3}  }{V_D(\xi,\rho_D(\xi,\eta)) \phi_1(\rho_D(\xi,\eta))}\,m(d\eta)\Bigg]\\
&\leq   c_0 \, (s+v)^2 \,\left(1+ e^{2s\lambda/3} e^{2 v/3}
\left(\int_{0}^{1}  \frac{u}{ \phi_1(u)} \,du
+ c_1  \right)\right)  \\
 & \leq   \frac{c_2}6
\, (s+v)^2 \, e^{2s\lambda/3} e^{2 v/3}.
\end{align*} With this estimate at hand,
one can follow the argument in (ii)  with $v=\frac{3}{8}\log \frac{R}{t}$ and $s=\frac{3}{4\lambda}\log \frac{R}{t}$
to
get
the desired assertion.
 \qed

Using Propositions \ref{P:2.41} and \ref{P:2.4c}, we have the following result.

\begin{proposition}\label{C:up}
Suppose that
$(\bJ_{\phi_1, \beta_*, \leq })$ holds for some $\beta_*\in (0, \infty]$.
Then the following hold.   \begin{itemize}
\item[{\rm(i)}] If $\beta_*\in (0,1]$, then, for any $a>0$ and any
 $l\ge2$
there are constants
$c_1,c_2, c_3, c_4, C_*, C^*>0$
such that for all $\lambda^2\ge t>a$, $x_0\in {\bar D}$ and $x\in B_{{\bar D}}(x_0,l\lambda)$,
\begin{equation}\label{eqn:kern41}
\int_{\{y\in {\bar D}: \rho_D(x,y)\ge  \lambda \}} q^{\langle \lambda \rangle, B(x_0,l\lambda)} (t,x,y)\,m(dy)\le
   \begin{cases}
c_1e^{-c_2\lambda^2/t}& \text { if } C_*\lambda^{2-\beta_*}\le t\le \lambda^2,
\\
c_1 e^{-c_2\lambda^{\beta_*}}& \text { if } C_*\lambda^{2-\beta_*} \ge t;
\end{cases}
\end{equation} and that for all
$0<t\le a$,
 $4(d_2+2)
\le l_0\le l-1$,
$\lambda\ge a/(4C^*(d_2+2))$, $x_0\in {\bar D}$ and $x\in B_{{\bar D}}(x_0,l\lambda)$,
\begin{equation}\label{eqn:kern40}\int_{\{y\in {\bar D}: \rho_D(x,y)\ge l_0\lambda\}} q^{\langle \lambda \rangle, B(x_0,l\lambda)} (t,x,y)\,m(dy)\le c_3 \left(\frac{\lambda}{\sqrt{t}}\right)^{-l_0/8} e^{-c_4\lambda^{\beta_*}}.\end{equation}

\item[{\rm (ii)}] If $\beta_*\in  (1, \infty]$, then, for every
 $C_* \in (0,1)$, $a>0$ and every  $l\ge 2$,
 there exist  $c_1,c_2, c_3, c_4, C_*, C^*>0$
 such that for all $x_0\in {\bar D}$,  $t>0$,  $ 0<\lambda \leq t/C_*$ and $x\in B_{{\bar D}}(x_0,l\lambda)$,
$$
\int_{\{y\in {\bar D}: \rho_D(x,y)\ge  \lambda\}} q^{\langle \lambda \rangle, B(x_0,l\lambda)} (t,x,y)\,m(dy)\le
c_1
 \exp \left(-\frac{c_2\lambda^2}t \right),
$$
and  for all $x_0\in {\bar D}$, $t>0$, $\lambda > t/C_*$ and $x\in B_{{\bar D}}(x_0,l\lambda)$,
$$
\int_{\{y\in {\bar D}: \rho_D(x,y)\ge  \lambda\}} q^{\langle \lambda \rangle, B(x_0,l\lambda)} (t,x,y)\,m(dy)\le
c_1
 \exp \left(-c_2\lambda \left(1+\log^+\frac{\lambda}{t}\right)^{\frac{\beta_*-1}{\beta_*}} \right);
$$and that for all $0<t\le a$,
$4(d_2+2)\le l_0\le l-1$, $\lambda\ge a/(4C^*(d_2+2))$, $x_0\in {\bar D}$ and $x\in B_{{\bar D}}(x_0,l\lambda)$,
$$\int_{\{y\in {\bar D}: \rho_D(x,y)\ge l_0\lambda\}} q^{\langle \lambda \rangle, B(x_0,l\lambda)} (t,x,y)\,m(dy)\le
c_3 \left(\frac{\lambda}{\sqrt{t}}\right)^{-l_0/8} e^{-c_4\lambda(\log \lambda/t)^{({\beta_*-1})/\beta_*}}.
$$
\end{itemize}
 \end{proposition}

\pf    (i) We first prove \eqref{eqn:kern41}. Let $x_0\in {\bar D}$ and $l\ge2$. Let $C_*$ be the constant in Proposition  \ref{P:2.41}(i). When  $ \lambda^2  \wedge (C_*\lambda^{2-\beta_*}) \ge t\ge a$, according to
Proposition \ref{P:2.41}(i) and
 \eqref{univd},
we have that for  $x\in B_{{\bar D}}(x_0,l\lambda)$,
\begin{align*}\int_{\{y\in {\bar D}: \rho_D(x,y)\ge \lambda\}} q^{\langle \lambda \rangle, B(x_0,l\lambda)} (t,x,y)\,m(dy)&\le \frac{c_1}{V_D(x_0,\lambda)} \left(\frac{\lambda}{\sqrt{t}}\right)^{d_2}\int_{\{y\in {\bar D}: \rho_D(x,y)\ge  \lambda\}} e^{-c_2\lambda^{\beta_*-1}\rho_D(x,y)}\,m(dy)\\
&\le \frac{c_3}{V_D(x,\lambda)} \left(\frac{\lambda}{\sqrt{t}}\right)^{d_2}\sum_{i=0}^\infty e^{-c_2\lambda^{\beta_*-1}( 2^{i}\lambda)}V_D(x, 2^{i+1}\lambda)\\
&\le c_4  \left(\frac{\lambda}{\sqrt{a}}\right)^{d_2}
\sum_{i=0}^\infty  e^{-c_22^{i}\lambda^{\beta_*}}2^{id_2}\le c_5e^{-c_6\lambda^{\beta_*}},
 \end{align*} where $c_5,c_6$ depend on $a$.  On the other hand, when $\lambda^2\ge t\ge C_*\lambda^{2-\beta_*}$, also by Proposition \ref{P:2.41}(i),
 \begin{align*}
&\int_{\{y\in {\bar D}: \rho_D(x ,y)\ge \lambda\}} q^{\langle \lambda \rangle, B(x_0 ,l\lambda)} (t,x,y)\,m(dy)\\
&\le  \frac{c_7}{V_D(x_0,\lambda)} \left(\frac{\lambda}{\sqrt{t}}\right)^{d_2}\int_{\{y\in \bar D: \lambda\le \rho_D(x ,y)\le  C_*^{-1} t \lambda^{\beta_*-1}\}}e^{-c_8\rho_D(x,y)^2/t}\,m(dy)\\
&\quad+ \frac{c_7}{V_D(x_0,\lambda)} \left(\frac{\lambda}{\sqrt{t}}\right)^{d_2}\int_{\{y\in \bar D:  \rho_D(x ,y)\ge C_*^{-1} t\lambda^{\beta_*-1}\}}e^{-c_8\lambda^{\beta_*-1}\rho_D(x,y)}\,m(dy)\\
&\le \frac{c_7}{V_D(x,\lambda)}\left(\frac{\lambda}{\sqrt{t}}\right)^{d_2}
\left( \sum_{i=0}^\infty e^{-c_8 (2^{i}\lambda)^2/t}V_D(x, 2^{i+1}\lambda)+ \sum_{i=0}^\infty e^{-c_8 2^{i}C_*^{-1} \lambda^{2\beta_*-2}t}
V_D(x, 2^{i+1}C_*^{-1} t\lambda^{\beta_*-1}) \right)
\\
&\le c_9 \left(\frac{\lambda}{\sqrt{t}}\right)^{d_2}
\left( \sum_{i=0}^\infty e^{-c_8 (2^{i}\lambda)^2/t}
 2^{i d_2}+ \sum_{i=0}^\infty e^{-c_8 2^{i}C_*^{-1} \lambda^{2\beta_*-2}t}
  2^{id_2} (t\lambda^{\beta_*-2})^{d_2}
  \right)
\\
&= c_9 \
\left(\frac{\lambda}{\sqrt{t}}\right)^{d_2}\sum_{i=0}^\infty e^{-c_8 2^{2i}\lambda^2/t}
 2^{i d_2}+ c_9\left(\lambda^{2\beta_*-2}t \right)^{d_2/2}\sum_{i=0}^\infty e^{-c_8 2^{i}C_*^{-1} \lambda^{2\beta_*-2}t}
  2^{id_2}
\\
&\le  c_{10}\left(e^{-c_{11}\lambda^2/t}+e^{-c_{12}\lambda^{2\beta_*-2}t}\right)\le c_{13} e^{-c_{14}\lambda^2/t}. \end{align*}

Next, we consider \eqref{eqn:kern40}. According to Proposition \ref{P:2.4c}(i), for every
 $l>4(d_2+2)$ and $a>0$, there exists $C^* \in (0,1)$
 such that for every $x_0\in {\bar D}$,   $4(d_2+2)\le l_0\le l-1$,
 $\lambda>a/(4C^*(d_2+2))$, $0<  t \le a$ and
 $x\in B_{{\bar D}}(x_0,l\lambda) $,
\begin{align}
&\int_{\{y\in {\bar D}: \rho_D(x,y)\ge l_0\lambda\}} q^{\langle \lambda \rangle, B(x_0,l\lambda)} (t,x,y)\,m(dy)\nonumber\\
 & \le
\frac{c_1}{V_D(x_0,\lambda)}\left(\frac{\lambda}{\sqrt{t}}\right)^{d_2}\int_{\{y\in {\bar D}: \rho_D(x,y)\ge l_0\lambda\}}\left(\frac{\rho_D(x,y)}{t}\right)^{-\rho_D(x,y)/(8\lambda)} e^{-c_2\rho_D(x,y)\lambda ^{\beta_*-1}}\,m(dy)\nonumber
\\
&\le \frac{c_3}{V_D(x,\lambda)}\left(\frac{\lambda}{\sqrt{t}}\right)^{d_2} \left(\frac{l_0\lambda}{ {t}}\right)^{-l_0/8}\int_{\{y\in {\bar D}: \rho_D(x,y)\ge l_0\lambda\}} e^{-c_2
\rho_D(x,y)\lambda ^{\beta_*-1}}\,m(dy)\nonumber\\
&= \frac{c_3l_0^{-l_0/8}\lambda ^{d_2/2}}{V_D(x,\lambda)}\left(\frac{\lambda}{ {t}}\right)^{-l_0/8+d_2/2}\int_{\{y\in {\bar D}: \rho_D(x,y)\ge l_0\lambda\}} e^{-c_2\rho_D(x,y)\lambda ^{\beta_*-1}}\,m(dy).\label{eq:R2-10-28}
\end{align}
Noting $\lambda/t\ge 1/(  4C^*(2+d_2))$ and $l\ge l_0+1$, we have
\begin{align*}
\eqref{eq:R2-10-28} \le& \frac{c_3l_0^{-l_0/8}\lambda ^{d_2/2}}{V_D(x,\lambda)}\left(\frac{1}{ 4C^*(2+d_2)}\right)^{-l_0/16+d_2/2}\left(\frac{\lambda}{ {t}}\right)^{-l_0/16}\int_{\{y\in {\bar D}: \rho_D(x,y)\ge l_0\lambda\}} e^{-c_2\rho_D(x,y)\lambda ^{\beta_*-1}}\,m(dy)\\
\le &c_4 \lambda ^{d_2/2+l_0/16} \left(\frac{\lambda}{\sqrt{t}}\right)^{-l_0/8}\frac{1}{V_D(x,\lambda)}\int_{\{y\in {\bar D}: \rho_D(x,y)\ge l_0\lambda\}} e^{-c_2\rho_D(x,y)\lambda ^{\beta_*-1}}\,m(dy)\\
\le &c_5 \lambda ^{d_2/2+l_0/16} \left(\frac{\lambda}{\sqrt{t}}\right)^{-l_0/8} e^{-c_6 \lambda^{\beta_*}}
\le c_7 \left(\frac{\lambda}{\sqrt{t}}\right)^{-l_0/8} e^{-c_8 \lambda^{\beta_*}}.
\end{align*}
Here, $c_4$ in the second inequality depends on $l_0$, and the third inequality follows from
 \eqref{univd}
 by splitting the integral into pieces.

\medskip

(ii) For $\beta_*\in (1,\infty)$, the proof is similar to that for (i) by applying Proposition \ref{P:2.41}(ii) (resp.\ Proposition \ref{P:2.4c}(ii)) instead of Proposition \ref{P:2.41}(i) (resp.\ Proposition \ref{P:2.4c}(i)). The details are omitted here.
From Proposition  \ref{P:2.41}(iii) and Proposition \ref{P:2.4c}(iii), we can also establish the desired assertion for the case $\beta_*=\infty$.
 \qed

 Recall that $ \tau_B^{(\lambda)}$ is the first exit time from the ball $B$
 by
 the process $ Y^{(\lambda)}$. Let $ Y^{(\lambda),U}$ be the killed process of $ Y^{(\lambda)}$ when it exits from the open set $U$, and $ \tau_B^{(\lambda),U}$ be the first exit time from the ball $B$ of the process $ Y^{(\lambda),U}$.
 We can check that  for any ball $B:=B(x,r)$ with radius $r$ and $x\in {\bar D}$, $\lambda>0$ and  $t>0$,
 $$
   \bP_x \big( \tau_{B(x,r)}^{(\lambda)}\le t \big)
 =\bP_x\big(\tau_{B(x,r)}^{(\lambda),B(x,\lambda+r)}\le t\big).
 $$
 Furthermore, by
the strong Markov property, for any
$x\in {\bar D}$ and $r>0$,
\begin{align} \label{e:key}
&  \int_{\{y\in {\bar D}:\rho_D(y,x)\ge r\}}  q^{\langle \lambda \rangle}(t,x,y)\,m(dy)  \nonumber \\
&= \bP_x \big(\rho_D( Y_t^{(\lambda)},x)\ge r \big)\le
   \bP_x \big( \tau_{B(x,r)}^{(\lambda)}\le t \big)
   = \bP_x \big( \tau_{B(x,r)}^{(\lambda),B(x,r+\lambda)}\le t \big) = \bP_x\big(\sup_{0<s\le t}\rho_D \big( Y_s^{(\lambda),B(x,r+\lambda)},x \big)\ge r\big)  \nonumber  \\
 &\le \bP_x\big(\rho_D( Y_{2t}^{(\lambda),B(x,r+\lambda)},x)\ge r/2\big)
+ \bP_x\big(\sup_{0<s\le t}\rho_D( Y_s^{(\lambda),B(x,r+\lambda)},x)\ge r,
\, \rho_D( Y_{2t}^{(\lambda),B(x,r+\lambda)},x)\le r/2\big)
\nonumber  \\
&\le \bP_x\big(\rho_D( Y_{2t}^{(\lambda),B(x,r+\lambda)},x)\ge r/2\big)
+ \bP_x\Big( \tau_{B(x,r)}^{(\lambda)}\le t ,
\, \rho_D( Y_{2t}^{(\lambda),B(x,r+\lambda)},Y^{(\lambda)}_{\tau_{B(x,r)}^{(\lambda)}})\ge r/2\Big)
\nonumber  \\
&\le  \bP_x\big(\rho_D( Y_{2t}^{(\lambda),B(x,r+\lambda)},x)\ge r/2\big)
+ \sup_{z\in B(x,r+\lambda)}\sup_{t\le s\le 2t}\bP_z\big(\rho_D( Y_s^{(\lambda),B(x,r+\lambda)},z)\ge r/2\big) \nonumber  \\
&\le   2\sup_{z\in B(x,r+\lambda)}\sup_{t\le s\le 2t}\bP_z\big(\rho_D( Y_s^{(\lambda),B(x,r+\lambda)},z)\ge r/2\big).
\end{align}

\medskip

The following is the main result of this section.
  Recall that $H_{\phi, \beta}(t, x, r)$ is defined  right
 after \eqref{eqn:4}.

\begin{theorem}\label{T1u}
Suppose that $(\bJ_{\phi_1, \beta_*, \leq })$ holds for some $\beta_*\in (0, \infty]$.  Then there
are  positive constants $c_1$ and $c_2$ that
depend only on
characteristic constants $(C_1, C_2, C_3, C_4)$ of $D$
and the constant parameters in
$(\bJ_{\phi_1, \beta_*, \leq})$ and \eqref{eqn:poly} for $\phi_1$
so that
$$
  q(t, x, y) \leq c_1
H_{\phi_1, \beta_*}
 (t, x, c_2 \rho_D(x, y))
\quad \hbox{for every } t>0 \hbox{ and } x, y \in \bar D.
$$
 \end{theorem}

\pf
We will  give the proof for the case when $\beta_*\in (0, 1]$. The other cases can be proved similarly.
Let $C>1$.  When  $t\in (0,1]$ and $\rho_D(x,y)\le C$, the desired estimate follows from
\eqref{e:hk00} with $\phi_*=\phi_1$,
since there is a constant $c_0:=c_0(C)$ such that $c_0   p^{(j)}(t, x,\rho_D(x,y))\ge
  p^{(j)}_{\phi_1}(t,x,\rho_D(x,y))$
for all $t>0$ and $x,y\in \bar D$ with $\rho_D(x,y)\le C$.
Next, set $k=32(2+d_2)$ and $C:= k/(4C^*(2+d_2))=8/C^*$, where $C^* >0$ is the constant in
 Proposition \ref{C:up}.
For $t\in (0,1]$ and $x, y\in \bar D$ with $\rho_D(x,y)\ge C$,  let
   $\lambda=\rho_D(x,y)/k$. According to \eqref{e:up-01}, it suffices to verify that for any $t\in (0,1]$ and $\rho_D(x,y)\ge C$.
\begin{equation}\label{e:pro1}  q^{\langle \lambda\rangle}(t,x,y)\le \frac{c_1t}{V_D(x,\rho_D(x,y))} e^{-c_2\rho_D(x,y)^{\beta_*} }.\end{equation} Indeed, letting $R=\rho_D(x,y)>0$,
by \cite[Lemma 7.2(2)]{CKW1},
it holds that for any $t\in (0,1]$ and $\rho_D(x,y)\ge C$,
$$ q^{\langle \lambda\rangle}(t,x,y)\le  q (t,x,y) e^{t/c(\lambda)}\le c_1q(t,x,y)\le \frac{c_2}{V_D(x,\sqrt{t})},$$
 where $c(\lambda)$ is a positive increasing function of  $\lambda$.
Here the second inequality is due to the facts that
$\lambda=\rho_D(x,y)/k\geq C/k=1/(4C^*(2+d_2))$
and  $t\in (0,1]$, and the third inequality is by \eqref{e:hk00}. Using this, we have
for any
$t>0$ and $x, y\in \bar D $ with $\rho_D(x,y)\ge C$,
 \begin{align*}
   q^{\langle \lambda\rangle}(t,x,y)&=\int_{{\bar D}}  q^{\langle \lambda\rangle}(t/2,x,z)  q^{\langle \lambda\rangle}(t/2,z,y)\,m(dy)\\
&\le \left(\int_{B_{{\bar D}}(x,R/2)^{c}}+\int_{B_{{\bar D}}(y,R/2)^{c}}\right)  q^{\langle \lambda\rangle}(t/2,x,z)  q^{\langle \lambda\rangle}(t/2,z,y)\,m(dy)\\
&\le \frac{c_2}{V_D(y,\sqrt{t})}\int_{B_{{\bar D}}(x,R/2)^{c}} q^{\langle \lambda\rangle}(t/2,x,z)\,m(dz)+\frac{c_2}{V_D(x,\sqrt{t})}\int_{B_{{\bar D}}(y,R/2)^{(c)}} q^{\langle \lambda\rangle}(t/2,y,z)\,m(dz)\\
&\le \frac{c_3}{V_D(x,\rho_D(x,y))}\left(\frac{\rho_D(x,y)}{\sqrt{t}}\right)^{d_2} \cdot \sup_{x\in {\bar D}} \sup_{z\in B(x,R/2+\lambda)}\sup_{t\le s\le 2t}\bP_z  \big( \rho_D(  Y_s^{(\lambda),B(x,R/2+\lambda)},z)\ge R/4 \big) \\
&\le \frac{c_3}{V_D(x,\rho_D(x,y))}\left(\frac{\rho_D(x,y)}{\sqrt{t}}\right)^{d_2} \left(\frac{\lambda}{\sqrt{t}}\right)^{-k/32} e^{-c_4 (R/k)^{\beta_*}}\\
&\le  \frac{c_5}{V_D(x,\rho_D(x,y))}
\frac t{\rho_D(x,y)^2} e^{-c_6 \rho_D(x,y)^{\beta_*}},
\end{align*}
where the third inequality follows from \eqref{e:1.8} and  \eqref{e:key},
and the fourth inequality follows from  \eqref{eqn:kern40} with $l=1+(k/2)$
and $l_0=k/4$.
Thus we have
$$   q^{\langle \lambda\rangle}(t,x,y)\le
\frac{c_7t}{V_D(x,\rho_D(x,y))}e^{-c_8 \rho_D(x,y)^{\beta_*}},$$
which proves \eqref{e:pro1}.

Finally,  suppose $t>1$. For any $c>0$, if $\rho_D(x,y)\le c \sqrt{t}$, then the estimate follows from
 \eqref{e:hk00}; if $\rho_D(x,y)\ge c\sqrt{t}$, we can follow the argument above and use \eqref{eqn:kern41} in place of \eqref{eqn:kern40}
  to get the desired result.
\qed

The upper bounds of $q(t,x,y)$ in Theorem \ref{T1}
follow directly from Theorems \ref{T:2.4} and \ref{T1u}
 for the case that $\beta_*=0_+$ and $\beta_*\in (0,\infty]$, respectively.

\section{Lower bound estimates  under $(\bJ_{\phi_1, 0_+, \leq})$ and $(\bJ_{\phi_2, \beta^*, \geq})$
with $\beta^*\in (0,\infty]$} \label{S:5}

  Throughout this section,   we assume that {\rm diam}$(D)=\infty$ and that $({\bf J}_{\phi_1, 0_+, \leq})$ holds.
 Recall from Subsections \ref{S:2.1} and \ref{S:3.1} that  $Y$ is the Feller process associated with the regular Dirichlet form $(\sE,\sF)$
 of \eqref{e:sF}--\eqref{e:sE} on $L^2(D; m)$,
 which has a jointly H\"older continuous transition density function $  q(t, x, y)$ that enjoys the properties
 stated in Theorems \ref{T:2.4} and Proposition \ref{P:3.1}.

In this section, we will assume the jumping kernel $J$ satisfies condition $(\bJ_{\phi_2, \beta^*, \geq })$ for some $\beta^*\in (0, \infty]$,
where the strictly increasing function $\phi_2$ can  be different from $\phi_1$ in the upper condition $({\bf J}_{\phi_1, 0_+, \leq})$ for $J(x,y)$.

\begin{proposition}\label{main50}
Assume that $({\bf J}_{\phi_1, 0_+, \leq})$  and $(\bJ_{\phi_2, \beta^*, \geq })$ hold for some $\beta^*\in (0, \infty]$.
 Then there are positive constants $c_k$, $k=1,\cdots, 4$, and $c_6$, depending only on
 the characteristic constants $(C_1, C_2, C_3, C_4)$ of $D$ and
constant parameters in  $({\bf J}_{\phi_1, 0_+, \leq})$  and $(\bJ_{\phi_2, \beta^*, \geq })$
as well as  in \eqref{eqn:poly} for $\phi_1$ and $\phi_2$,
and   on the constant  $C_*$ below $($for $c_1, \dots, c_4$$)$ and on both $C_*$ and $c_5$ $($for $c_6)$, so that the following hold.
\begin{enumerate}
\item [\rm (i)] Suppose that $\beta^* \in (0, \infty]$.
Then for every $C_*>0$, there exist $c_1,c_2>0$
  such that for every $t >0$ and $x,y\in \bar D$ with $\rho_D(x,y)\ge  C_* t^{1/2}$,
$$
  q(t, x,y)\geq \frac{c_1}{V_D(x, \sqrt{t})}\exp \left( -\frac{c_2 \rho_D(x,y)^2}t \right).
$$

\item[\rm (ii)] Suppose that  $\beta^* \in (0, \infty)$.
Then for every $C_*>0$, there exist $c_3, c_4>0$ such that for every $t>0$ and $x,y\in \bar D$ with $\rho_D(x,y) \ge
C_*t^{1/2}$,
\begin{equation}\label{ew218}
  q(t, x,y)\geq \frac{c_3
t }{V_D(x, \rho_D(x,y)) \phi_2(\rho_D(x,y))}e^{-c_4 \rho_D(x,y)^{\beta^*}}.
\end{equation}

\item[\rm (iii)] Suppose that the jumping kernel $  J(x,y)$ satisfies
$(\bJ_{\phi_2, \infty, \geq })$,  or equivalently,  \eqref{e:1.13} with $\phi_2$ in place of $\phi$ there.
 Then for any $c_5\in (0,1)$ and $C_*>0$, there is a constant $c_6>0$
 such that for every $t >0$ and $x,y\in \bar D$
$$
  q(t, x,y)\geq
  \frac{c_6t}{V_D(x,\rho_D(x,y))\phi_2 (\rho_D(x,y))}
  \quad \hbox{whenever }   C_*t^{1/2} \leq \rho_D(x,y) \leq c_5.
$$
\end{enumerate}
\end{proposition}

\pf (i)
We mainly follow the proof of \cite[Proposition 5.5]{CKW4}.
Without loss of generality, we may and do
 assume that  $C_*$ is less than  the constant $c_2$ in Proposition \ref{P:3.1}.
Let $r=\rho_D(x,y)$.
Let $N\ge
2$ be an integer such that
\[
9r^2/(C_*^2t) \le N < 9r^2/(C_*^2t)+1 \quad \text{ so that } C_* \sqrt{2t/N} \le 3r/N\le C_* \sqrt{t/N}.
\]
Let $x_0=x, x_1,\cdots, x_{N}=y$ such that
$\rho_D(x_i,x_{i+1})\le
r/N$,
and set
$B_i=B_{\bar D} (x_i, r/N)$ for all  $1\le i\le N$. Then, by Proposition \ref{P:3.1}
 and \eqref{univd},
\begin{align*}
  q(t,x,y) &\ge \int_{
 B_1}  q(t/N,x,y_1)  q(t/N,y_1,y_2)\,m(dy_1)\cdots \int_{B_{N-1}}  q(t/N,y_{N-1},y)\,m(dy_{N-1})\\
& \ge  \frac{c_1}{V_D\big(x,\sqrt{t/N}\big)} \int_{
 B_1}\frac{c_1}{V_D\big(y_1,\sqrt{t/N}\big)}  \,m(dy_1)\cdots \int_{B_{N-1}} \frac{c_1}{V_D\big(y_{N-1},\sqrt{t/N}\big)}\,m(dy_{N-1}) \\
& \ge  \frac{c_1}{V_D\big(x,\sqrt{t/N}\big)} \prod_{i=1}^{N-1}\frac{c_1V_D\big(x_i, C_* \sqrt{2t/N}/3)}{V_D\big(x_i, r/N+  \sqrt{t/N}\big)}
\ge   \frac{c_1}{V_D\big(x,\sqrt{t}\big)} c_2^{N-1} \ge \frac{c_1}{V_D\big(x,\sqrt{t}\big)} e^{-9C_*^{-2}\ln (1/c_2)r^2/t},
\end{align*}
which yields the first assertion.

(ii) and (iii) follow from the L\'evy system
(see \eqref{levy-1}) of the process $Y$,
condition $(\bJ_{\phi_2, \beta^*, \geq })$,
and Proposition \ref{P:3.1}; see the proof of \cite[Proposition 5.4]{CKW1}.
We remark that the restriction  $\rho_D(x,y)\leq c_5$ with $c_5\in (0,1)$ is imposed  in case (iii),
where $  J(x,y)=0$ for all $\rho_D(x,y)> 1$,
  in order to use the L\'evy system argument.
\qed

Recall that $p^{(c)}(t,x,r)$, $ p^{(j)}_{ \phi} (t,x, r)$ and $  p^{(j)}_{\phi, \beta} (t,x, r)
$ are defined by \eqref{e:pc}, \eqref{eqn:5} and \eqref{eqn:4}, respectively.
Applying Propositions \ref{P:3.1} and \ref{main50},  and following  the proof of
\cite[Proposition 5.5]{CKW4},
we can
obtain

\begin{corollary}\label{main5111}
Suppose that
$({\bf J}_{\phi_1, 0_+, \leq})$
 and $(\bJ_{\phi_2, \beta^*, \geq })$ hold with $\beta^*\in (0, \infty]$. Then there
are positive constants $c_i$,
$1\leq i\leq 8$, that
depend on
 the characteristic constants $(C_1, C_2, C_3, C_4)$ of $D$
 and the constant parameters in
 $({\bf J}_{\phi_1, 0_+, \leq})$
  and $(\bJ_{\phi_2, \beta^*, \geq })$  as well as in \eqref{eqn:poly} for $\phi_1$ and $\phi_2$,
so that the following hold.

\begin{enumerate}

\item [\rm (i)]
When $\beta^* \in (0, \infty)$,
$$
  q(t, x,y)\geq  \frac{c_1}{V_D(x,\sqrt{t})} \wedge \left( p^{(c)}(t, x, c_2\rho_D(x,y))+
 p^{(j)}_{\phi_2, \beta^*} (t,x, c_3\rho_D(x,y))
 \right)
\quad \hbox{for } t>0 \hbox{ and } x, y \in \bar D.
$$

 \item [\rm (ii)]
When $\beta^*=\infty$,
 $$
   q(t,x,y)\ge  \frac{c_4 }{V_D(x,\sqrt{t})} \wedge \left( p^{(c)}(t, x,c_5\rho_D(x,y))+p_{\phi_2}^{(j)}(t, x,c_6\rho_D(x,y)) \right)
$$
holds for every
$t \in (0,c_7]$ and $x,y\in \bar D$ with $\rho_D(x,y)\le c_8$.
\end{enumerate}
\end{corollary}

 To  further refine
  lower bound estimates  of the heat kernel  $q(t, x, y)$
  for the case that $\beta^*\in (1,\infty]$, we need the following estimates.

\begin{proposition}\label{main33}
Assume that $({\bf J}_{\phi_1, 0_+, \leq})$  and $(\bJ_{\phi_2, \beta^*, \geq })$ hold for some $\beta^*\in (0, \infty]$.
Then  for any $C^*\in (0, 1)$ and $c\in (0, 1)$,
there are positive constants $c_i$, $1\leq i\leq 4$, that
depend only on
the characteristic constants $(C_1, C_2, C_3, C_4)$ of $D$,
 the constant parameters in $({\bf J}_{\phi_1, 0_+, \leq})$  and $(\bJ_{\phi_2, \beta^*, \geq })$ and in \eqref{eqn:poly} for $\phi_1$ and $\phi_2$
as well as on  the constants $C^*$ $($for $c_1$ and $c_2)$ and on both  $C^*$ and $c$ $($for $c_3$ and $c_4$),
so that the following hold.
\begin{itemize}
\item[{\rm(i)}]
When  $\beta^* \in (1, \infty)$,   for all $t>0$ and $x,y\in \bar D$ with
$\rho_D(x,y)\geq t/C^*$,
\begin{equation}\label{ew2181}
  q(t, x,y)\geq \frac{c_1 t}{V_D(x,\rho_D(x,y)) } \exp\left(-c_2\rho_D(x,y) \left(\log \frac{\rho_D(x,y)}{t} \right)^{({\beta^*}-1)/{\beta^*}}\right).
\end{equation}

\item[{\rm (ii)}]  When $\beta^*=\infty$,  for all  $t>0$ and $x,y\in \bar D$ with $\rho_D(x,y) \ge c\vee ( t/C^*)$,
$$
  q(t, x,y)\geq  \frac{c_3 t}{V_D(x,\rho_D(x,y)) } \left( \frac t{\rho_D(x,y)} \right)^{c_4 \rho_D(x,y)}.
$$
\end{itemize}
 \end{proposition}

\proof (i) Let $r:=\rho_D(x,y)$.
 Note that $\exp\big(-cr(\log (r/t))^{(\beta^*-1)/\beta^*}\big)\ge \exp(-cr^{\beta^*})$ is
equivalent to $t\ge r\exp(-r^{\beta^*})$. Then, according to Propositions \ref{P:3.1} and \ref{main50}, it  suffices to consider the case $C^*r\ge t\ge r\exp(-r^{\beta^*})$. In this case we have $r(\log r/t)^{-1/\beta^*}\ge 1$.
Let $l \ge
2$ be an integer such that
\[
   r\big(\log (r/t) \big)^{-1/\beta^*} < l \le
   r\big(\log (r/t)\big)^{-1/\beta^*}+1 \le 2r\big(\log (r/t)\big)^{-1/\beta^*},
\]
and let $x=x_0,x_1,\cdots,x_l=y\in {\bar D}$ be such that $r/(2l)\le \rho_D(x_i, x_{i+1})\le 2r/l$
for $i= 0,\cdots, l-1$.
We observe that
$$
\frac{t}{2r} (\log (r/t))^{1/\beta^*}\le \frac{t}{l} \le \frac{t}{r} (\log (r/t))^{1/\beta^*}
\le \sup_{s \ge 1/C^*} s^{-1} (\log s)^{1/\beta^*} =: t_0 <\infty$$
and
$$
\frac{r}{2l} \ge
\frac14  (-\log {C^*})^{1/\beta^*} =: r_0 >0.$$
Thus for all $(y_i,y_{i+1}) \in B_{{\bar D}}(x_i, r/(8l)) \times B_{{\bar D}}(x_{i+1},r/(8l))$,
$3r/l\ge \rho_D(y_i, y_{i+1}) \ge r/(4l)$ and so
$\rho_D(y_i, y_{i+1})  \ge r/(4l)\ge   r_0/2\ge (r_0/(2t_0^{1/2}))\cdot ({t}/{l})^{1/2}.
$
According to Corollary \ref{main5111}(i),
for all $(y_i,y_{i+1}) \in B_{{\bar D}}(x_i,r/(8l)) \times B_{{\bar D}}(x_{i+1},r/(8l))$,
\begin{align} \label{chalb}
    q(t/l, y_i,y_{i+1})
  &\geq    \frac{c_1t/l}{V_D(y_i, {r/l})   \phi_2(r/l)}\exp(-c_2(r/l)^{\beta^*})\ge \frac{c_3t/l}{V_D(y_i, {r/l})}\exp(-2c_2(r/l)^{\beta^*})
 \nn\\
 & \ge\frac{(c_3/2)(t/r)}{V_D(y_i, {r/l})}
( \log(r/t))^{1/\beta^*}
 (t/r)^{2c_2}\ge\frac{(c_3/2)(t/r)^{1+2c_2}}{V_D(y_i, {r/l})}
( \log(C_*^{-1}))^{1/\beta^*}
 \ge \frac{c_4(t/r)^{c_5}}{V_D(y_i, {r/l})},
  \end{align}
where $c_5>1$.  Let $B_i=B_{{\bar D}}(x_i,r/(8l))$ for all $1\le i\le l$.
Using
\eqref{chalb} and  \eqref{univd},
we have
\begin{align*}
   q(t, x,y)
 &\geq   \int_{
   B_1}\dots\int_{B_{l-1}}   q(t/l, x,y_1)\dots
     q(t/l, y_{l-1},y)\,m(dy_1)\dots \,m(dy_{l-1})  \\
          & \geq      \frac{c_4(t/r)^{c_5}}{V_D(x, {r/l})}
      \prod_{i=1}^{l-1}c_4(t/r)^{c_5}   \frac{V_D(x_i, {r/l})}{V_D(x_i,9r/(8l))}  \geq    \frac{c_6}{V_D(x, {r/l})}\left(c_7(t/r)^{c_5}\right)^{l}\\
           & \geq      \frac{c_6}{V_D(x, {r})}
             \exp\left(-c_8 r(\log (r/t))^{-1/\beta^*}(1+ \ln (r/t)) \right)\\
& \ge   \frac{c_{6}}{V_D(x,r)}\exp \left(-c_{9}r(\log (r/t))^{(\beta^*-1)/\beta^*}\right).
\end{align*}
This yields (\ref{ew2181}) by considering the cases that $t\in (0,1]$ and $t>1$ respectively.

(ii) Let $c,C^*\in (0,1)$. For any $t>0$ and $r:=\rho_D(x,y)\ge  c\vee (t/C^*)$, we take $c_0=10/c$,
and let $l \ge
2$ be an integer such that $c_0 r\le  l \le c_0r+1\le 2 c_0r$. Choose
$x=x_0,x_1,\cdots,x_l=y\in {\bar D}$ such that $1/(4c_0)\le r/(2l)\le \rho_D(x_i, x_{i+1})\le 2r/l
\le 2/c_0$ for $i= 0, \cdots, l-1$.
It is clear that
$$
\frac{t}{l}\le \frac{C^*r}{l} \le \frac{C^* }{c_0} \quad \text{and}\quad
 \frac{1}{8c_0}\le  \frac{r}{4l} \le \frac{1}{c_0}.
 $$
  Then, since $\rho_D(y_i,y_{i+1}) \le (9/4)r/l \le 9/c_0=9c/10<1$ for $(y_i,y_{i+1}) \in B_{{\bar D}}(x_i,r/(8l)) \times B_{{\bar D}}(x_{i+1},r/(8l))$,
 according to Proposition
\ref{main50}(iii), we
have for all $(y_i,y_{i+1}) \in B_{{\bar D}}(x_i,r/(8l)) \times B_{{\bar D}}(x_{i+1},r/(8l))$
\begin{equation}\label{chalb_1}
    q(t/l, y_i,y_{i+1})\geq
\frac{c_1t/l}{V_D(y_i, r/l)\phi_2 ( r/l)}\ge \frac{c_2t/l}{V_D(y_i,r/l)}.
\end{equation}
 Let $B_i=B_{{\bar D}}(x_i,r/(8l))$ for all  $1\le i\le l$.
 Using
  \eqref{chalb_1} and  \eqref{univd}
 we
have
\begin{align*}
  q(t, x,y)&\geq \int_{
B_1}\dots\int_{B_{l-1}}   q(t/l, x,y_1)\dots
  q(t/l, y_{l-1},y)\,m(dy_1)\dots
\,m(dy_{l-1})\\
&\geq   \frac{c_2(t/l)}{V_D(x,r/l)}(c_3(t/l))^{l-1}
\,\ge \frac{c_4}{V_D(x,r)}
(c_5t/l)^{2c_0r}
\, \ge\,  \frac{c_6}{V_D(x,r)} (t/r)^{c_7r}.
\end{align*}
This completes the proof. \qed

Now, we can state the main result in this section, which provides lower bounds of
$  q(t,x,y)$.  Recall that $H_{\phi, \beta}(t, x, r)$ is defined  right
  after \eqref{eqn:4}.

\begin{theorem}\label{general}
Suppose that  ${\rm diam} (D)=\infty$, and that
$({\bf J}_{\phi_1, 0_+, \leq})$  and $(\bJ_{\phi_2, \beta^*, \geq })$ hold for some $\beta^*\in (0, \infty]$. Then there
are positive constants
$c_1$ and $c_2$, which
depend only on
the characteristic constants $(C_1, C_2, C_3, C_4)$ of $D$
and the constant parameters
$({\bf J}_{\phi_1, 0_+, \leq})$
and $(\bJ_{\phi_2, \beta^*, \geq })$
as well as in \eqref{eqn:poly} for $\phi_1$ and $\phi_2$,
so that
\begin{equation}
  q(t, x, y) \geq c_1
H_{\phi_2, \beta^*}
(t, x, c_2 \rho_D(x, y))
\quad \hbox{for all }
t>0
 \hbox{ and } x, y \in \bar D.
\end{equation}
\end{theorem}

\noindent {\bf Proof}.
(i)
 When $\beta^*\in (0,1]$,
  the assertion for the case that $t\in (0,1]$ follows from Corollary \ref{main5111} and Proposition \ref{main50}(ii).
 When $t>1$, near diagonal lower bounds of $  q(t,x,y)$ is proved by Proposition \ref{P:3.1}. From \eqref{ew218} we see that
 for every $C_*>0$, $t>0$ and $x,y\in \bar D$ with $\rho_D(x,y)\ge C_*t^{1/2}$,
\begin{equation}\label{e:llesi}\begin{split}  q(t,x,y)\ge& \frac{c_1t}{V_D(x,\rho_D(x,y))} e^{-c_2\rho_D(x,y)^{\beta^*}}\ge \frac{c_3t}{V_D(x,\sqrt{t})}\left(\frac{t}{\rho_D(x,y)^2}\right)^{d_2/2} e^{-c_2\rho_D(x,y)^{\beta^*}}\\
\ge & \frac{c_4}{V_D(x,\sqrt{t})} e^{-c_5\rho_D(x,y)^{\beta^*}},\end{split}\end{equation} where in the second inequality we used
 \eqref{univd}.
 Thus, off diagonal lower bounds of $  q(t,x,y)$ for $t\ge 1$ are a consequence of  Proposition \ref{main50}(i) and \eqref{e:llesi}.

(ii) When $\beta^*\in (1,\infty)$,
we can get from \eqref{ew2181} that for any $C^*\in (0,1)$, $t\in (0,1)$ and $x,y\in \bar D$ with $\rho_D(x,y)\ge t/C^*$,
$$  q(t,x,y)\ge \frac{c_1t}{V_D(x,\rho_D(x,y))\phi_2(\rho_D(x,y))}    \exp\left(-c_2\rho_D(x,y) \Big(1+\log^+ \frac{\rho_D(x,y)}{t} \Big)^{({\beta^*}-1)/{\beta^*}}\right),$$  where the constants $c_1,c_2$ may be different from those in \eqref{ew2181}. Moreover, we can also obtain form \eqref{ew2181} and
 \eqref{univd}
 that for any $t\in [1,\infty)$ and $x,y\in \bar D$ with $\rho_D(x,y)\ge t/C^*\ge t^{1/2}/C^*$ with any $C^*\in (0,1)$,
\begin{align*}  q(t,x,y)\ge&\frac{c_3 t}{V_D(x,\sqrt{t})}\left(\frac{t^{1/2}}{\rho_D(x,y)}\right)^{d_2}\, \exp\left( -c_4  \rho_D(x,y) \Big(1+\log^+\frac{\rho_D(x,y)}{t} \Big)^{({\beta^*}-1)/{\beta^*}} \right)\\
\ge& \frac{c_5 t}{V_D(x,\sqrt{t})}\exp\left( -c_6 \rho_D(x,y) \Big(1+\log^+\frac{\rho_D(x,y)}{t} \Big)^{({\beta^*}-1)/{\beta^*}}\right).\end{align*}
With these two estimates above at hand, we can obtain lower bounds of $  q(t,x,y)$ for the case that $\beta^*\in (1,\infty)$ from
Proposition \ref{P:3.1}, Proposition \ref{main50},
Corollary \ref{main5111}, and  Proposition \ref{main33}.

(iii) The $\beta^*=\infty$ case    be verified
by a similar argument as that  for (ii).
\qed

The  lower bounds of $q(t,x,y)$ in Theorem \ref{T1} follows directly from
 \eqref{e:Jcom} and Theorem \ref{general}.

\section{Appendix}\label{S:6}

\subsection{Preliminary estimates}\label{S:4.1}
\begin{lemma}\label{intelem} Assume that
$\varphi$ is a strictly increasing function on $[0, \infty)$
such that $\varphi (0)=0$, $\varphi (1)=1$ and
 \begin{equation}\label{eqn:polyn}
 c_1 \Big(\frac Rr\Big)^{\varrho_1} \,\leq\,
\frac{\varphi (R)}{\varphi  (r)}  \ \leq \ c_2 \Big(\frac
Rr\Big)^{\varrho_2}
\quad \hbox{for every } 0<r<R<\infty,
\end{equation} where
$0<\varrho_1\le \varrho_2<2$
and $c_1,c_2>0$.
Then there exist  constants $c_i \in (0, \infty)$,
$i=3,4,5,$
depending on  $\varrho_1, \varrho_2$,  $c_1$ and $c_2$,
such that the following are true.
\begin{itemize}
\item[ {\rm(i)}]Let $\gamma>0$ and $\beta\in [0,\infty)$.
Then, for all $r>0$,  $$\sup_{\eta\in {\bar D}} \int_{{D} \setminus B_{{\bar D}}(\eta,r)}\frac{e^{-\gamma \rho_D(\xi,\eta)^{\beta}}}{ V_D( \eta,\rho_D (\eta, \xi))\,\varphi (\rho_D (\eta,\xi))}\,m (d\xi)\,\le\,c_3 e^{-\gamma r^\beta} /\varphi (r).$$
\item[{\rm(ii)}] Let $\gamma>0$ and $\beta\in [0,\infty)$.
 Then, for all $r>0$,
 \begin{equation}\label{e:6.2}
 \sup_{\eta\in {D}}\int_{B_{{\bar D}}(\eta,r)}\frac{e^{-\gamma \rho_D(\xi,\eta)^{\beta}   }\r(\eta,\xi)^2}{  V_D(\eta, \rho_D (\eta, \xi))\varphi(\rho_D (\eta,\xi))}\,m (d\xi)\le c_4 \int_0^r\frac{s}{\varphi (s)} e^{-\gamma (s/2)^{\beta}} \, ds.
 \end{equation}
In particular, when $\beta\in (0,\infty)$,
$$\sup_{\eta\in {D}}
\int_{D}
\frac{e^{-\gamma \rho_D(\xi,\eta)^{\beta}   }\r(\eta,\xi)^2}{  V_D(\eta, \rho_D (\eta, \xi))\varphi(\rho_D (\eta,\xi))\,}\,m (d\xi)\le
c_4
\int_0^\infty\frac{s}{\varphi (s)} e^{-\gamma (s/2)^{\beta}} \, ds<\infty.$$

\item[{\rm(iii)}]
Let $\gamma>0$ and $\beta\in [0,\infty)$. Then, for all $r>0$, $$\sup_{\eta\in {D}}\int_{B_{{\bar D}}(\eta,r)}\frac{e^{\gamma \rho_D(\xi,\eta)^{\beta}   }\r(\eta,\xi)^2}{
  V_D( \eta,\rho_D (\eta, \xi))\varphi(\rho_D (\eta,\xi))}\,m (d\xi)\le
 c_5
  \int_0^r \frac{s}{\varphi (s)}  e^{\gamma (2s)^{\beta}}\,ds.$$
     \end{itemize}
\end{lemma}

\pf
(i)
By  \eqref{univd},
we have for any $\eta\in {D}$,
\begin{align*}
&\int_{{D} \setminus B_{{\bar D}}(\eta,r)}\frac{e^{-\gamma \rho_D(\xi,\eta)^{\beta}}}{V_D( \eta,\eta \rho_D (\eta, \xi))\varphi (\rho_D (\eta,\xi))  }\,m (d\xi)\\
=& \sum_{i=0}^\infty \int_{\{ \xi \in D: 2^i r\le \rho_{{\bar D}}(\xi,\eta)<2^{i+1} r\}}\frac{e^{-\gamma \rho_D(\xi,\eta)^{\beta}}}{ V_D(\eta, \rho_D (\eta, \xi))\varphi (\rho_D (\eta,\xi)) }\,m (d\xi)\\
\le &  c_1\sum_{i=0}^\infty \frac{e^{-\gamma (2^ir)^{\beta}}}{\varphi (2^ir)\,  V_D(\eta,2^ir)}\left(
   V_D(\eta,2^{i+1}r)- V_D(\eta,2^ir)\right)\\
\le & c_1e^{-\gamma r^{\beta}}\sum_{i=0}^\infty \frac{1}{\varphi(2^ir)\, V_D(\eta,2^ir)}
(c_2-1) V_D(\eta,2^ir)\\
\le &  \frac{c_3e^{-\gamma r^{\beta}}}{\varphi (r)} \sum_{i=0}^\infty\frac{\varphi (r)} {\varphi (2^ir)} \,\le \, \frac{c_4e^{-\gamma r^{\beta}}}{\varphi (r)} \sum_{i=0}^\infty2^{-i\varrho_1} \,\le \, \frac{c_5e^{-\gamma r^{\beta}}}{\varphi (r)},
\end{align*}
where the lower bound in (\ref{eqn:polyn}) was used in the second to the last inequality.

\bigskip

\noindent
 (ii)
By  \eqref{univd}
and (\ref{eqn:polyn}),
for any $\eta\in {D}$, we have
\begin{align*}
 &\int_{B_{{\bar D}}(\eta,r)}\r(\eta,\xi)^2\frac{e^{-\gamma \rho_D(\xi,\eta)^{\beta}   }}{  V_D( \eta, \rho_D (\eta, \xi))\varphi(\rho_D (\eta,\xi))
}\,m (d\xi)\\
\le&   c_1\sum_{i=0}^\infty (2^{-i}r)^2    e^{-\gamma (2^{-i-1}r)^{\beta}   }\frac{  V_D(\eta, 2^{-i}r)-V_D(\eta, 2^{-i-1}r )}
{V_D(\eta, 2^{-i-1}r )\varphi(2^{-i-1 }r )}\\
\le & c_1\sum_{i=0}^\infty (2^{-i}r)^2    e^{-\gamma (2^{-i-1}r)^{\beta}   }\frac{(c_2-1) V_D(\eta, 2^{-i-1}r)}
{V_D(\eta, 2^{-i-1}r )\varphi ( 2^{-i-1}r )}\\
= &c_3\sum_{i=0}^\infty \frac{( 2^{-i}r)^2}{\varphi (2^{-i}r)}    e^{-\gamma (2^{-i-1}r)^{\beta}   }\,\le\, c_4\int_0^r\frac{s}{\varphi (s)} e^{-\gamma (s/2)^{\beta}} \, ds.
\end{align*}
This establishes \eqref{e:6.2}.
Taking $r\to \infty$ in \eqref{e:6.2} yields  the   second assertion in (ii) .

(iii)  Its proof  is similar to that of (ii) and  is thus omitted.
\qed

Recall that $(\EE^{0, \rf},
\FF)$
 defined in \eqref{EErf} is a strongly local regular
Dirichlet form on $L^2({D}; m)$,
 and that (VD) and (PI(2)) hold
for $(\EE^{0, \rf},
\FF
)$ on $({\bar D}, \rho_D, m)$.
The diffusion process $Z$ associated with
$(\EE^{0, \rf}, \FF^{0, \rf}_D)$ admits a jointly
continuous transition density function $p_D^N(t, x, y)$ on
$(0, \infty)\times {\bar D}\times {\bar D}$, which enjoys the estimate in \eqref{e:1.4}.
Let $\alpha\in (0,2)$. We consider
$\alpha/2$-stable subordinator $S:=(S_t)_{t\ge0}$ independent of $Z$.
Let $(Z_{S_t})_{t\ge0}$ be the subordinated diffusion.
It is easy to see that
 (e.g. see \cite{Ok})
$(Z_{S_t})_{t\ge0}$ has a transition density $ p_\alpha (t,x,y)$ given by
$$
 p_\alpha (t,x,y)=\int_0^\infty  p_D^{N}(s, x,y)\,\P(S_t \in ds),\quad x,y\in \bar D, t>0,$$ and the jumping function of $(Z_{S_t})_{t\ge0}$ is given by
\begin{equation}\label{jumping function1}
j_\alpha(x,y):= \int_0^{\infty}p_D^{N}(s, x,y)\mu(s)\,ds,\,\quad x,y\in D,
\end{equation}
where $\mu(t)=c(\alpha) t^{-1-\alpha/2}$ is the L\'evy density of $S$.
Further, we have the following estimate for $j_\alpha(x,y)$.

\begin{lemma}\label{p:Ialpha1} There is a constant  $c_0\ge 1$, depending on $\alpha\in (0, 2)$ and the constants in \eqref{e:1.4},
such that for every $x,y \in {D}$,
$$\frac{c_0^{-1}}{V_D(x,\rho_D (x, y))\rho_D (x, y)^{\alpha}}\le j_\alpha(x,y)\le \frac{c_0}{V_D(x,\rho_D (x, y))\rho_D (x, y)^{\alpha}} .$$
\end{lemma}
\pf
By \eqref{jumping function1} and \eqref{e:1.4},
$j_\alpha(x,y) \le c_1 k^\alpha (x,c_2 \rho_D (x, y))$, where
\begin{align*}
k^\alpha(x,r):&= \int_0^{\infty}\frac{1}{V_D(x,\sqrt{s})} \exp\left(-\frac{r^2}{s}\right) s^{-1-\alpha/2}\,ds\\
&= \int_0^{r^2}\frac{1}{V_D(x,\sqrt{s})} \exp\left(-\frac{r^2}{s}\right) s^{-1-\alpha/2}\,ds
+ \int_{r^2}^\infty\frac{1}{V_D(x,\sqrt{s})} \exp\left(-\frac{r^2}{s}\right) s^{-1-\alpha/2}\,ds.
\end{align*}

With the change of variable $t=r^2/s$, we have
\begin{align}\label{e:jchvar}
\int_{r^2}^\infty\frac{1}{V_D(x,\sqrt{s})} \exp\left(-\frac{r^2}{s}\right) s^{-1-\alpha/2}\,ds
=r^{-\alpha}\int_{0}^1\frac{1}{V_D(x,r/\sqrt{t})} \exp\left(-t\right) t^{-1+\alpha/2}\,dt.
\end{align}
Since  $V_D(x,r) \le V_D(x,r/\sqrt{t})$ for $x\in {\bar D}$, $r>0$ and $t \in (0, 1]$, the above is less than or equal to
\begin{align*}
r^{-\alpha}\frac 1{V_D(x,r)} \int_{0}^1 \exp\left(-t\right) t^{-1+\alpha/2}\,dt \le \frac{c_3}{r^{\alpha}V_D(x,r)}.
\end{align*}
On the other hand, by
\eqref{univd},
for $s \le r^2$,
$$
\frac{1}{V_D(x,\sqrt{s})} = \frac1{V_D(x,r)} \frac{V_D(x,r)}{V_D(x,\sqrt{s})} \le \frac{c_4}{V_D(x,r)} \frac{r^{d_2}}{s^{d_2/2}}.
$$
Thus,
 \begin{align*}
\int_0^{r^2}\frac{1}{V_D(x,\sqrt{s})} \exp\left(-\frac{r^2}{ s}\right) s^{-1-\alpha/2}\,ds
&\le \frac{c_4}{V_D(x,r)}  \int_0^{r^2}\exp\left(-\frac{r^2}{ s}\right) r^{d_2} s^{-1-\alpha/2-d_2/2}\,ds\\
&\le \frac{c_4}{V_D(x,r)} \left(\sup_{a\ge1}  e^{-a} a^{1+ \alpha/2+d_2/2}\right)
 \int_0^{r^2}  r^{-2-\alpha}\,ds \le \frac{c_5}{r^{\alpha}V_D(x,r)}.
\end{align*}
 Combining both
estimates above and using
 \eqref{univd},
we obtain the upper bound.

Next, we consider the lower bound.
By \eqref{jumping function1} and \eqref{e:1.4},
$j_\alpha(x,y) \ge c_6 k^\alpha (x,c_7\rho_D (x, y))$.
According to
\eqref{univd},
for all $x\in {\bar D}$, $r>0$ and  $t \in (0, 1]$,
$$
\frac{1}{V_D(x,r/\sqrt{t})} = \frac1{V_D(x,r)} \frac{V_D(x,r)}{V_D(x,r/\sqrt{t})} \ge \frac{c_8}{V_D(x,r)} t^{d_2/2}.
$$
Thus, using \eqref{e:jchvar}, we have
 \begin{align*}
k^\alpha(x,r)\ge \frac{c_8}{r^{\alpha}V_D(x,r)} \int_{0}^1 \exp\left(-t\right) t^{-1+(\alpha+d_2)/2}\,dt \ge \frac{c_{9}}{r^{\alpha}V_D(x,r)}.
\end{align*}
Hence, the lemma follows from
\eqref{univd}.
\qed

In the following, let
$(\EE^{(\alpha)}, \FF^{(\alpha)})$ be the Dirichlet form associated with the subordinated diffusion $(Z_{S_t})_{t\ge0}$,
where  $S:=(S_t)_{t\ge0}$ is an $\alpha/2$-stable subordinator  independent of $Z$ with $\alpha\in (0,2)$.
Then,
$$
 \EE^{(\alpha)}
(f,f)=\frac{1}{2}\int_{{D}} \int_{{D}} (f(y)-f(x))^2
 j_{\aa} (x, y) \, m(dx)\, m(dy),\quad f\in
  \FF^{(\alpha)}.
 $$
 We have by \cite[Theorem 2.1(i)]{Ok}

\begin{lemma}\label{p:Ialpha1--}
 $\FF^{0, \rf}_D \subset  \FF^{(\alpha)}$, and, there is a constant $c_1>0$ such that
 $$
\EE^{(\alpha)}_1 (f,f):=\EE^{(\alpha)} (f,f)+\|f\|_2^2 \le c_1     \EE^{0, \rf}_1 (f,f)
\quad \hbox{for every } f\in \FF^{0, \rf}_D .
 $$
 \end{lemma}

\subsection{UJS}\label{S:4.2}

Following  \cite{BBK,CKK2,CKK3},  for $R\in (0, \infty]$, we say the jumping kernel $J(x, y)$
satisfies the $ {\bf UJS}_R $ condition if   there is a constant $c>0$ such that
for $m$-a.e. $x, y\in {  D}$,
$$
 J(x,y) \le \frac{c}{V_D(x,r)}
\int_{B_{{\bar D}} (x,r)} J(z,y)\, m(dz) \quad
\hbox { whenever }  r\le  \rho_D (x, y)/2 <R. \eqno({\bf UJS}_R)
 $$
 When $R=\infty$,  the ${\bf UJS_{\infty}}$ condition will simply be called
 the {\bf UJS} condition.

The {\bf UJS} condition plays an important role in the study of parabolic Harnack inequalities for non-local Dirichlet forms; see \cite{CKK2, CKW2}.
If $J(x, y)$ satisfies both $(\bf J_{\phi, \beta_*, \leq })$ and  $(\bf J_{\phi, \beta^*, \geq })$ with $\beta_*\leq \beta^*$ in $[0, \infty]\cup \{0_+\}$,
then it is easy to check that ${\bf UJS}_{1}$
condition holds. It follows from the proofs of
\cite[Theorem 5.2]{CKK2},
 \cite[Theorem 3.8]{CKW2} and  \cite[Theorem 1.18]{CKW4} that
parabolic Harnack inequalities for finite ranges  hold as well.
See \cite[page 1071]{CKK2}  for the definition of parabolic Harnack inequalities for finite ranges. We note that
the {\bf UJS} condition holds if the jumping kernel $J(x, y)$ has two-sided bounds $(\bf J_{\phi,  \infty, \leq })$ and  $(\bf J_{\phi, \infty , \geq })$.
In this case, the scale-invariant parabolic Harnack inequality holds for full ranges  holds.

Suppose there are positive constants
$c_1$, $c_2$,  $\b \in (0,\infty)$,  $\theta_1\geq \theta_2 \geq 0$ and a strictly increasing function $\phi$ on $[0, \infty)$
satisfying $\phi (0)=0$, $\phi (1)=1$ and \eqref{eqn:poly}
so that
\begin{equation} \label{e:bm1}
   c_1\frac{\exp \left(-\theta_1 \rho_D(x, y)^\beta \right) }{V_D(x, \rho_D(x, y) \phi (\rho_D(x, y)) }  \le J(x, y)
 \leq  c_2 \frac{\exp  \left(-\theta_2 \rho_D(x, y)^\beta   \right)}{V_D(x, \rho_D(x, y) \phi (\rho_D(x, y)) }
  \end{equation}
for every $(x,y) \in D\times D \setminus
 {\rm diag} $.
 The next lemma shows that,
for the jumping kernel $J(x,y)$ satisfying   \eqref{e:bm1}  with $ \theta_1= \theta_2$,
then {\bf UJS} holds.
 However,  when $ \theta_1\not=  \theta_2$ in \eqref{e:bm1},
then {\bf UJS} does not hold in general.  (See \cite[Example 2.4]{CKK3}.)

\begin{lemma}\label{lem2-1ujs}
 Suppose that the jumping kernel
$J(x,y)$ satisfies \eqref{e:bm1} with
$ \theta_1= \theta_2$. Then {\bf UJS} holds.

\end{lemma}
\pf
 The proof is
a simple modification of  that of \cite[Lemma 2.1]{CKK3}.
 Let $x,y\in \bar D$. Suppose $2r\le \rho_D(x,y)$.
Let
$$ A_{x, y,r}:=\left\{ z \in B_{{\bar D}}(x,r) :\rho_D(z,y) \le \rho_D(x,y) \right\} ,
$$ and define $    \tilde \phi(r):= \phi(r) e^{\theta_1 r^\beta} $.
Note that
\begin{align*}
\int_{B_{{\bar D}}(x,r)}  \frac{m(dz)}{V_D(z,\rho_D(z,y))
  \tilde \phi(\rho_D(z,y))}
&\ge
\int_{A_{x, y,r}}  \frac{m(dz)}{V_D(z,\rho_D(z,y))  \tilde \phi(\rho_D(z,y))}\ge
\frac{c_2m(A_{x, y,r})}{V_D(x,\rho_D(x,y)) \tilde \phi(\rho_D(x,y))}.
\end{align*}
Let  $\gamma$ be a continuous curve with $\gamma(0)=x$ and  $\gamma(1)=y$
such that $\hbox{length}(\gamma)\le  \rho_D(x,y)+r/4$.
Choose $w \in \gamma$ such that $\rho_D(x, w)=r/2$.
Then
$$
 \rho_D(x,y)+\frac{r}{4}\ge  \hbox{length}(\gamma)= \hbox{length}(\gamma_{x,w})+\hbox{length}(\gamma_{w,y})\ge  \frac{ r}{2} +\rho_D(w,y),
 $$
 so that
 $\rho_D(w,y) \le \rho_D(x,y) - 4^{-1} r$.
 Thus, for every $z \in B_{{\bar D}}(w,  4^{-1} r)$, we have
 $\rho_D(z,y) \le \rho_D(w,z) +\rho_D(w,y) \le \rho_D(x,y)$, and so
 $B_{{\bar D}}(w,  r/4) \subset A_{x, y,r}$.
 Therefore, by
\eqref{e:1.8},
$$m(A_{x, y,r}) \ge m(B_{{\bar D}}(w, r/4))=V_D(w,r/4) \ge c_3 V_D(x,r).$$
 Using
\eqref{e:bm1} with $ \theta_1= \theta_2$,
we conclude that there exists a constant $c_5>0$ such that for every $r > 0$ and $ 2r\le \rho_D(x,y)$, we have
$$
\int_{B_{{\bar D}} (x,r)} J(z,y)\, m(dz)  \,\ge\,
\frac{c_3\kappa_1^{-1} V_D(x,r)}{V_D(x,\rho_D(x,y))   \tilde \phi(\rho_D(x,y))}\,\ge\,
c_5V_D(x,r)J(x,y).$$
The proof is complete.
\qed

\subsection{Heat kernel estimates: the case
$0=\beta_*<\beta^*\le \infty$ }\label{Case3}

 Note that the lower bound heat kernel estimates on $q(t, x,y)$ in
 Theorem \ref{general} are obtained under the condition $({\bf J}_{\phi_1, 0_+, \leq})$ in addition to
 $(\bJ_{\phi_2, \beta^*, \geq })$ with  $\beta^*\in (0, \infty]$.
 In this subsection, we comment that it seems hard to replace  condition $({\bf J}_{\phi_1, 0_+, \leq})$ by
 a weaker condition
 condition $({\bf J}_{\phi_1, \leq })$.

Suppose that $D$ is unbounded. We first  observe, under the condition $(\bJ_{\phi_1, 0_+, \leq 1})$ which is stronger than $({\bf J}_{\phi_1, \leq })$,  what estimates one can obtain.
 Let   $(\EE,\FF)$ be  the non-local Dirichlet form   given by
\eqref{e:sF}--\eqref{e:sE} with jumping kernel $J(x,y)$ satisfying
$(\bJ_{\phi_1, 0_+, \leq 1})$ and $(\bJ_{\phi_2, \beta^*, \geq})$ for some  $\beta^*\in (0, \infty]$.
By  Proposition \ref{l:regular},  $(\EE,\FF)$ is   a regular Dirichlet form on $L^2(D; m)$.
   By Remark \ref{Add:oneremark}, Theorem \ref{T:2.4} and Proposition \ref{P:3.1},
there is a conservative
 Feller process $Y$ on $\bar D$ associated with the regular Dirichlet form $(\EE,\FF)$ on $L^2(D; m)$
 that starts from every point in $\bar D$. Moreover, $Y$ has a jointly H\"older continuous transition density function
 $q(t, x, y)$ on $(0, \infty) \times \bar D \times \bar D$ with respect to the measure $m$,
  and there are constants $c_1,c_2>0$ so  that
 for all $x,y\in {\bar D} $ with $\rho_D(x,y)\le 1$ and
   $t\in (0, 1]$,
$$
 q(t,x,y)\le c_1
\left(\frac{1}{V_D(x, \sqrt{t}) } \wedge \left( p^{(c)}(t,x, c_2\rho_D(x,y))+ p_{\phi_1}^{(j)}(t,x,  \rho_D(x,y))
  \right)\right),
$$ where $p^{(c)}(t,x,r)$ and $ p^{(j)}_{\phi_1}(t,x,r)$ are defined by \eqref{e:pc} and \eqref{eqn:5},  respectively,
and there are constants
$c_3,c_4>0$ such that for any
 $t\in (0,1]$ and any $x,y\in {\bar D}$ with $\rho_D(x,y)\le c_2t^{1/2}$,
$$ q(t,x,y)\ge \frac{c_3}{V_D(x,\sqrt{t})}.$$
Following
the proof of \cite[Theorem 1.13]{CKW4} and the arguments in Section \ref{S:5}, we can obtain that for all $x,y\in {\bar D}$ with $\rho_D(x,y)\le 1$
 and $t\in (0,1]$,
$$ c_3H_{\phi_2, \beta^*}   ( t, x, c_4 \rho_D(x, y)  ) \le  q(t,x,y) \leq c_5 H_{\phi_1, 0} (t, x, c_5 \rho_D(x, y)).$$

However,
 it seems to be hard to obtain good explicit  estimates for $q(t,x,y)$  especially for large times
when the jumping kernel $J(x,y)$ satisfies the lower bound condition
  $(\bJ_{\phi_2, \beta^*, \geq})$  for some $\beta^*\in (0, \infty]$ and merely the upper bound condition $(\bJ_{\phi_1,   \leq 1})$.
One of the reasons is that the behavior of the on-diagonal estimates for $q(t,x,y)$ is completely different between the case $\beta_*=\beta^*=0$ and $0<\beta_*\le \beta^*\le \infty$ as shown in Theorems \ref{T:1.4} and \ref{T1}.

\bigskip

\noindent {\bf Acknowledgement.}
The research of Zhen-Qing Chen is partially supported by Simons Foundation Grant 520542.
The research of Panki Kim is supported by
 the National Research Foundation of Korea (NRF) grant funded by the Korea government (MSIP)
(No.\ 2016R1E1A1A01941893).\ The research of Takashi Kumagai is supported
by JSPS KAKENHI Grant Number JP17H01093 and by the Alexander von Humboldt Foundation.\
 The research of Jian Wang is supported by the National
Natural Science Foundation of China (Nos.\ 11831014 and 12071076), the Program for Probability and Statistics: Theory and Application (No.\ IRTL1704) and the Program for Innovative Research Team in Science and Technology in Fujian Province University (IRTSTFJ).

\vskip 0.3truein

{\bf Zhen-Qing Chen}

Department of Mathematics, University of Washington, Seattle,
WA 98195, USA

E-mail: zqchen@uw.edu

\bigskip

{\bf Panki Kim}

Department of Mathematical Sciences
and Research Institute of Mathematics,
Seoul

National University, Seoul 08826,
  Republic of Korea

E-mail: pkim@snu.ac.kr

\bigskip

{\bf Takashi Kumagai}

Research Institute for Mathematical Sciences,
Kyoto University, Kyoto 606-8502, Japan

E-mail: kumagai@kurims.kyoto-u.ac.jp

\bigskip

{\bf  Jian Wang}

College of Mathematics and Informatics, \newline
\indent Fujian Key Laboratory of Mathematical
Analysis and Applications (FJKLMAA),\newline
\indent Center for Applied Mathematics of Fujian Province (FJNU),\newline
\indent
Fujian Normal University, 350007 Fuzhou, P.R. China.

 Email: jianwang@fjnu.edu.cn

\end{document}